\documentclass[a4paper,10pt]{article}
\usepackage[french]{babel}
\usepackage[T1]{fontenc}
\usepackage[latin1]{inputenc}
\usepackage{graphics,enumerate}
\usepackage{amsmath,amsfonts,amstext,amssymb,color,epsfig,bbm,lscape}
\usepackage{color}
\usepackage{slashbox}

\newenvironment{enumeratea}
{\bgroup\begin{enumerate}}
{\end{enumerate}\egroup}

\def\C{\mathbb{C}}
\def\R{\mathbb{R}}
\def\Z{\mathbb{Z}}
\def\N{\mathbb{N}}

\def\tr{\text{tr}\,}
\def\H{\mathbb{H}}
\newtheorem{theo}{Th\'eor\`eme}[section]
\newtheorem{coro}[theo]{Corollaire}
\newtheorem{lem}[theo]{Lemme}

\newtheorem{prop}[theo]{Proposition}
\newenvironment{demo}{{\it D\'{e}monstration. }}{\hfill $\square$}

\begin{document}

\begin{center}\begin{LARGE}\textbf{Restriction des séries discrètes de
      $SU(2,1)$ à un sous-groupe exponentiel maximal et à un
      sous-groupe de Borel}\end{LARGE}\end{center}

  \begin{center} Gang LIU  \end{center}

\tableofcontents
\newpage
\null
\newpage
\begin{center} \section{Introduction} \end{center}

La théorie des représentations et l'analyse harmonique constituent un
domaine important des mathématiques, qui est relié à de nombreux
autres domaines des mathématiques et de la physique, comme la théorie
des nombres, les formes automorphes, la géométrie algébrique et la
géométrie différentielle, l'analyse fonctionnelle, la mécanique
quantique, etc.

Soit $G$ un groupe de Lie. Du point de vue de la théorie des
représentations, deux questions ont suscité l'intérêt de nombreux
mathématiciens.

La première concerne la description du dual unitaire $\widehat{G}$ de
$G$, i.e. de l'ensemble des classes d'équivalence de représentations
unitaires irréductibles de $G$. Il s'agit de paramétrer $\hat{G}$ à
l'aide d'un ensemble de paramètres, constitué d'objets géométriques
naturellement liés à $G$, et de donner pour chacun de ces paramètres
une construction d'un représentant de la classe d'équivalence de
représentations associée.

La deuxième, qui relève de l'analyse harmonique, concerne la
décomposition d'une représentation unitaire de $G$ en \flqq somme\frqq~
de représentations unitaires irréductibles.

En fait, on sait que ces deux questions ne peuvent avoir de réponse
satisfaisante que si le groupe $G$ est de type I. C'est le cas des
groupes algébriques réels et donc des groupes que nous rencontrerons
dans ce travail.

Un cas particulier de la deuxième question et auquel nous nous
intéressons ici est connu sous le nom de \flqq branching problem\frqq~
ou \flqq problème de branchement\frqq. Soit $G$ un groupe algébrique
réel et $H\subset G$ un sous-groupe algébrique. Soit
$\pi\in\widehat{G}$. Puisque $H$ est de type $\textrm{I}$, la
restriction de $\pi$ à $H$ que l'on note $\pi\vert_{H}$ peut se
décomposer de la manière (unique) suivante

$$\pi\mid_H= \int_{\widehat H}^\oplus  m_\pi (\tau) \tau \ d\mu_\pi (\tau)$$ où

$\mu_\pi$ est une mesure borélienne sur $\widehat{H}$ (par rapport à
la topologie de Fell)

$m_\pi: \widehat{H} \to \N\cup\{\infty\}$ est une fonction mesurable,
appelée {\it fonction de multiplicité}. Apporter une réponse au
problème de branchement consiste à donner une description explicite de
la mesure $\mu_{\pi}$ et de la fonction de multiplicité $m_{\pi}$.\\

En relation avec le problème de branchement, T. Kobayashi a introduit
dans [15] le concept de {\it $H$-admissibilité} :~ on
dit que $\pi\vert_{H}$ est {\it $H$-admissible}, si le support de la
mesure $\mu_\pi$ est discret et $m_\pi(\tau)<\infty$ pour tout
$\tau\in\widehat{H}$. Autrement dit $\pi\vert_{H}$ est $H$-admissible
si et seulement si $\left.\pi \right\vert_{H}=\widehat{\oplus}_{i\in
  I} k_i.\tau _i $ o\`{u} $I$ est un ensemble dénombrable, $k_i \in \N
$ et les $\tau _i\in \widehat {H} $ sont deux \`{a} deux
diff\'{e}rentes.

Il est intéressant d'étudier les cas où une représentation unitaire
irréductible d'un groupe est admissible pour un de ses
sous-groupes et de résoudre le problème de branchement correspondant.\\

Revenons un moment à la question de la construction de représentations
unitaires irréductibles d'un groupe de Lie.

La {\it méthode des orbites} permet de construire \flqq beaucoup\frqq~
de telles représentations pour un groupe de Lie $G$ de type
$\textrm{I}$. Initiée par Kirillov dans sa thèse, elle a été
développée par de nombreux mathématiciens, parmi lesquels
Auslander-Kostant, pour les groupes résolubles, Duflo pour
les groupes réductifs et Duflo pour les groupes de Lie généraux.

Soit $G$ un groupe de Lie de type I d'algèbre de Lie
$\mathfrak{g}$. Duflo établit une correspondance entre l'ensemble des
$G$-orbites dans $\mathfrak{g}^{*}$, l'espace vectoriel dual de
$\mathfrak{g}$, qui sont {\it admissibles et bien polarisables} et
$\widehat{G}$. Cette correspondance associe à chaque orbite admissible
et bien polarisable un ensemble d'éléments de $\widehat{G}$ paramétré
par un ensemble de {\it données d'admissibilité} attaché à l'orbite et
elle est injective, en ce sens que les ensembles associés à deux
orbites distinctes sont disjoints. Lorsque le groupe $G$ est
nilpotent, on retrouve la correspondance bijective établie par
Kirillov dans sa thèse entre l'ensemble des orbites coadjointes de $G$
et son dual unitaire. Lorsque $G$ est résoluble simplement connexe, on
retrouve les résultats d'Auslander-Kostant (à une translation près
dans la paramétrisation)~: dans ce cas, la méthode permet encore de
décrire complètement le dual unitaire. Par contre, dans le cas général
on n'obtient pas tout le dual unitaire de $G$, mais une partie
suffisamment \flqq grosse\frqq~pour supporter la mesure de
Plancherel. Nous donnerons plus de détails, pour ce qui concerne les
groupes algébriques, dans le chapitre 2.\\

La question qui se pose alors est de savoir si l'on peut apporter une
réponse au problème de branchement dans le cadre de la méthode des
orbites. C'est le cas pour les groupes résolubles exponentiels comme
le montrent les travaux de Lipsman ([22]) et Fujiwara ([10]). Cependant, dans cette situation et pour les séries
discrètes les cas de $H$-admissibilité sont rares (voir Kouki
[18]).

Lorsque $G$ est un groupe de Lie compact, toute représentation
unitaire irréductible de $G$ est $H$ admissible. L'étude du problème
de branchement dans cette situation a conduit à un ensemble de
travaux autour de la conjecture \flqq La quantification commute avec
la réduction\frqq, due à Guillemin et Sternberg, que nous présentons
brièvement ci-dessous.

Soit $(M,\omega, K, \Phi)$ un $K$-espace hamiltonien, où $(M,\omega)$
est une variété symplectique compacte, $K$ est un groupe de Lie
compact d'algèbre de Lie $\mathfrak{k}$, et $\Phi$ désigne l'application
moment. On se donne un fibré de Kostant-Souriau $\mathcal{L}$ sur $M$
qui est un fibré en droites hermitien $K$-équivariant de forme
courbure $-i\omega$. On peut alors construire une quantification notée
$Q^{\mathcal{L}}(M)$ qui est l'indice équivariant d'un opérateur
différentiel elliptique $K$-équivariant associé aux objets considérés
et ainsi une représentation formelle de $K$, cela signifie que
$Q^{\mathcal{L}}(M)$ s'écrit $\Sigma_{\tau\in
  \widehat{K}}n_{\tau}\tau$, avec les $n_{\tau}\in \Z$ et tous nuls
sauf un nombre fini. De plus si $K$ est trivial ou agit trivialement
(sur $M$ et $\mathcal{L}$), alors $Q^{\mathcal{L}}(M)$ s'interprète
comme un nombre.  Fixons un tore maximal $\mathbb{T}$ (d'algèbre de
Lie $\mathfrak{t}$) de $K$ et un ensemble de racines positives de
$\mathfrak{t}$ dans $\mathfrak{k}$ 
:~ on peut donc identifier $\widehat{K}$ à l'ensemble des plus hauts
poids dans $\mathfrak{t}^*\subset \mathfrak{k}^*$ (on identifie
$\mathfrak{t}^{*}$ à l'orthogonal du sous-espace
$[\mathfrak{t},\mathfrak{k}]$ de $\mathfrak{k}$). Soit $\mu$ un plus
haut poids, et $\mathcal{O}_{\mu}=K.\mu$ l'orbite coadjointe de $\mu$
sous $K$ dans $\mathfrak{k}^*$, alors la variété réduite associée
$\Phi^{-1}(\mathcal{O}_{\mu})/K$ est naturellement munie d'un fibré de
Kostant-Souriau induit par celui considéré . Soit
$\tau_{\mu}\in\widehat{K}$ la représentation de plus haut poids
$\mu$. La conjecture de Guillemin-Sternberg dit que le calcul de la
multiplicité $n_{\tau_{\mu}}$, se ramène à quantifier
$\Phi^{-1}(\mathcal{O}_{\mu})/K$. Plus précisément le nombre
$n_{\tau_{\mu}}$ est exactement le résultat de la quantification de la
variété réduite $\Phi^{-1}(\mathcal{O}_{\mu})/K$ et du fibré associé,
tout du moins si la variété réduite est une variété lisse ou un \flqq
orbifold\frqq. Leur conjecture qui est connue sous le slogan : "La
quantification commute avec la réduction", est résolue par Meinrenken
et Sjamaar.

 Il existe une autre formulation de cette conjecture dans le cadre
voisin de la quantification $\text{Spin}_{c}$. Disons simplement que
dans ce cadre, le fibré de Kostant-Souriau est remplacé par
$\tilde{L}$, son \flqq translaté\frqq~ par la racine carrée du fibré
en droite défini par le déterminant d'une structure spinorielle
complexe sur la variété $M$. Le résultat de la quantification est
alors noté $\mathcal{Q}_{\text{spin}}^{\tilde{L}}(M)$. On reprend les
notations précédentes concernant le groupe compact $K$. Soit
$\Delta^{+}$ un ensemble de racines positives de $\mathfrak{t}_{\C}$
dans $\mathfrak{k}_{\C}$, $\delta$ la demi-somme des éléments de
$\Delta^{+}$ et $\mu\in\mathfrak{t}^{*}$ un poids dominant. Alors la
structure spinorielle et le fibré $\tilde{L}$ sur $M$ induisent des
données du même type sur la variété réduite
$\Phi^{-1}(\mathcal{O}_{\mu+\delta})/K$. Ici encore, la multiplicité
de la représentation $\tau_{\mu}$ dans
$\mathcal{Q}_{\text{spin}}^{\tilde{L}}(M)$ est donnée par la
quantification $\text{Spin}_{c}$ de la variété réduite et du fibré
associé.

Cette conjecture a influencé des recherches de relations
entre représentations d'un groupe compact $K$ opérant de manière
hamiltonienne sur une variété symplectique $M$ pas forcément compacte
et l'application moment. Tian-Zhang et Paradan ont obtenu des
démonstrations directes pour le cas général qui s'adapte à certaines
situations nouvelles: variétés à bord, variétés non compactes (dans la
situation où l'application moment est propre) etc. Notamment Paradan
([24, 25, 26]) a traité le cas où $M$ est une orbite coadjointe
associée à une série discrète d'un groupe de Lie réductif $G$ (pas
forcément compact), et $K\subset G$ est un sous-groupe compact
maximal. Dans ce cas $M$ n'est pas forcément compacte, cependant
Paradan a montré que dans cette situation, la quantification commute
également avec la réduction. Il a aussi réussi à réinterpréter la
formule de Blattner (une identité combinatoire) dans le cadre de la
géométrie hamiltonienne.\\

Inspiré par les divers travaux cités dessus, Duflo a formulé une
conjecture concernant le problème de branchement pour les séries
discrètes d'un groupe algébrique réel. Pour des raisons de commodité,
nous allons l'énoncer dans le cas d'un groupe algébrique réductif
connexe.

Soit donc $G$ un groupe de Lie réel connexe {\it algébrique réductif}
d'algèbre de Lie $\mathfrak{g}$, et $H\subset G$ un sous-groupe
algébrique d'algèbre de Lie $\mathfrak{h}$. Soit $\mathfrak{g}^*$
(resp. $\mathfrak{h}^*$) le dual linéaire de $\mathfrak{g}$
(resp. $\mathfrak{h}$). Supposons que $\pi$ est une \textbf{série
  discrète} de $G$. Il est donc bien connu que la s\'{e}rie
discr\`{e}te $\pi$ correspond \`{a} une (unique) $G$-orbite coadjointe
"\textbf{admissible et fortement régulière}" dans $\mathfrak{g}^*$ (au
sens de Duflo). Notons donc $\mathcal{O}_{\pi}$ l'orbite liée à
$\pi$. Munie de la structure symplectique de Kirillov-Kostant-Souriau
$\omega$, $(\mathcal{O}_{\pi},\omega)$ est une variété symplectique et
comme le stabilisateur d'un point de l'orbite $\mathcal{O}_{\pi}$ est
un tore, elle admet un unique fibré de Kostant-Souriau . Le
sous-groupe $H$ agit dans $\mathcal{O}_{\pi}$ par l'action coadjointe,
de sorte que $(\mathcal{O}_{\pi},\omega,H, \text{p})$ devient un
$H$-espace hamiltonien, l'application moment associée étant la
projection naturelle $\text{p}$ de $\mathcal{O}_{\pi}$ dans
$\mathfrak{h}^*$.

Maintenant, on \'{e}nonce la conjecture de Duflo (on garde les
notations et hypothèses du paragraphe précédent):

\textbf{La conjecture de Duflo} : Notons $\mathfrak{h}_{fr}^*$ la
réunion des $H$-orbites coadjointes fortement r\'{e}guli\`{e}res dans
$\mathfrak{h}^*$. Alors :

(i) Pour que la restriction $\pi$ \`{a} $H$,
$\pi\vert_{H}$, soit $H$-admissible, il faut et il suffit que la
projection (c'est-à-dire l'application moment de
$(\mathcal{O}_{\pi},\omega,H, \text{p})$
) $$\text{p}:\mathcal{O}_{\pi}\longrightarrow
\text{p}(\mathcal{O}_{\pi})\subseteq \mathfrak{h}^* $$ $$ f\longmapsto
\left. f \right \vert_{\mathfrak{h}}$$ soit\textbf{ faiblement
  propre}.

(ii) Si $\pi\vert_{H}$ est $H$-admissible, toute $\tau _i $ qui figure
dans la d\'{e}composition $\left.\pi \right\vert_{H}=\sum k_i.\tau _i$
correspond \`{a} une (unique) $H$-orbite coadjointe admissible et
fortement régulière $\Omega_{\tau _i}$ de $\mathfrak{h}^*$(au sens de
Duflo), avec $\Omega_{\tau _i}$ contenue dans
$\text{p}(\mathcal{O}_{\pi})$.

(iii) Toujours dans le cas où $\pi$ est $H$-admissible, les
multiplicités $k_{i}$ doivent pouvoir s'exprimer à l'aide de la
quantification de la variété réduite correspondante. Autrement dit le
slogan \textbf{ quantification commute avec réduction} reste valable
dans ce cadre.\\

Ici la notion \og faiblement propre\fg~ signifie que l'image
r\'{e}ciproque de tout compact contenu dans
$\text{p}(\mathcal{O}_{\pi})\bigcap \mathfrak{h}_{fr}^*$ par la
projection p est compact. Dans cet article, on utilisera aussi la
notion classique \og \textbf{propre sur l'image} \fg~ selon laquelle
l'image r\'{e}ciproque par la projection p de tout compact contenu
dans $\text{p}(\mathcal{O}_{\pi})$ est un compact. Pour simplifier, dans
cet article, nous dirons "\textbf{propre}" au lieu de \og
\textbf{propre sur l'image} \fg. Il est \'{e}vident que si la
projection est propre, elle est aussi faiblement propre, alors que
l'inverse est g\'{e}n\'{e}ralement faux. \\

Dans ce travail, on s'int\'{e}resse  au cas
o\`{u} $G=SU(2,1)$ et $H$ est un sous-groupe de Borel ou un
sous-groupe exponentiel maximal de $G$.

Maintenant, on énonce les principaux résultats de cet article.\\

Donc soit $G=SU(2,1)$ d'algèbre de Lie $\mathfrak{g}$, $B=MAN$ un
sous-groupe de Borel (de $G$) d'algèbre de Lie $\mathfrak{b}$ et
$B_1=AN$ le sous-groupe exponentiel maximal associé d'algèbre de Lie
$\mathfrak{b}_1$. Soit $\pi$ une série discrète de $G$ et
$\mathcal{O}_{\pi}$ la $G$-orbite coadjointe qui lui est associée par
Duflo. Soit $\text{p}: \mathcal{O}_{\pi}\longrightarrow
\text{p}(\mathcal{O}_{\pi})\subset \mathfrak{b}^* $
(resp. $\text{p}_1: \mathcal{O}_{\pi}\longrightarrow
\text{p}_1(\mathcal{O}_{\pi})\subset \mathfrak{b}_1^* $), la
projection naturelle de $\mathcal{O}_{\pi}\subset \mathfrak{g}^*$ dans
$\mathfrak{b}^*$ (resp. $\mathfrak{b}_1^*$). Donc comme expliqué plus
haut, muni de la structure symplectique de Kirillov-Kostant-Souriau
$\omega$, $(\mathcal{O}_{\pi},\omega)$ est une variété
symplectique. $B$ (resp. $B_1$) agit dans $\mathcal{O}_{\pi}$ par
l'action coadjointe, de sorte que $(\mathcal{O}_{\pi},\omega,B,
\text{p})$ (resp. $(\mathcal{O}_{\pi},\omega,B_1, \text{p}_1)$)
devient un espace hamiltonien, où l'application moment associée est
$\text{p}$ (resp. $\text{p}_1$). De plus, $(\mathcal{O}_{\pi},\omega)$
admet un unique fibré de Kirillov-Kostant-Souriau (car on l'a déjà
expliqué, le stabilisateur d'un point de l'orbite $\mathcal{O}_{\pi}$
est un tore).

On obtient les résultats suivants:

1) (Théorème 5.4) La projection (l'application moment pour
$(\mathcal{O}_{\pi},\omega,B_1, \text{p}_1)$) $\text{p}_1:
\mathcal{O}_{\pi}\longrightarrow \text{p}_1(\mathcal{O}_{\pi})\subset
\mathfrak{b}_1^* $ est propre (sur l'image) ou faiblement propre, si
et seulement si la série discrète $\pi$ est holomorphe ou anti-holomorphe.

2) (Théorème 5.5) La projection (l'application moment pour
$(\mathcal{O}_{\pi},\omega,B, \text{p})$)
$\text{p}:\mathcal{O}_{\pi}\longrightarrow
\text{p}(\mathcal{O}_{\pi})\subset \mathfrak{b}^* $ est faiblement
propre. De plus elle est propre (sur l'image) si et seulement si $\pi$
est holomorphe ou anti-holomorphe.

3) (Théorème 6.3) Si $\pi$ est holomorphe ou anti-holomorphe,
$\pi\vert_{B_1}$ (resp. $\pi\vert_{B}$) est $B_1$-admissible
(resp. $B$-admissible), et on a une décomposition explicite de
$\pi\vert_{B_1}$ (resp. $\pi\vert_{B}$). Le théorème 6.3 énonce, dans
le cas particulier du groupe $SU(2,1)$, un résultat dû à Rossi-Vergne
([29]) dans le cas général d'un groupe réductif
connexe admettant des séries discrètes holomorphes.

4) Ainsi, lorsque $\pi$ est holomorphe (ou anti-holomorphe), d'après
le théorème 5.3 (resp. la proposition 6.5) qui décrit les
$B_1$-orbites (resp. $B$-orbite) coadjointes fortement régulières et
admissibles (au sens de Duflo) contenues dans
$\text{p}_1(\mathcal{O}_{\pi})$ (resp. $\text{p}(\mathcal{O}_{\pi}$)),
les assertions (i) et (ii) de la conjecture de Duflo sont vérifiées
pour $(G,\pi, B_1)$ et $(G,\pi, B)$. On constate alors que, parmi les
séries dicrètes de $B$ associées à une orbite contenue dans
$\text{p}(\mathcal{O}_{\pi})$, seulement une sur trois intervient
dans la décomposition $\pi\vert_{B}$.

5) (Théorème 6.8) Soit $\pi$ une série discrète ni holomorphe ni
anti-holomorphe. En utilisant la réalisation de $\pi$ dans un espace
de $L^2$-cohomologie due à Nara\-simhan-Okamoto ([23])
   ainsi que les résultats de Hersant ([11]), on obtient
un encadrement pour les multiplicités des représentations unitaires
irréductibles de $B$ dans $\pi\vert_{B}$: la multiplicité $n_{\tau}$
de chaque $\tau\in \widehat{B}$ (le dual unitaire de $B$) est égale à
la dimension du
sous-espace des solutions de carré intégrable sur $]0,+\infty[$ (par
    rapport à la mesure de Haar $ dt/t$) d'un système différentiel
    linéaire ordinaire ayant deux points singuliers (en $0$ et
    $+\infty$), dont un de deuxième espèce. En étudiant le
    comportement asymptotique de ses solutions, on obtient une borne
    inférieure et une borne supérieure pour toutes les
    multiplicités. Ceci nous permet de montrer que $\pi\vert_{B}$ est
    $B$-admissible et $\pi\vert_{B_1}$ n'est pas
    $B_1$-admissible. Donc d'après 1) (le théorème 5.4), les
    assertions (i) et (ii) de la conjecture de Duflo sont vérifiées
    pour $(G,\pi, B_1)$ avec $\pi$ ni holomorphe ni
    anti-holomorphe. Donc d'après 4), elles sont vérifiées pour
    $(G,\pi, B_1)$ et pour toute série discrète $\pi$. De plus, cette
    méthode permet de retrouver les résultats de 3) obtenus dans le
    cas des séries discrètes holomorphes.

Il est à noter que le fait que les séries discrètes ni holomorphes ni
anti-holomorphes ne sont pas $B_{1}$-admissibles a été démontré pour
$SU(n,1)$ par Rosenberg-Vergne ([28]), également en
utilisant leur réalisation dans un espace de $L^{2}$-coho\-mologie.

6) (Théorème 6.10) En combinant les travaux de Fabec ([9]) et ceux de
Kraljevic ([19, 20]) sur les séries principales de $SU(2,1)$ avec le
théorème 6.8 et le théorème 6.3 (dû à Rossi et Vergne), on parvient à
décomposer explicitement $\pi\vert_{B}$ pour $\pi$ ni holomorphe ni
anti-holomorphe.

7) En comparant le théorème 6.10 avec la proposition 6.11 (qui
décrit les $B$-orbites coadjointes fortement régulières et admissibles
 contenues dans $\text{p}(\mathcal{O}_{\pi})$), on
montre que les assertions (i) et (ii) de la conjecture de Duflo sont
vérifiées pour $(G,\pi, B)$.  Ainsi, d'après 4) et 5), les assertions
(i) et (ii) de la conjecture de Duflo sont vérifiées pour $(G,\pi, B)$
et $(G,\pi, B_1)$ pour toute série discrète $\pi$. On constate
également que, parmi les séries dicrètes associées à une orbite
contenue dans $\text{p}(\mathcal{O}_{\pi})$, seulement une sur trois
intervient dans la décomposition $\pi\vert_{B}$.

8) (Proposition 6.12) Puisque l'on obtient une décomposition explicite
de $\pi\vert_{B}$ pour $\pi$ ni holomorphe ni anti-holomorphe, en
comparant le théorème 6.8 et le théorème 6.10, on peut trouver la
dimension exacte du sous-espace des solutions de carré intégrable de
chaque système différentiel concerné que l'on évoque dans 4). Ce
résultat semble difficile à obtenir par des méthodes
\flqq directes\frqq.

9) (Chapitre 7) En modifiant certains paramètres, on montre que dans
notre cas les assertions (i) et (ii) de la conjecture de Duflo sont
encore vérifiées si on utilise la paramétrisation de Blattner pour les
séries discrètes et celle d'Auslander-Kostant pour les représentations
unitaires irréductibles de $B_{1}$ et $B$.

10) (Propositions 8.1, 8.2 et 8.3) On montre que les variétés réduites
pour l'espace hamiltonien $(\mathcal{O}_{\pi},\omega,B, \text{p})$
sont des points. On vérifie qu'une représentation de $B$ associée à
une orbite contenue dans $\text{p}(\mathcal{O}_{\pi})$ apparaît dans
$\pi\vert_{B}$ si et seulement son caractère central coïncide avec
celui de $\pi$ ($B$ et $G$ ont même centre).

On interprète ce résultat dans le cadre de la quantification
géométrique et on montre qu'il est l'analogue de celui démontré par
Paradan dans ([26], section 2.4 "Quantization of points", Theorem
2.16), dans lequel il considère un groupe de Lie compact muni d'une
action hamiltonienne sur une variété symplectique avec application
moment propre et calcule la quantification de la variété réduite
lorsqu'elle est réduite à un point.

Ceci montre que l'assertion (iii) de la
conjecture de Duflo est vérifiée et explique pourquoi seulement une
série discrète de $B$ sur trois parmi celles associées à une orbite
contenue dans $\text{p}(\mathcal{O}_{\pi})$ intervient dans la
décomposition $\pi\vert_{B}$.

11) (Propositions 8.1 et 8.4) Lorsque $\pi\vert_{B_1}$ est $B_1$-admissible,
on montre que la variété réduite est une sous-variété symplectique de
$\mathcal{O}_{\pi}$, difféomorphe à une sphère de dimension $2$, et
que la multiplicité s'exprime également comme la quantification de
cette dernière, c'est à dire comme l'intégrale de sa forme volume
naturelle. Ceci montre que l'assertion (iii) de la conjecture de Duflo
est également vérifiée dans ce cas.\\

\noindent \textbf{Remarques et Remerciement}. La présentation de ce texte est essentiellement la version originale de ma thèse ayant été soutenue le 5 juillet 2011 à l'Université de Poitiers. J'aimerais souligner que récemment, j'ai obtenu des résultats (bien) plus généraux concernant la conjecture de Duflo et résolu une autre conjecture formulée dans ma thèse (celle suivant le théorème 5.4). Cependant, je voudrais présenter ici ma thèse sans modification.  Malgré tout, une méthode générale et efficace pour traiter la restriction des séries discrètes d'un groupe de Lie simple à un sous-groupe parabolique semble très difficile à trouver.

Je tiens à remercier P. Torasso de m'avoir proposé ce
sujet et m'avoir initié dans le domaine de la théorie des
représentations et l'analyse harmonique non commutative. Ma profonde
reconnaissance va aussi à M. Duflo, dont les travaux ont joué un rôle
très important dans ce travail, et avec qui j'ai eu des discussions
importantes. Je voudrais aussi remercier autres membres du jury de ma thèse, D. Vogan et M. Vergne d'avoir 
fait certaines remarques éclairantes. Je remercie C. Sabbah pour son
aide déterminante dans l'étude des équations différentielles
ordinaires. Enfin, je remercie T. Kobayashi pour des discussions bien
utiles. Durant la préparation de ma thèse, j'ai bénéficié d'une
allocation de recherche du ministère français de l'enseignement
supérieur et de la recherche.

Finalement, un grand merci à B. Krötz pour me soutenir scientifiquement et moralement après ma thèse.

\pagebreak

\section{La méthode des orbites et la théorie de Duflo}\label{méthode-orbites}

Dans ce chapitre, nous allons rappeler quelques élements dont nous
aurons besoin de la construction de Duflo de représentations unitaires
irréductibles. \\

 \textbf{Orbites fortement r\'{e}guli\`{e}res} : \\

Rappelons qu'un groupe de Lie presque algébrique est la donnée d'un
triplet $(G,j,\mathbf{G})$ où $\mathbf{G}$ est un groupe algébrique
affine défini sur $\mathbb{R}$, $G$ est un groupe de Lie et $j$ est un
morphisme de groupes de Lie de $G$ dans le groupe des points réels
$\mathbf{G}_{\mathbb{R}}$ de $\mathbf{G}$ dont le noyau est discret
central dans $G$ et dont l'image est un sous-groupe ouvert de
$\mathbf{G}_{\mathbb{R}}$ dense dans $\mathbf{G}$ pour la topologie de
Zariski (voir [7], p. 199). On notera plus simplement $G$ le groupe
presque algébrique $(G,j,\mathbf{G})$.

Soit $G$ un groupe de Lie presque algébrique sur $\mathbb{R}$
d'alg\`{e}bre de Lie $\mathfrak{g}$. On dit qu'une forme linéaire
$f\in\mathfrak{g}^{*}$ est régulière si son orbite coadjointe $G.f$
est de dimension maximale ou si, de manière équivalente, son
stabilisateur $G(f)$ dans $G$ est de dimension minimale. Supposons que
$f\in \mathfrak{g}^*$ est une forme lin\'{e}aire régulière.  Alors, il
est bien connu que l'alg\`{e}bre de Lie $\mathfrak{g}(f)$ de $G(f)$
est algébrique ab\'{e}lienne. On d\'{e}signe par $\mathfrak{s}(f)$ son
unique facteur réductif qui est l'ensemble des $X\in \mathfrak{g}(f)$
pour lesquels $\text{ad}X$ est semi-simple. On dit que $f$ est
fortement r\'{e}guli\`{e}re, si $\mathfrak{s}(f)$ est de dimension
maximale lorsque $f$ parcourt l'ensemble des formes régulières. Il est
clair que $f$ est r\'{e}guli\`{e}re (fortement r\'{e}guli\`{e}re), si
et seulement s'il en est de même pour tous les éléments de son orbite
coadjointe. Nous dirons qu'un orbite est régulière (resp. fortement
régulière) si c'est l'orbite du forme régulière (resp. fortement
régulière). \\

\textbf{Formes lin\'{e}aires admissibles et bien polarisables} \\

Soient $G$ un groupe de Lie d'alg\`{e}bre de Lie $\mathfrak{g}$ et
$f\in \mathfrak{g}^*$. Consid\'{e}rons la forme bilin\'{e}aire
altern\'{e}e $B_f$ sur $\mathfrak{g}$ d\'{e}finie
par $$B_f(X,Y)=f([X,Y]),\ \ \ X,\ Y\in \mathfrak{g}.$$

On sait que $\mathfrak{g}(f)$ est le noyau de $B_f$, de sorte que
$\mathfrak{g}/\mathfrak{g}(f)$ est un espace symplectique, et $G(f)$
laisse invariant $B_f$. En 1972, Duflo a introduit le groupe
$G(f)^{\mathfrak{g}}$, un rev\^{e}tement \`{a} deux feuillets de
$G(f)$ associé à $(\mathfrak{g},B_{f})$, qui joue un r\^{o}le central
dans la m\'{e}thode des orbites (pour plus de d\'{e}tails, voir [7],
p.153 et suivantes).

Comme nous en aurons besoin plus loin, nous présentons la construction
de ce revêtement dans un cadre plus général : on se donne un espace
vectoriel $V$ de dimension finie sur $\mathbb{R}$ muni d'une forme
bilinéaire alternée $\beta$ dans lequel agit un groupe de Lie $H$ en
laissant $\beta$ invariante. On considère l'espace symplectique
$V/\ker\beta$, le groupe symplectique $Sp(V/\ker\beta)$ et le groupe
métaplectique $Mp(V/\ker\beta)$, revêtement connexe à deux feuillets
du précédent. Le revêtement à deux feuillets de $H$, dit revêtement
métaplectique associé à $(V,\beta)$, est alors le produit fibré de $H$
par $Mp(V/\ker\beta)$ au dessus de $Sp(V/\ker\beta)$ : on le note
$H^{V}$.

On prolonge $\beta$ par bilinéarité complexe  en une forme bilinéaire
alternée sur l'espace vectoriel complexifié $V_{\mathbb{C}}$ de $V$,
encore notée $\beta$.  On appelle lagrangien complexe de
$V_{\mathbb{C}}$ de $V$ un sous-espace isotrope maximal de
$V_{\mathbb{C}}$ pour $\beta$.

Si $l\subseteq V_{\mathbb{C}}$ est un lagrangien complexe stable sous
l'action de l'algèbre de Lie $\mathfrak{h}$ de $H$, on pose
$$
\rho_{l}(X)=\frac{1}{2}tr X_{l/\ker\beta}\mbox{, }X\in\mathfrak{h}.
$$
\\

\noindent \textbf{Remarque}. $G(f)^{\mathfrak{g}}$ n'est pas
forc\'{e}ment connexe, m\^{e}me si $G(f)$ est connexe. Cependant, si $G(f)$
est connexe, alors, $G(f)^{\mathfrak{g}}$ a au plus 2 composantes
connexes. \\

On appelle $\varepsilon$ l'élément non trivial du noyau de la
projection de $G(f)^{\mathfrak{g}}$ sur $G(f)$. Notons $X(f)$
l'ensemble des classes de repr\'{e}sentations unitaires irréductibles
$\tau$ de $G(f)^{\mathfrak{g}}$ qui v\'{e}rifient les
propri\'{e}t\'{e}s suivantes:
$$\tau(\varepsilon)=-id\ \ \ \ \ \ \ \ \ \ \ \ \ \ \ \ \ \ \ \ \ \ \ \ \ \ \ \ \ \ \ \ \ \ \ \ \ \ \ \ \  (1)$$

$$\text {la diff\'{e}rentielle} \ d\tau \ \text {est un multiple de}
\ if\vert_{g(f)} \ \ \ (2)$$ On dit que $f$ est {\it admissible}, si
$X(f)$ est non vide. Soit $\eta$ un caractère d'un sous-groupe fermé
central $\Gamma$ de $G$, on dit qu'un élément $\tau\in X(f)$ est {\it
  $\eta$-admissible}, si $\tau\vert_{\Gamma}$ est un multiple de
$\eta$.\\

Concernant l'admissibilit\'{e} de $f$, il y a un lemme tr\`{e}s utile
d\^{u} \`{a} Duflo (voir [7], p. 154 Remarque 2) :

\begin{lem}\label{admissibilité} Supposons qu'il existe un sous-espace lagrangien
  complexe $l\subseteq \mathfrak{g}_{\C} $ pour $B_{f}$, qui est
  stable par $\mathfrak{g}(f)$, alors, $f$ est admissible si et
  seulement s'il existe un caract\`{e}re de ${G(f)}_0$ (la composante
  neutre de $G(f)$) dont la diff\'{e}rentielle est $\rho_l +
  if\vert_{\mathfrak{g}(f)}$. \end{lem}

Dans [7], Duflo a introduit la notion de forme linéaire bien
polarisable. Soit $f$ une forme linéaire sur $\mathfrak{g}$. On
appelle {\it polarisation} en $f$ une sous-algèbre de Lie $\mathfrak{b}$ de
$\mathfrak{g}_{\C}$ qui soit un sous-espace lagrangien complexe pour
$B_{f}$.

Soit $\mathfrak{b}$ une polarisation en $f$ et soit $B$ le sous-groupe
analytique d'algèbre de Lie $\mathfrak{b}$ du groupe de Lie simplement
connexe d'algèbre de Lie $\mathfrak{g}_{\C}$. On dit que la
polarisation $\mathfrak{b}$ vérifie la {\it condition de Pukanszky} si
on a $B.f=f+\mathfrak{b}^{\perp}$ (ici $\mathfrak{b}^{\perp}$ désigne
l'orthogonal de $\mathfrak{b}$ dans $\mathfrak{g}_{\C}^{*}$).

La forme linéaire $f$ est dite {\it bien polarisable} si elle admet une
polarisation résoluble et vérifiant la condition de Pukanszky.

 Notons $\mathfrak{g}_{ap}^*\subset \mathfrak{g}^* $, l'ensemble des
 formes linéaires admissibles et bien polarisables. Alors on a
 $G.\mathfrak{g}_{ap}^*=\mathfrak{g}_{ap}^*$. Ceci nous permet de
 définir les {\it orbites coadjointes admissibles et bien
   polarisables} d'une façon évidente. Notons $\widehat{G}$, le dual
 unitaire de $G$. Soit $f\in \mathfrak{g}_{ap}^*$ et $\tau\in X(f)$.
 Au couple $(f,\tau)$, Duflo associe une représentation unitaire
 irréductible $\pi_{f,\tau}$ de $G$.  De plus, si $f_1\in G.f$, il
 existe un et un seul élément $\tau_1$ dans $X(f_1)$ tel que
 $\pi_{f,\tau}\cong \pi_{f_1,\tau_1}$. Ainsi, étant donnée une orbite
 coadjointe admissible et bien polarisable $\mathcal{O}$, il existe un
 ensemble $X(\mathcal{O})$ qui se met canoniquement en bijection avec
 $X(f)$, pour tout $f\in \mathcal{O}$. On définit alors
 $\Upsilon_{a.p}:=\{(\mathcal{O}, \tau): \mathcal{O} \ \text{est
   admissible et bien polarisable},\ \tau\in X(\mathcal{O})\}$. Alors,
 la correspondance $$\Upsilon_{a.p}\longrightarrow
 \widehat{G}$$ $$(\mathcal{O},\tau)\longmapsto \pi_{f,\tau}$$ est bien
 définie, où pour tout $\mathcal{O}$, le point $f\in \mathcal{O}$ une
 fois choisi est fixé, de sorte que $ X(\mathcal{O})$ est identifié à
 $X(f)$. Duflo a montré que cette correspondance est injective.

Suposons que $G$ soit un groupe presque algébrique. Alors Duflo a
montré que $\pi_{f,\tau}$ est une série discrète modulo un sous-groupe
$\Gamma$ du centre de $G$ si et seulement si $G(f)/\Gamma$ est compact
et toute série discrète de $G$ modulo $\Gamma$ est isomorphe à une
telle $\pi_{f,\tau}$ (pour tout ceci voir [7], p. 203, théorème 3).

Si $G$ est connexe réductif et $f$ correspond à une série discrète,
alors $G(f)$ est connexe et $X(f)$ a un unique élément $\tau_{f}$ qui est un
caractère. Nous noterons simplement $T_{f}=T_{f,\tau_{f}}$.

Il est connu que les orbites fortement r\'{e}guli\`{e}res sont
\textbf{bien polarisables}. Pour $G$ réductif, $f$ est fortement
r\'{e}guli\`{e}re si et seulement si $\mathfrak{g}(f)$ est une
sous-algèbre de Cartan. De plus Duflo a d\'{e}montr\'{e} le
r\'{e}sultat suivant: Soit $G$ un groupe presque
alg\'{e}brique. Notons $\widehat{G}_{fr}$ l'ensemble des classes des
repr\'{e}sentations unitaires irréductibles de $G$ li\'{e}es aux
orbites fortement r\'{e}guli\`{e}res admissibles. Alors le
compl\'{e}mentaire de $\widehat{G}_{fr}$ dans $\widehat{G}$ est
n\'{e}gligeable par rapport \`{a} la mesure de Plancherel. Donc toute
s\'{e}rie discr\`{e}te $\pi$ de $G$ est liée \`{a} une orbite
fortement r\'{e}guli\`{e}re admissible, puisque par rapport \`{a} la
mesure de Plancherel, $\{\pi\}$ a une mesure non nulle. De plus,
d'après ce qui précède, un élément de $\widehat{G}_{fr}$ est une
série discrète si et seulement s'il correspond à une orbite dont le
stabilisateur d'un point quelconque est compact.

Nous dirons qu'un forme linéaire $f\in\mathfrak{g}^{*}$ est {\it de
  type compact} si $G(f)$ est compact. Nous dirons qu'une orbite
coadjointe est {\it de type compact}, si elle contient une forme
linéaire de type compact. On voit donc que les séries discrètes d'un
groupe presque algébrique sont les représentations associées par la
construction de Duflo aux orbites fortement régulières admissibles et
de type compact.\\

\begin{center}\section{Quelques propri\'{e}t\'{e}s de $SU(2,1)$
  }\label{SU21}\end{center}

 Dans tout ce qui suit, on note $G$ le groupe $SU(2,1)$ et
 $\mathfrak{g}$ son algèbre de Lie.\\

Désignons par $I_{2,1}$ la matrice dans la base canonique de la forme
hermitienne sur $\mathbb{C}^{3}$, $\langle
z,t\rangle=z_{1}\overline{t}_{1}+z_{2}\overline{t}_{2}-z_{3}\overline{t}_{3}$.
Alors
$$G=\left\{\left. g \in SL(3,\C) \right\vert {^t
  \overline{g}}I_{2,1}g=I_{2,1} \right\}, $$
$$
\mathfrak{g}=\left\{\left. X\in M_{3}(\mathbb{C}) \right\vert {^t
  \overline{X}}I_{2,1}+I_{2,1}X=0\mbox{ et } tr X=0 \right\}.
$$ Le groupe de Lie $G$ est simple de dimension 8, son complexifié est
$SL(3,\C)$. Le centre  de $G$ est le groupe à trois éléments
\begin{equation*}
Z_{G}=\{\xi Id\,\vert\, \xi^{3}=1\}.
\end{equation*}

Notons $\theta$ l'involution de Cartan sur $G$ et
$\mathfrak{g}$, telle que $\theta(g)= {^t\overline{g^{-1}}}$ pour
$g\in G$ et $\theta(X)=-{^t\overline{X}}$ pour
$X\in\mathfrak{g}$. Désignons par $K$ le sous-groupe des points fixes
de $\theta$ dans $G$, par $\mathfrak{k}$ son algèbre de Lie et par
$\mathfrak{p}$ le sous-espace propre de $\theta$ pour la valeur prore
$-1$ dans $\mathfrak{g}$.  Alors $K$ est un sous-groupe compact
maximal de $G$ et on a les décompositions de Cartan de $G$ et
$\mathfrak{g}$ :
$$
G=K\exp\mathfrak{p},
$$
$$
\mathfrak{g}=\mathfrak{k}\oplus\mathfrak{p}.
$$

Le centre de $z(\mathfrak{k})$ de $\mathfrak{k}$ est de dimension
$1$. Il est engendré par la matrice
$$Z= i\left(
\begin{array}{ccc}
1 & 0 & 0 \\
0& 1 &0\\
0 & 0  & -2 \\
\end{array}
\right) $$
et son centralisateur dans $\mathfrak{g}$ est $\mathfrak{k}$. Si $Z_{K}$
désigne le centre de $K$, on a
$$
Z_{K}=\exp\R Z.
$$

Soit $\mathfrak{k}'$ l'algèbre dérivée de $\mathfrak{k}$. Elle est
isomorphe à $\mathfrak{su}(2)$ :
$$
\mathfrak{k}'=\left\{
\left.
\left(
\begin{array}{cc}
A & 0 \\
0   & 0 \\
\end{array}
\right)
\right \vert
A \in \mathfrak{su}(2)\right\}$$
et on a $\mathfrak{k}=\mathfrak{k}'\oplus z(\mathfrak{k})$.
Le sous-groupe analytique $K'$ d'algèbre de Lie $\mathfrak{g}$ de $G$
est le sous-groupe dérivé de $K$. Il est isomorphe à $SU(2)$ et on a
$$
K'=\left\{
\left.
\left(
\begin{array}{cc}
U & 0 \\
0   & 1 \\
\end{array}
\right)
\right \vert
U \in SU(2)\right\}
\quad\text{et}\quad K=K'Z_{K}.
$$
Enfin, on a
$$
\mathfrak{p}=\left\{
\left.
\left(
\begin{array}{cc}
0 & X \\
 ^t \overline{X}  & 0 \\
\end{array}
\right)
\right \vert
X \in M_{2,1}(\C)\right\}.$$

Soit $\mathfrak{t}$ la sous-algèbre de Lie de $\mathfrak{g}$
constituée des matrices diagonales. C'est une sous-algèbre de Cartan
de $\mathfrak{g}$ contenue dans $\mathfrak{k}$. C'est donc également
une sous-alg\`{e}bre de Cartan de $\mathfrak{k}$. Désignons par
$diag(h_1,h_2,h_3)$ la matrice diagonale de $M_{3}(\mathbb{C})$ de
valeurs propres $h_{1},h_{2},h_{3}$, dans la base canonique de
$\mathbb{C}^{3}$. Alors, on a
$$\mathfrak{t}=\left\{ \left.  diag(ih_{1},ih_{2},ih_{3}) \right
\vert h_1 , h_2, h_{3} \in \R\mbox{ et }h_{1}+h_{2}+h_{3}=0
\right\},$$
$$
\mathfrak{t}_{\mathbb{C}}=\left\{ \left.  diag(h_{1},h_{2},h_{3}) \right
\vert h_1 , h_2, h_{3} \in \C\mbox{ et }h_{1}+h_{2}+h_{3}=0
\right\}.
$$

Le syst\`{e}me des racines de $\mathfrak{t}_{\C}$ dans
$\mathfrak{g}_{\C}$ est
$$ \Sigma(\mathfrak{g}_\C,\mathfrak{t}_\C)=\left\{ \left. \alpha_{kl}
\right \vert 1 \leqslant k\neq l\leqslant 3\right\},$$ où
$$ \alpha_{kl}(diag(h_1,h_2,h_3))=h_k-h_l.$$
Le sous-espace radiciel de $\mathfrak{g}$ pour la racine $\alpha_{kl}$
est $\mathfrak{g}_{kl}=\C E_{kl}$, où $E_{kl}\in M_{3}(\C)$ désigne la
matrice élémentaire d'indice $kl$. La coracine de $\alpha_{kl}$ est la
matrice $H_{kl}=E_{kk}-E_{ll}$. Posons $H=iH_{12}$. Alors
$\mathfrak{t}':=\mathfrak{t}\cap\mathfrak{k}'=\R H$,
$\mathfrak{t}=\mathfrak{t}'\oplus z(\mathfrak{k})$
et $(H,Z)$ est une base de $\mathfrak{t}$.

De plus, $\left\{\alpha_{12},\alpha_{21}\right\}$  est l'ensemble des
racines compactes et son complémentaire
$ \Sigma(\mathfrak{g}_\C,\mathfrak{t}_\C)\setminus
\left\{\alpha_{12},\alpha_{21}\right\}$ est l'ensemble des
  racines non-compactes.

Maintenant, si $ S= \left(
\begin{array}{ccc}
0 & 0&1 \\
 0  & 0&0 \\
1&0&0
\end{array}
\right)$, $\mathfrak{a}=\R S$ est une sous-alg\`{e}bre ab\'{e}lienne
maximale de $\mathfrak{p}$. L'ensemble des racines $\Delta$ de
$\mathfrak{a}$ dans $\mathfrak{g}$ est $\pm\beta$, $\pm2\beta$ avec
$\beta( tS)=t, \ t\in \R$. Les sous-espaces radiciels correspondants
sont

$$\mathfrak{g}_{\beta}=\R E_1\bigoplus \R E_1',\ \ \mathfrak{g}_{2\beta}=\R E_2,$$

$$\mathfrak{g}_{-\beta}=\theta(\mathfrak{g}_{\beta}),\ \ \mathfrak{g}_{-2\beta}=\theta(\mathfrak{g}_{2\beta}),$$
où "$\theta$" est l'involution de Cartan (au niveau alg\'{e}brique) et

$$
E_1=\left(
\begin{array}{ccc}
0& -1&0 \\
 1 & 0&-1\\
0&-1&0
\end{array}
\right)
\quad E_1'=\left(
\begin{array}{ccc}
0& -i&0 \\
 -i & 0&i\\
0&-i&0
\end{array}
\right)
\quad E_2=i\left(
\begin{array}{ccc}
2& 0&-2\\
 0 & 0&0\\
2&0&-2
\end{array}
\right).
$$

Notons $\mathfrak{n}=\mathfrak{g}_{\beta}+\mathfrak{g}_{2\beta}$ et
$N=\exp(\mathfrak{n})$.  Alors les décompositions d'Iwasawa
correspondantes sont

$$\mathfrak{g}=\mathfrak{k}\oplus\mathfrak{a}\oplus\mathfrak{n}\  \text{et} \   G=KAN.$$

Il est clair que $[E_1,\ E_2]=[E_1',\ E_2]=0$ et
$[E_1,\ E_1']=E_2$. Donc $\mathfrak{n}$ est une alg\`{e}bre de
Heisenberg de dimension 3 et $N$ est un groupe de Heisenberg de
dimension 3.

$$\text{Si}\ W=\frac{i}{3}\left(
\begin{array}{ccc}
1& 0&0 \\
 0 & -2&0\\
0&0&1
\end{array}
\right) ,$$

\noindent alors $\mathfrak{m}=\R W$ est le centralisateur de
$\mathfrak{a}$ dans $\mathfrak{k}$ . De plus $[W,\ E_1]=E_1'$,
$[W,\ E_1']=-E_1$, et $[W,\ E_2]=0$. Donc
$\mathfrak{b}=\mathfrak{m}\oplus \mathfrak{a}\oplus \mathfrak{n}$ est
une sous-alg\`{e}bre de Borel. Notons $M=\exp(\mathfrak{m})$ qui est
le centralisateur du sous-groupe $A$ dans $K$. Ainsi, $B=MAN$ est le
sous-groupe de Borel associé à $\mathfrak{b}$. Pour simplifier, on
note $\mathfrak{b}_1=\mathfrak{a}\oplus \mathfrak{n}$ et $B_1=AN$.

On vérifie que $Z_{G}$ est également le centre de $B$.\\

\pagebreak
\begin{center}{\section{Rappels sur les séries discrètes des
groupes de Lie simples}}\end{center}

\noindent\begin{theo}\label{series-discrètes} Soient $G$ un groupe de
Lie lin\'{e}aire semi-simple, $K\subset G$ un sous-groupe compact
maximal, $\mathfrak{g}$ et $\mathfrak{k}$ leurs algèbres de Lie
respectives. On suppose que $K$ et $G$ ont le même rang. Soit donc
$\mathfrak{t}\subseteq \mathfrak{k}\subseteq \mathfrak{g}$ une
sous-alg\`{e}bre de Cartan de $\mathfrak{g}$. Notons $\bigtriangleup
=\Sigma (\mathfrak{g}_{\C},\ \mathfrak{t}_{\C})$ et $\bigtriangleup_K
=\Sigma (\mathfrak{k}_{\C},\ \mathfrak{t}_{\C})$.

 Supposons que $\lambda \in (i\mathfrak{t})^*$ est non-singulière
 (autrement dit $ \left\langle \lambda,\ \alpha \right\rangle \neq0$
 pour toute $\alpha \in \Delta$ ). On considère l'ensemble de racines
 positives (pour $\Delta$), $\Delta^+=\left\{\left.\alpha \in \Delta
 \right\vert \left\langle\lambda,\ \alpha \right\rangle>0
 \right\}$. Alors, $\Delta_{K}^+ = \Delta^+\bigcap\Delta_K$ est un
 ensemble de racines positives pour $\Delta_K$. Notons $\delta_G$
 (resp. $\delta_K$) la demi-somme des éléments de $\Delta^+$
 (resp. $\Delta_K^+ $).

Si $\lambda+\delta_G$ est analytiquement int\'{e}grable (i.e. est la
différentielle d'un caractère du tore maximal $T$), alors il existe
une unique s\'{e}rie discr\`{e}te $\pi_\lambda$ de $G$ qui poss\`{e}de les
propri\'{e}t\'{e}s suivantes:

(a)  $\pi_\lambda$ a le caract\`{e}re infinit\'{e}simal
$\zeta_\lambda$, où $\zeta_\lambda$ est le caractère du centre de
l'algèbre enveloppante de $\mathfrak{g}$ défini par $\lambda$.

(b) Dans la restriction de $\pi_\lambda$ \`{a} $K$ que l'on note
$\pi_\lambda \vert_K $, le $K$ type avec le plus haut poids (par
rapport \`{a} $\Delta_K^+$) $$\Lambda=\lambda+\delta_G-2\delta_K$$
figure une seule fois.

(c) si $\Lambda'$ est le plus haut poids d'un $K$ type qui
appara\^{\i}t dans $\pi_\lambda \vert_K $, alors, $\Lambda'$ est de la
forme $$\Lambda'= \Lambda+\sum_{\alpha \in \Delta^+} n_\alpha \alpha
\ \ \ \ \text{ici}\ n_\alpha \ \text{sont des entiers positifs ou
  nuls}.$$

De plus, deux telles $\pi_{\lambda_1}$ et $\pi_{\lambda_2}$ sont
\'{e}quivalentes si et seulement si $\lambda_1$ et $\lambda_2$ sont
conjugu\'{e}es par $W_K$, où $W_K$ est le groupe de Weyl pour
$\Delta_K$. Et, toute s\'{e}rie discr\`{e}te de $G$ est obtenue de
cette mani\`{e}re \end{theo}

\noindent \textbf{Remarque 1}. (1) La paramétrisation des séries
discrètes  est due à Harish-Chandra. Les assertions b) et c) du
théorème sont une partie de la conjecture de Blattner, démontrée par
Hecht et Schmid.

(2) Dans le th\'{e}or\`{e}me, le
param\`{e}tre $\lambda$ est appelé \textbf{le param\`{e}tre de
  Harish-Chandra}, et le param\`{e}tre $\Lambda$ est appelé \textbf{le
  param\`{e}tre de Blattner}.

(3) Supposons que $G$ est simple et que le centre de $K$ est de
dimension supérieure ou égale à $1$. Alors $\pi_{\lambda}$ est une
série discrète holomorphe si et seulement si
$\mathfrak{p}^{+}_{\C}=\oplus_{-\alpha\in
  \Delta^+_{\text{n}}}\mathfrak{g}^{\alpha}_{\C}$ est une sous-algèbre
(abélienne) de $\mathfrak{g}_{\C}$, où $\Delta^+_{\text{n}}$ est
l'ensemble des racines non compactes positives (relativement à
$\lambda$). De plus $\mathfrak{g}^{\alpha}_{\C}$ est l'espace radiciel
de $\alpha$. Dans ce cas $\pi_{\lambda}$ est holomorphe pour la
structure complexe définie par $\mathfrak{p}^{+}_{\C}$ sur $G/K$.

(4) Si $f\in \mathfrak{g}^*$ est une forme linéaire fortement
régulière, d'après ce que l'on a expliqué dans le chapitre 2,
$\mathfrak{g}(f)$ est une sous-algèbre de Cartan. Supposons $f$ de
type compact. Sans perte de généralité, on peut supposer
$\mathfrak{g}(f)=\mathfrak{t}$.  Donc $f\in \mathfrak{t}^*\subset
\mathfrak{g}^*$ (on utilise la décomposition en somme directe
$\mathfrak{g}=\mathfrak{t}\oplus[\mathfrak{t},\mathfrak{g}]$ pour
identifier $\mathfrak{t}^*$ à un sous-epace de $\mathfrak{g}^*$). On
peut donc définir l'ensemble de racines positives $\Delta^{+}$ relatif
à $if\in (i\mathfrak{t})^* $ comme on a fait pour $\lambda$ dans le
théorème précédent. Soit $\delta$ (resp. $\delta_{K}$) la demi-somme
des racines dans $\Delta^{+}$
(resp. $\Delta^{+}_{K}:=\Delta^{+}\cap\Delta_{K}$).

Ainsi $f$ est admissible pour $G$ si et seulement si $\lambda=if$ est
un paramètre de Harish-Chandra. Supposons que nous soyons dans ce
cas. Alors, la représentation $T_{f}$ de la paramétrisation de Duflo est la
série discrète $\pi_{\lambda}$.

Soit $\chi_{f}$ la caractère de $T$ de
différentielle $if+\delta-2\delta_{K}$. On voit facilement que le caractère
central de la série discrète $T_{f}=\pi_{\lambda}$ est la restriction
de $\chi_{f}$ au centre $Z_{G}$ de $G$.

(5) Soit $f\in \mathfrak{g}^*$ est une forme linéaire fortement
régulière de type compact. Avec les notations de (4), on dit que $f$
est dans \textbf{le cône holomorphe}, si
$\mathfrak{p}^{+}_{\C}=\oplus_{-\alpha\in
  \Delta^+_{\text{n}}}\mathfrak{g}^{\alpha}_{\C}$ est une sous-algèbre
(abélienne) de $\mathfrak{g}_{\C}$, où $\Delta^+_{\text{n}}$ est
l'ensemble des racines non compactes dans $\Delta^{+}$. Donc il est
clair que si $f$ est de plus admissible, la série discrète $T_{f}$ est
holomorphe si et seulement si $f$ est dans le cône holomorphe.\\

Maintenant, on revient sur $G=SU(2,1)$ et on reprend du chapitre 2 le
concernant. On veut trouver à quelles conditions un élément
$\lambda\in\mathfrak{t}^{*}$ correspond \`{a} une s\'{e}rie
discr\`{e}te et à quelles conditions cette dernière est holomorphe.

Il est clair que $W_K=\left\{s_{\alpha_{12}},\ id\right\}$.  Donc,
d'apr\`{e}s le th\'{e}or\`{e}me, on peut toujours supposer que
$\alpha_{12}$ est dans $\Delta^+$. Les ensembles de racines positives
contenant $\alpha_{12}$ sont $$\Delta_1^+=
\left\{\alpha_{12},\ \alpha_{32},\ \alpha_{31}\right\},\ \Delta_2^+=
\left\{\alpha_{12},\ \alpha_{23},\ \alpha_{13}\right\},
\ \Delta_3^+=\left\{\alpha_{12},\ \alpha_{32},\ \alpha_{13}\right\}$$
D'après la remarque (3) précédente, on peut déduire que
$\Delta_1^+$ et $\Delta_2^+$ correspondent aux s\'{e}ries
discr\`{e}tes holomorphes et $\Delta_3^+$ à celles non
holomorphes. \\

\noindent \textbf{Remarque 2}.  (1) Comme $\delta_{G}=
1/2(\alpha_{12}+\alpha_{32}+\alpha_{31})= \alpha_{32} $ est un poids,
$\lambda$ correspond à une série discrète si et seulement
si $\lambda$ est analytiquement int\'{e}grable.

(2) Par convention, on appelle \og s\'{e}ries discr\`{e}tes
holomorphes \fg, celles qui correspondent \`{a} l'un des ensembles de
racines positives $\Delta_1^+$ ou $\Delta_2^+$, et \og s\'{e}ries
discr\`{e}tes anti-holomorphes \fg, celles qui correspondent à
l'autre. Dans la suite, pour nous conformer à l'article de Rossi et
Vergne ([29]) que l'on va utiliser, on appellera \og s\'{e}ries
discr\`{e}tes holomorphes \fg, celles correspondant \`{a}
$\Delta_1^+$.\\

D'apr\`{e}s le th\'{e}or\`{e}me, les \og$\lambda$ \fg~  qui
correspondent \`{a} $\Delta_1^+$ sont celles qui v\'{e}rifient que $
\left\langle \lambda,\ \alpha \right\rangle >0$ pour toute $\alpha \in
\Delta_1^+$ et $\lambda+\delta_{G}$ est analytiquement
int\'{e}grable.  On en d\'{e}duit donc
que les \og$\lambda$ \fg~ qui correspondent à $\Delta_1^+$ sont celles
qui v\'{e}rifient
 $$\lambda(H_{12})=n_1 \in \N^+ ,\ \lambda(H_{32})=n_2 \in \N^+
,\ \lambda(H_{31})=n_3 \in\N^+.$$ Mais, comme $H_{32}=H_{12}+H_{31}$,
et $H_{12}$ et $H_{31}$ sont lin\'{e}airement ind\'{e}pendants, la
condition pr\'{e}c\'{e}dente est \'{e}quivalente \`{a}
$$\lambda(H_{12})=n_1 \in \N^+,\ \lambda(H_{31})=n_3 \in \N^+ \ \ \text{o\`{u}}\ \text{le choix de}\  n_1 \ \text{et}\ n_3 \ \text{est ind\'{e}pendant}.$$

\pagebreak

\begin{center}\section{\'Etude des orbites coadjointes de $G$, $B$ et
    $B_{1}$}\end{center}


Nous commençons par rappeler ce que nous entendons par propreté et
faible de propreté pour la projection d'une orbite coadjointe. On se
donne un groupe de Lie presque algébrique connexe $G$ d'algèbre de Lie
$\mathfrak{g}$ et $H\subset G$ un sous-groupe presque algébrique
d'algèbre de Lie $\mathfrak{h}$.  On désigne par $p$ la projection
naturelle de $\mathfrak{g}^{*}$ sur $\mathfrak{h}^{*}$.  Soit
$\mathcal{O}\subset\mathfrak{g}^{*}$ une orbite coadjointe sous $G$.

Nous dirons que la restriction de la projection $p$ à $\mathcal{O}$
est {\it propre}, si elle est propre sur l'image, c'est à dire si,
pour tout compact $L\subset p(\mathcal{O})$,
$p^{-1}(L)\cap\mathcal{O}$ est compact.

Nous dirons que la restriction de la projection $p$ à $\mathcal{O}$
est {\it faiblement propre},  si,
pour tout compact $L\subset p(\mathcal{O})\cap\mathfrak{h}_{fr}^{*}$,
$p^{-1}(L)\cap\mathcal{O}$ est compact.\\

 Pour simplifier, dans tout ce qui suit, sauf indication contraire, on
 garde toutes les notations concernant $SU(2,1)$ du chapitre précédent
 (i.e. on note $G=SU(2,1)$ etc). On désigne par $\text{p}$
 (resp. $\text{p}_1$) la projection naturelle de $\mathfrak{g}^{*}$
 sur $\mathfrak{b}^{*}$ (resp. $\mathfrak{b}_{1}^{*}$).\\

 Maintenant, on va énoncer le plan de ce chapitre. \\

 1) Détermination des orbites coadjointes fortement régulières et
 celles fortement régulières et admissibles de $\mathfrak{b}_1^*$ et
 détermination des données d'admissibilité.

 2) La même chose pour $\mathfrak{b}^*$.

 3) \'Etude de la faible propreté (et de la propreté) pour la restriction
 de $\text{p}_1$ à une orbite fortement régulière de type compact.

 4) La même chose pour $p$.\\

 \subsection{Orbites coadjointes fortement régulières (et admissibles) de
$\mathfrak{b}_1^*$ }

 Maintenant, $S,\ E_1,\ E_1',\ E_2$ sont comme dans le chapitre 3, ils
 constituent une base de $\mathfrak{b}_1$. Notons
 $S^*,\ E_1^*,\ E_1'^*,\ E_2^*$ la base duale correspondante dans
 $\mathfrak{b}_1^*$.

 \begin{prop}\label{orbites-B1}
 Notons $\Omega^-=\{b\in\mathfrak{b}_1^*\mid b(E_2)<0\} $ et
 $\Omega^+=\{b\in\mathfrak{b}_1^*\mid b(E_2)>0\}$. Alors
 $\Omega^{\pm}=B.(\pm E_2^*)$ sont les seules $NA$-orbites coadjointes
 fortement r\'{e}guli\`{e}res dans $\mathfrak{b}_1^*$, et elles sont
 également admissibles.
 \end{prop}
\noindent\begin{demo} Notons
${\mathfrak{b}_1^*}'=\{b\in\mathfrak{b}_1^*\mid b(E_2)\neq0\}$. Comme
$\R E_{2}$ est un idéal de $\mathfrak{b}_{1}$, il est clair que
${\mathfrak{b}_1^*}'$ est un ouvert $B_1$-invariant (pour l'action
coadjointe). Soit $b\in{\mathfrak{b}_1^*}'$. Nous allons voir que
$\mathfrak{b}_1(b)=0$. Soit $n=b_{\vert\mathfrak{n}}$. Comme l'idéal
$\mathfrak{n}$ de $\mathfrak{b}_{1}$ est une alg\`{e}bre de Lie de
Heisenberg de dimension 3 de centre $E_{2}$, on a $N.n=n+(\R
E_2)^\bot$, où $(\R E_2)^\bot=\left\{\left.g \in n^*\right\vert
g(E_2)=0\right\}$. On peut supposer que $b=\mu S^* +\gamma E_2^*$
avec $\gamma \neq 0$. Dans ce cas, la matrice de la forme $B_{b}$ dans
la base $S,\ E_1,\ E_1',\ E_2$ est $$\left(
\begin{array}{cccc}
0& 0&0&2\gamma \\
 0 & 0&\gamma&0\\
0&-\gamma&0&0\\
-2\gamma&0&0&0
\end{array}
\right) $$

On en déduit donc que les $B_1$-orbites coadjointes dans
${\mathfrak{b}_1^*}'$ sont ouvertes. Comme elles sont connexes, par
argument de connexité, il y en deux : $\Omega^-$ et $\Omega^+$. De
plus comme $\Omega^{\pm}$ est simplement connexe, si $b\in
\Omega^{\pm}$ on déduit que $$B_1\longrightarrow
\Omega^{\pm}$$ $$x\longmapsto x.b$$ est un difféomorphisme, de sorte
que le stabilisateur $B_1(b)$ est trivial. Donc $\Omega^-$ et
$\Omega^+$ sont fortement r\'{e}guli\`{e}res et admissibles, d'où le
résultat.
\end{demo}

\subsection{Orbites coadjointes fortement régulières (et admissibles)
  de $\mathfrak{b}^*$ }

$W$, $S$, $E_1$, $E_1'$ et $E_2$ est une base de $\mathfrak{b}$ et on
note $W^*$, $S^*$, $E_1^*$, $E_1'^*$ et $E_2^*$ la base duale de
$\mathfrak{b}^*$.

\begin{prop}\label{orbites-B}
Pour $r\in \R$, notons $\Omega_{r,-}=B.(rW^*-E_2^*)\subset
\mathfrak{b}^*$ et $\Omega_{r,+}=B.(rW^*+E_2^*)\subset
\mathfrak{b}^*$, alors

1) Les $B$-orbites fortement r\'{e}guli\`{e}res dans $\mathfrak{b}^*$
sont les $\Omega_{r,-}$ et $\Omega_{r,+}$. Elles sont deux à deux
distinctes.

2) Les $B$-orbites fortement r\'{e}guli\`{e}res et admissibles dans
$\mathfrak{b}^*$ sont les $\Omega_{r,-}$ et $\Omega_{r,+}$, avec
$r+\frac{1}{2} \in \Z/3$.
\end{prop}
\noindent\begin{demo}On va commencer par d\'{e}terminer toutes les
orbites fortement r\'{e}guli\`{e}res de $B=MAN$:
$\mathfrak{b}=\mathfrak{m}\oplus\mathfrak{a}\oplus\mathfrak{n}$. L'algèbre
de Lie $\mathfrak{b}$ est de dimension $5$ et elle contient
$\mathfrak{b}_{1}$ comme un idéal. Il suit alors de la proposition
\ref{orbites-B1} de la section précédente, que si $f\in
\mathfrak{b}^*$ vérifie $f(E_2)\neq 0$, on a $\text{dim}B.f=4$ . Soit
donc $f\in \mathfrak{b}^*$ telle que $f(E_2)\neq 0$. Quitte à
translater $f$ par un élément de $B_{1}$, on peut supposer que
$f=rW^*+\varepsilon E_2^*$ avec $\varepsilon\in \{\pm 1\}$. Il est
alors immédiat que $\mathfrak{b}(f)\subset\mathfrak{m}$ et donc que
$\mathfrak{b}(f)=\mathfrak{m}$, pour des raisons de dimension. On en
déduit que les orbites $\Omega_{r,-}$ et $\Omega_{r,+}$ sont fortement
r\'{e}guli\`{e}res. Maintenant si $f(E_2)=0$, alors $\R
E_2\subset\mathfrak{b}(f)$. Donc $f$ ne peut pas être fortement
r\'{e}guli\`{e}re. Par suite, $\mathfrak{b}_{fr}^*=\{f\in
\mathfrak{b}^*\mid f(E_2)\neq 0 \}$ et $\Omega_{r,-}$ et
$\Omega_{r,+}$, $r\in \R$ sont les seules $B$-orbites fortement
r\'{e}guli\`{e}res dans $\mathfrak{b}^*$. Comme $B=MB_{1}$ et comme le
stabilisateur d'un point de $\Omega_{\pm}$ dans $B_{1}$ est trivial,
on en déduit que les formes linéaires $rW^*+\varepsilon E_2^*$ et
$r'W^*+\varepsilon' E_2^*$ avec $r,r'\in\R$ et
$\varepsilon,\varepsilon'\in \{\pm 1\}$, ne peuvent être dans la même
$B$-orbite que si $r=r'$ et $\varepsilon=\varepsilon'$. Par suite, les
orbites $\Omega_{r,\pm}$ sont deux à deux distinctes.

Maintenant, on va d\'{e}terminer toutes les $B$-orbites fortement
r\'{e}guli\`{e}res admissibles dans $\mathfrak{b}^*$. On a
d\'{e}j\`{a} vu que les $B$-orbites fortement r\'{e}guli\`{e}res (dans
$\mathfrak{b}^*$) sont les $\Omega_{r,\pm} \ {r\in \R}$ . Soit donc
$f=rW^*+\varepsilon E_2^*\in b^*$. On veut chercher les conditions,
pour lesquelles $f$ est admissible. Par des calculs directs, on voit
que le stabilisateur de $f$ dans $B$ est $B(f)=\exp(\R W)$ et que
$\xi=\C W\oplus \C (E_1+\varepsilon iE_1')\oplus \C E_2$ est une
polarisation positive pour la forme bilin\'{e}aire altern\'{e}e
$B_f$. Donc $\rho_\xi(W)=\frac{1}{2}\text{tr}(W_{\xi/\C
  W})=-\varepsilon\frac{i}{2}$ (car $[W,E_1+\varepsilon
  iE_1']=-\varepsilon i(E_1+\varepsilon iE_1')$ et $[W,E_2]=0$). Donc
en appliquant le lemme \ref{admissibilité}, on d\'{e}duit que $f$ est
admissible si et seulement s'il existe un caract\`{e}re $\chi$ de
$\exp(\R W)$ dont la diff\'{e}rentielle est $\rho_\xi+if\vert_{\R
  W}$. Comme $(\rho_\xi+if\vert_{\R
  W})(wW)=iw(r-\varepsilon\frac{1}{2}),\ w\in \R$, si un tel caractère
$\chi$ existe, on a
$\chi(\exp(wW))=e^{i(r-\varepsilon\frac{1}{2})w},\ w\in \R$. Or il est
clair que $\exp(wW)=1$ si et seulement si $w\in 6\pi \Z$. Donc un tel
caract\`{e}re existe, si et seulement si $6\pi
i(r-\varepsilon\frac{1}{2})\in 2 \pi i \Z$, c'est-à-dire
$r+\frac{1}{2} \in \Z/3$. On en d\'{e}duit que \textbf{les
  $B$-orbites fortement r\'{e}guli\`{e}res et admissibles dans
  $\mathfrak{b}^*$ sont les $\Omega_{r,-}$ et $\Omega_{r,+}$, avec
  $r+\frac{1}{2} \in \Z/3$}, d'où le résultat.
\end{demo}\\

Pour ce dont on a besoin ultérieurement, on va décrire explicitement
les orbites $\Omega_{r,\pm}$. D'abord on a $B.(rW^*+\varepsilon
E_2^*)B_{1}M(rW^*+\varepsilon E_2^*)=B_{1}.(rW^*+\varepsilon
E_2^*)$. Soit donc $b=\exp(xE_1+yE_1')\exp(zE_2).\exp(sS)$ dans $B_{1}$,
avec $ x,\ y$, $z,\ s\in \R$. Comme $[W,E_1]=E_1'$, $[W,E_1']=-E_1$,
$[W,E_2]=[W,S]=0$, on en d\'{e}duit que
$\exp(-xE_1-yE_1').W=W+[W,xE_1+yE_1']+\frac{[[W,xE_1+yE_1'],xE_1+yE_1']}{2}=W+xE_1'-yE_1-(\frac{x^2+y^2}{2})E_2$. Donc
$b^{-1}.W=W+xe^{-s}E_1'-ye^{-s}E_1-(\frac{x^2+y^2}{2}e^{-2s})E_2$ et
$b.(rW^*+\varepsilon
E_2^*)(W)=r-\varepsilon(\frac{x^2+y^2}{2})e^{-2s}$, ici
$\varepsilon\in \{1,-1\}$. D'autre part, on vérifie que l'on a
$b.E_{2}^{*}=(2z e^{-2s})S^*+ye^{-2s} E_1^*-xe^{-2s}E_1'^*+ e^{-2s}E_2^*$.
On en d\'{e}duit que

\begin{equation*}
\begin{split}
b(rW^*+\varepsilon
E_2^*)&=(2\varepsilon z e^{-2s})S^*+\varepsilon ye^{-2s}
E_1^*-\varepsilon xe^{-2s}E_1'^*+\varepsilon
e^{-2s}E_2^*\\ &\quad+(r-\varepsilon(\frac{x^2+y^2}{2}e^{-2s}))W^*.
\end{split}
\end{equation*}
Ceci montre que l'on a
\begin{equation*}
\Omega_{r,\pm}=\{wW^*+sS^*+xE_1^*+yE_1'^*+zE_2^*\,\vert
x,\ y,\ s\in \R ,\ \pm z>0, w=r-\frac{x^2+y^2}{2z} \}.
\end{equation*}
Donc on voit facilement que deux formes lin\'{e}aires de
$\mathfrak{b}^*$,
$f_i=w_iW^*+s_iS^*+x_iE_1^*+y_iE_1'^*+z_iE_2^*,\ i=1,\ 2$ avec
$z_1,z_2 \neq0$ sont dans la m\^{e}me orbite, si et seulement si
$z_1z_2>0$ et
$w_1+\frac{x_1^2+y_1^2}{2z_1}=w_2+\frac{x_2^2+y_2^2}{2z_2}$. Donc
\textbf{$w_1W^*+\varepsilon_1 E_2^*$ et $w_2W^*+\varepsilon_2E_2^*$
  avec $\varepsilon_1, \varepsilon_2\in \{\pm 1\} $ sont dans la
  m\^{e}me orbite, si et seulement si $\varepsilon_1=\varepsilon_2$ et
  $w_1=w_2$, et $f=wW^*+sS^*+xE_1^*+yE_1'^*+zE_2^*\in \mathfrak{b}^*$
  avec $z\neq0$, est dans la m\^{e}me orbite que
  $(w+\frac{x^2+y^2}{2z})W^*+\varepsilon E_2^*$ avec $\varepsilon z>
  0$}.

Une conséquence importante est que l'on obtient un difféomorphisme :

$$\Psi: \R \times AN\times \{\pm1\}\longrightarrow \mathfrak{b}_{fr}^*$$

$$(r,x,\varepsilon)\longmapsto x.(rW^*+\varepsilon E_2^*).$$

Pour $f\in \mathfrak{b}_{fr}^*$, on écrit
$\Psi^{-1}(f)=(r(f),x(f),\varepsilon(f))$. Alors on a
$r(f)=w+\frac{x^2+y^2}{2z}$ et $\varepsilon(f)=\frac{z}{\mid z\mid}$
pour $f=wW^*+sS^*+xE_1^*+yE_1'^*+zE_2^*\in \mathfrak{b}_{fr}^*$.

 \subsection{La faible propreté (et la propreté) de la projection $\text{p}_1$}\label{propreté-p1}

Supposons que $f_0\in \mathfrak{g}^*$ est fortement régulière de type
compact (pas forcément admissible). Donc d'après ce que l'on a
expliqué, sans perte de généralité, on peut supposer que
$\mathfrak{g}(f_0)=\mathfrak{t}$, auquel cas $f_{0}$ est de la forme
$f_0=f_0(H)H^*+f_0(Z)Z^*\in \mathfrak{t}^*\subset\mathfrak{g}^*$, avec
$|f_0(H)|\neq0 \ \text{et}\ |f_0(H)| \neq | f_0(Z)|$ (ici,
$(H^{*},Z^{*})$ désigne la base duale de la base de $\mathfrak{t}$
introduite dans la section \ref{SU21}). 
D'autre part, on voit facilement que $f_0$ est dans le cône holomorphe
(que l'on a défini dans le chapitre précédent) si et seulement si
$|f_0(H)|\neq0 \ \text{et}\ |f_0(H)| < | f_0(Z)|$.

Maintenant, on va déterminer $\text{p}_1(G.f_0)$. Comme $G=NAK$, il est
\'{e}vident que
$\text{p}_1(\mathcal{O}_{\pi_\lambda})=\text{p}_1(NAK.f_0)=NA.{\text{p}_1(K.f_0)}$
et $K.f_0=f_0(H)K.H^*+f_0(Z)Z^*\subset \mathfrak{g}^*$.

 On identifie $\mathfrak{g}^*$ à $\mathfrak{g}$
 par $$\mathfrak{g}\longrightarrow \mathfrak{g}^*$$ $$X\longmapsto
 (Y\longmapsto-\frac{1}{2}\text{tr}(XY)),$$ où $X,\ Y \in
 \mathfrak{g}$.  Puisque la forme "$\text{tr}(XY)$" est
 proportionnelle à la forme de Killing, et $\text{tr}(H^2)=-2$ ,
 $\text{tr}(Z^2)=-6$, on en déduit que sous cette identification,
 $H^*\cong H$ et $Z^*\cong \frac{Z}{3}$.

Posons $$ \ F=\left( \begin{array}{ccc}
  0&1&0\\-1&0&0\\0&0&0 \end{array} \right),
\ V=\left( \begin{array}{ccc} 0&i&0\\i&0&0\\0&0&0 \end{array}
\right).$$ Alors, on a les relations $[H,F]=2V$, $[F,V]=2H$ et
$[V,H]=2F$. De plus, $(H,F,V)$ est une base orthonormée de l'algèbre
dérivée $\mathfrak{k}'$ de $\mathfrak{k}$. Comme le sous-groupe dérivé
$K'$ de $K$ est isomorphe à $SU(2)$ et  les orbites coadjointes de
$SU(2)$ sont de dimension $2$, on voit que
$$K.H^*\cong\left\{\left.x_1H+x_2F+x_3V\right\vert
x_1^2+x_2^2+x_3^2=1 \right\}.$$

Maintenant, $S,\ E_1,\ E_1',\ E_2$ sont comme dans le chapitre
précédent; ils constituent une base de $\mathfrak{b}_1$. Notons
$S^*,\ E_1^*,\ E_1'^*,\ E_2^*$ la base duale correspondante dans
$\mathfrak{b}_1^*$. Par des calculs directs, on montre que
\begin{equation}\label{image-base-p1}
H\vert_{\mathfrak{b}_1}=E_2^*\quad F\vert_{\mathfrak{b}_1}=-E_1^*
\quad V\vert_{\mathfrak{b}_1}=-E_1'^*\quad
Z^*\vert_{\mathfrak{b}_1}=\frac{Z}{3}\vert_{\mathfrak{b}_1}=E_2^*.
\end{equation}
Donc
$$\text{p}_1(K.H^*)=\left\{\left.x_1E_2^*+x_2E_1^*+x_3E_1'^* \right\vert x_1^2+x_2^2+x_3^2=1 \right\} $$ et
$$ \text{p}_1(Z^*)=E_2^*.$$
Ainsi, on a $$\text{p}_1(K.f_0)=\left\{\left.(f_0(H)x_1+f_0(Z))E_2^*+f_0(H)(x_2E_1^*+x_3E_1'^*) \right\vert x_1^2+x_2^2+x_3^2=1 \right\}.$$

\textbf{On en déduit donc
  que} $$\text{p}_1(G.f_0)\subset{\mathfrak{b}_1^*}'={\mathfrak{b}_1}_{fr}^*\Longleftrightarrow|f_0(H)|
<|f_0(Z)| $$

$$\Longleftrightarrow f_0 \ \text{est dans le cône holomorphe},$$ où
${\mathfrak{b}_1^*}'=\{b\in\mathfrak{b}_1^*\mid
b(E_2)\neq0\}$. Surtout dans ce cas, $\text{p}_1(G.f_0)$ est une
$B_1$-orbite coadjointe ouverte dans $\mathfrak{b}_1^*$.

Maintenant on va démontrer un théorème général qui permet de conclure
la faible propreté (et la propreté) de la projection $\text{p}_1$.

 \begin{theo}\label{propreté-G1}
Soit $G_1$ un groupe de Lie simple linéaire connexe d'algèbre de Lie
$\mathfrak{g}_1$. Soit $G_1=K_1A_1N_1$ une d\'{e}composition
d'Iwasawa, où $K_1$ est un sous-groupe compact maximal de $G_{1}$. On
suppose que $A_1N_1$ admet une orbite coadjointe ouverte dans
$(\mathfrak{a}_1\oplus\mathfrak{n}_1)^*$ (o\`{u} $\mathfrak{a}_1$ et
$\mathfrak{n}_1$ sont les alg\`{e}bres de Lie de $A_1$ et $N_1$
respectivement ). Soit $g\in\mathfrak{g}_1^*$ tel que le stabilisateur
$G_1(g)=T_1$ soit un tore compact. Alors pour que la projection
$\text{p}_{1}: \mathcal{O} := G_1.g \longrightarrow
\text{p}_{1}(\mathcal{O})\subset(\mathfrak{a}_1\oplus\mathfrak{n}_1)^*$
soit faiblement propre, il faut et il suffit que
$\text{p}_{1}(\mathcal{O})$ soit une $A_1N_1$-orbite coadjointe
ouverte (dans $(\mathfrak{a}_1\oplus\mathfrak{n}_1)^*$). De plus au
cas où $\text{p}_{1}$ est faiblement propre, elle est aussi propre.
   \end{theo}
\noindent\begin{demo} D'abord, d'après le théorème 9.1 de l'appendice,
$\text{p}_{1}(\mathcal{O})$ contient une $A_1N_1$-orbite coadjointe
ouverte de $(\mathfrak{a}_1\oplus\mathfrak{n}_1)^*$. De plus, un
groupe exponentiel n'ayant pas de sous-groupe fini non trivial,
$A_1N_1$ est diff\'{e}omorphe \`{a} toutes ses orbites coadjointes
ouvertes. Maintenant si $\text{p}_{1}(\mathcal{O})\triangleq \Omega$
est une orbite coadjointe ouverte, fixons $f_1\in \mathcal{O}\subset
\mathfrak{g}_1^* $ et $h_1 \in \Omega \subset
(\mathfrak{a}_1\oplus\mathfrak{n}_1)^*$. Alors $$A_1N_1\longrightarrow
A_1N_1.h_1$$ $$a\longmapsto a.h_1$$ est un
diff\'{e}omorphisme. Supposons que $L\subseteq \Omega$ est un compact,
et $x=b_1k_1f_1\in {p}_{1}^{-1}(L)$ o\`{u} $b_1 \in A_1N_1,\ k_1\in
K_1 $. Donc $b_1.\text{p}_{1}(k_1f_1)\in L$, et il est clair que
$\text{p}_{1}(K_1f_1)$ est un compact dans $\Omega$. Donc selon ce qui
pr\'{e}c\`{e}de, $\text{p}_{1}(K_1f_1)$ peut s'\'{e}crire
$\text{p}_{1}(K_1f_1)=\Theta_1.h_1$, où $\Theta_1$ est un compact dans
$A_1N_1$.  De m\^{e}me, $L=\Theta_2.h_1$ avec $\Theta_2$ un compact
dans $A_1N_1$. On en donc d\'{e}duit que $b_1\Theta_1\subset
\Theta_2$. Donc $b_1 \in \Theta_2.\Theta_1^{-1}$, ainsi, on d\'{e}duit
que $\text{p}_{1}^{-1}(L)\subseteq \Theta_2.\Theta_1^{-1}.K_1f_1 $ qui
est un compact dans l'orbite $\mathcal{O}$. Donc
$\text{p}_{1}^{-1}(L)$ est aussi compact, puisque
$\text{p}_{1}^{-1}(L)$ est fermé dans $\mathcal{O}$. Donc
$\text{p}_{1}$ est propre, donc faiblement propre.

Maintenant supposons que $\text{p}_{1}(\mathcal{O})$ n'est pas une
$A_1N_1$-orbite coadjointe ouverte de
$(\mathfrak{a}_1\oplus\mathfrak{n}_1)^*$, et $\Omega=A_1N_1.h_0 $ est
une orbite ouverte contenue dans $\text{p}_{1}(\mathcal{O})$ ici,
$h_0=\text{p}_{1}(g)$ avec $g\in \mathcal{O}$. Il est
clair que $\text{p}_{1}(\mathcal{O})=A_1N_1.\text{p}_{1}(K_1.g)$. Donc
on en d\'{e}duit que $\text{p}_{1}(K_1.g)\nsubseteq \Omega$. D'autre
part, on a $\overline{\Omega}\bigcap
\text{p}_{1}(\mathcal{O})\supsetneqq \Omega$, o\`{u}
$\overline{\Omega}$ est l'adh\'{e}rence de $\Omega$ dans
$(\mathfrak{a}_1\oplus\mathfrak{n}_1)^*$. En effet, sinon
$\overline{\Omega}\bigcap\text{p}_{1}(\mathcal{O})=\Omega$, donc
$\Omega$ est ouvert et ferm\'{e} dans $\text{p}_{1}(\mathcal{O})$ qui
est connexe (car $G_1$ est suppos\'{e} connexe), donc
$\Omega=\text{p}_{1}(\mathcal{O})$. Ceci est une contradiction avec
l'hypoth\`{e}se \og $\text{p}_{1}(\mathcal{O})$ n'est pas une orbite
ouverte \fg. Donc on en d\'{e}duit qu'il existe une suite $\{k_n\}$
dans $K_1$ telle que $\text{p}_{1}(k_n.\lambda)\in \Omega$ et
$\text{p}_{1}(k_n.g)\longrightarrow h\notin \Omega$.  Il est clair que
$\text{p}_{1}(k_n.g)$ s'\'{e}crit d'une fa\c{c}on unique
$\text{p}_{1}(k_n.g)=l_n.h_0$, avec $l_n\in A_1N_1$. Donc on en
d\'{e}duit que la suite $\{l_n\}$ n'est pas born\'{e}e dans $A_1N_1$
(c'est-\`{a}-dire $\{l_n\}$ n'est pas contenue dans un compact de
$A_1N_1$). Donc $\{(l_n)^{-1}\}$ n'est pas born\'{e}e dans
$A_1N_1$. D'autre part, il est clair que $(l_n)^{-1}k_n.g\in
\text{p}_{1}^{-1}(\{h_0\})$. Or il est bien connu que $\mathcal{O}$
est diff\'{e}omorphe \`{a} $A_1N_1\times K_1/T_1$, o\`{u} $T_1$ est le
stabilisateur $G_1(g)$ qui est d'après l'hypothèse, un tore de
$K_1$. Donc on en d\'{e}duit que $\text{p}_{1}^{-1}(\{h_0\})$ n'est
pas compact, donc $\text{p}_{1}$ n'est pas faiblement propre.  Le
th\'{e}or\`{e}me est donc bien démontré.

\end{demo}\\

Maintenant, on revient à $G=SU(2,1)$.  Il est clair que les
conditions du théorème précédent sont vérifiées pour $G$ et
$f_0\in\mathfrak{g}^*$ fortement régulière. Donc selon ce que l'on a
obtenu, on a le théorème suivant:

\begin{theo}\label{propreté-G}
Soit $f_0\in \mathfrak{t}^*\subset\mathfrak{g}^*$ une forme linéaire
fortement régulière. On pose $\mathcal{O}_{f_0}=G.f_0$. Alors la
projection $\text{p}_1: \mathcal{O}_{f_0}\longrightarrow
\text{p}_1(\mathcal{O}_{f_0})\subset \mathfrak{b}_1^* $ est propre ou
faiblement propre, si et seulement si $f_0$ est dans le cône
holomorphe.
\end{theo}

Gardons les notations du théorème \ref{propreté-G1} et supposons de
plus que $G_1/K_1$ soit un espace hermitien symétrique. Alors il est
bien connu que $A_1N_1$ admet des orbites coadjointes ouvertes. Dans
ce cas Rossi et Vergne ([29]) ont demontr\'{e} le r\'{e}sultat
suivant: si $\pi$ est une s\'{e}rie discr\`{e}te holomorphe ou
anti-holomorphe de $G_1$, alors $\pi$ est $A_1N_1$-admissible (dans le
chapitre suivant, on va d\'{e}tailler plus sur les r\'{e}sultats de
Rossi-Vergne). D'autre part Rosenberg et Vergne ([28]) ont
demontr\'{e} que si $\pi$ est une s\'{e}rie discr\`{e}te ni holomorphe
ni anti-holomorphe (de $G_1$), alors $\pi$ n'est pas
$A_1N_1$-admissible. La mise en perspective de ces résultats, du
th\'{e}or\`{e}me \ref{propreté-G1} et la conjecture de Duflo
nous amène à poser une autre conjecture:\\

\textbf{Conjecture }: En gardant les notations du théorème 5.3, on
suppose que $G_1/K_1$ est un espace hermitien symétrique et $\pi$ est
une s\'{e}rie discr\`{e}te de $G_1$ avec $G_1$-orbite coadjointe
associée $\mathcal{O}_{\pi}\subset \mathfrak{g}_{1}^*$. Alors nous
conjecturons que \textbf{ $\pi$ est holomorphe ou anti-holomorphe, si
  et seulement si $\text{p}_{1}(\mathcal{O}_{\pi})$ est une
  $A_1N_1$-orbite coadjointe ouverte dans
  $(\mathfrak{a}_1\oplus\mathfrak{n}_1)^*$}. Plus généralement si
$\mathcal{O}\subset\mathfrak{g}_{1}^*$ est l'orbite coadjointe d'un
élément $f\in \mathfrak{g}_{1}^*$ dont le stabilisateur est un tore
compact, alors les assertions suivantes sont équivalentes

(\textrm{i}) $\text{p}_{1}$ est propre

(\textrm{ii}) $\text{p}_{1}$ est faiblement propre

(\textrm{iii}) $\text{p}_{1}(\mathcal{O})$ est une $A_1N_1$-orbite coadjointe.

(\textrm{iv})  $\text{p}_{1}(\mathcal{O})$ est une $A_1N_1$-orbite coadjointe ouverte.

(\textrm{v}) $f$ est dans le cône holomorphe.

\noindent \textbf{Remarque}. Cette conjecture a été prouvée par l'auteur récemment, comme souligné dans la remarque de l'introduction.

\pagebreak
\subsection {La faible propreté (et la propreté) de la projection $\text{p}$}\label{propreté-p}

\begin{theo}

Soit $f_0\in \mathfrak{t}^*\subset\mathfrak{g}^*$ une forme linéaire fortement régulière et $\mathcal{O}_{f_0}=G.f_0$. Alors la projection $\text{p}: \mathcal{O}_{f_0}\longrightarrow \text{p}(\mathcal{O}_{f_0})\subset \mathfrak{b}^* $ est faiblement propre. De plus elle est propre si et seulement si $f_0$ est dans le cône holomorphe.

\end{theo}
\noindent\begin{demo} Rappelons que
$\mathfrak{b}=\mathfrak{m}\oplus\mathfrak{a}\oplus\mathfrak{n}$, et
$W$, $S$, $E_1$, $E_1'$ et $E_2$ est une base de $\mathfrak{b}$, où
$S$, $E_1$, $E_1'$ et $E_2$ sont comme dans la section précédente et

$$ W=\frac{i}{3}\left(
\begin{array}{ccc}
1& 0&0 \\
 0 & -2&0\\
0&0&1
\end{array}
\right)\in \mathfrak{m} .$$

Notons $W^*$, $S^*$, $E_1^*$, $E_1'^*$ et $E_2^*$ la base duale de
$\mathfrak{b}^*$. Soit $k\in K$. Comme $\mathfrak{p}$ est
$K$-invariant, $\text{p}(K.f_0)$ est contenu dans l'orthogonal de
$\mathfrak{p}\cap \mathfrak{b}=\mathfrak{a}=\R S$. On a
donc $$\text{p}(k.f_0)=\langle k.f_0,W\rangle W^*+\langle
k.f_0,E_2\rangle E_2^*+\langle k.f_0,E_1\rangle E_1^*+\langle
k.f_0,E_1'\rangle E_1'^*.$$ De sorte que $\text{p}(k.f_0)\in
\mathfrak{b}_{fr}^*\Leftrightarrow \langle k.f_0,E_2\rangle\neq
0$. Soit $K'=\left\{k\in K \vert \langle k.f_0,E_2\rangle\neq
0\right\}$, c'est un ouvert de $K$. Maintenant, rappelons que l'on a
obtenu un difféomorphisme :

$$\Psi: \R \times AN\times \{\pm1\}\longrightarrow \mathfrak{b}_{fr}^*$$

$$(r,x,\varepsilon)\longmapsto x.(rW^*+\varepsilon E_2^*),$$
et que, pour $f=wW^*+sS^*+xE_1^*+yE_1'^*+zE_2^*\in \mathfrak{b}_{fr}^*$,
$\Psi^{-1}(f)=(r(f),x(f),\varepsilon(f))$, avec
$r(f)=w+\frac{x^2+y^2}{2z}$ et $\varepsilon(f)=\frac{z}{\mid
  z\mid}$.

Soit donc $k\in K'$, on pose $r(k)=r(\text{p}(k.f_0))$,
$x(k)=x(\text{p}(k.f_0))$ et
$\varepsilon(k)=\varepsilon(\text{p}(k.f_0))$. Alors $$r(k)=\langle
k.f_0,W\rangle+\frac{\langle k.f_0,E_1\rangle^2+\langle
  k.f_0,E_1'\rangle^2}{2\langle k.f_0,E_2\rangle}.$$

Maintenant soit $H$, $F$, $V$ et $Z$ comme dans la section
précédente. 
Alors, on a
$$
k.f_0=\langle k.f_0,H\rangle H^*+\langle k.f_0,F\rangle F^*+\langle
k.f_0,V\rangle V^*+\langle f_0,Z\rangle Z^*\in\mathfrak{g}^*,$$ avec
$\langle k.f_0,H\rangle^2+\langle k.f_0,F\rangle^2+\langle
k.f_0,V\rangle^2=\langle f_0,H\rangle^2$ et ici pour $X\in
\mathfrak{g}$, $\langle H^*, X\rangle=-\frac{1}{2}\text{tr}(HX)$,
$\langle F^*, X\rangle=-\frac{1}{2}\text{tr}(FX)$, $\langle V^*,
X\rangle=-\frac{1}{2}\text{tr}(VX)$ et $\langle Z^*,
X\rangle=-\frac{1}{6}\text{tr}(ZX)$.

D'autre part, par des calculs directs, on a $\langle H^*,
W\rangle=\frac{1}{2}$, $\langle F^*, W\rangle=\langle V^*, W\rangle=0$
et $\langle Z^*, W\rangle=-\frac{1}{6}$, de sorte que
$\langle k.f_0,W\rangle=\frac{1}{2}\langle
k.f_0,H\rangle-\frac{1}{6}\langle f_0,Z\rangle$. Enfin, il suit des
relations \ref{image-base-p1} que   $\langle
k.f_0,E_2\rangle=\langle k.f_0,H\rangle+\langle f_0,Z\rangle$,
$\langle k.f_0,E_1\rangle=-\langle k.f_0,F\rangle$ et $\langle
k.f_0,E_1'\rangle=-\langle k.f_0,V\rangle$.

Donc on a $$r(k)=\langle k.f_0,W\rangle+\frac{\langle
  k.f_0,E_1\rangle^2+\langle k.f_0,E_1'\rangle^2}{2\langle
  k.f_0,E_2\rangle}$$

$$=\frac{1}{2}\langle k.f_0,H\rangle-\frac{1}{6}\langle
f_0,Z\rangle+\frac{\langle f_0,H\rangle^2-\langle
  k.f_0,H\rangle^2}{2(\langle k.f_0,H\rangle+\langle f_0,Z\rangle)}$$

$$=\frac{1}{3}\langle f_0,Z\rangle+\frac{\langle
  f_0,H\rangle^2-\langle f_0,Z\rangle^2}{2(\langle
  k.f_0,H\rangle+\langle f_0,Z\rangle)},$$

soit $$r(k)=\frac{1}{3}\langle f_0,Z\rangle+\frac{\langle f_0,H\rangle^2-\langle f_0,Z\rangle^2}{2\langle k.f_0,E_2\rangle}\ \ \ (**).$$

Maintenant, pour $\delta>0$, posons $K_{\delta}=\{k\in K \vert \mid\langle k.f_0,E_2\rangle\mid\geqslant\delta\}$, c'est un compact de $K$.

Soit $\Gamma\subset\mathfrak{b}_{fr}^*$ un compact. Alors il existe $u<v\in\R$ et $L\subset AN$ compact tels que $\Gamma\subset\Psi([u,v]\times L\times \{\pm1\})$. Soit $k\in K$, alors si $AN\text{p}(k.f_0)\cap\Gamma\neq\emptyset$, on a $k\in K'$ et $r(k)\in [u,v]$. Compte tenu de la formule $(**)$ donnant $r(k)$, cela montre qu'il existe $\delta>0$ tel que l'ensemble $K_{\Gamma}=\{k\in K \vert AN\text{p}(k.f_0)\cap\Gamma\neq\emptyset \}$ est contenu dans $K_{\delta}$. Soit $x\in AN$ et $k\in K$, alors on a $$\text{p}(xk.f_0)\in\Gamma\Leftrightarrow k\in K_{\Gamma} \ \ \text{et} \ \ xx(k)(r(k)W^*+\varepsilon(k)E_2^*)\in\Gamma$$

$$\Rightarrow k\in K_{\delta}, \ r(k)\in [u,v]\ \ \text{et}\ \ xx(k)\in L$$
$$\Rightarrow k\in K_{\delta},  \ \text{et}\ \ x\in L(x(K_{\delta}))^{-1}.$$

On voit donc que $\text{p}^{-1}(\Gamma)\subset L(x(K_{\delta}))^{-1}K_{\delta} f_0 $. Or $L$ et $K_{\delta}$ étant compacts et $\text{p}^{-1}(\Gamma)$ étant fermé, on en déduit que $\text{p}^{-1}(\Gamma)$ est compact, donc $\text{p}$ est faiblement propre.

Maintenant, si $\mid f_0(Z)\mid >\mid f_0(H)\mid$, on a $\text{p}(G.f_0)\subset \mathfrak{b}_{fr}^*$. Donc $\text{p}$ est propre. Si $\mid f_0(Z)\mid <\mid f_0(H)\mid$, alors on voit facilement que $K \backslash K'\neq \emptyset$ et si $k\in K \backslash K'$, on vérifie directement que $\exp\R E_2. (k.f_0)\subset \text{p}^{-1}(\text{p}(k.f))$ de sorte que $\text{p}$ n'est pas propre. Donc la démonstration est bien achevée.

\end{demo}

\pagebreak
\begin{center}\section{Représentations}\end{center}

Toutes les notations concernant $G=SU(2,1)$ des chapitres précédents
sont gardées pour ce chapitre.\\

Soit $\pi_{\lambda}$ la série discrète de $G$ avec le paramètre de
Harish-Chandra $\lambda\in i\mathfrak{t}^*$. Selon le chapitre 3, on
peut supposer que $\lambda$ correspond à un ensemble de racines
positives $\Delta_j^+$ défini dans le chapitre 3, où $j\in \{1, 2,
3\}$. D'apr\`{e}s la th\'{e}orie de Duflo (dans ce cas l\`{a} c'est
aussi la th\'{e}orie de Harish-Chandra pour les groupes de Lie
r\'{e}ductifs), $\pi_\lambda$ correspond \`{a} la forme lin\'{e}aire
$f_0=-i\lambda\in \mathfrak{t}^*\subset\mathfrak{g}^*$, où l'inclusion
"$\mathfrak{t}^*\subset\mathfrak{g}^*$" est relative à la
décomposition $\mathfrak{g}=\mathfrak{t}\oplus
[\mathfrak{t},\mathfrak{g}]$ (et aussi à la forme de Killing). Dans ce
cas, comme on a déjà évoqué dans le chapitre 2, il y a un seul élément
$\tau$ dans $X(f_0)$, et $T^{G}_{f_0,\tau}\cong \pi_{\lambda} $. Pour
simplifier, on désigne $T^{G}_{f_0,\tau}$ par $T^{G}_{f_0}$ et on
désigne l'orbite coadjointe associée $G.f_0$ par
$\mathcal{O}_{\pi_\lambda}$. De plus on note la complexifiée de $f_0$
encore $f_0$. Donc on a $f_0(E_{ij})=0$, et
$f_0\vert_{\mathfrak{t}}=-i\lambda$, o\`{u} les $E_{ij},\ i\neq j (1
\leqslant i,\,j\leqslant 3)$ sont des matrices \'{e}l\'{e}mentaires
dans $M_3(\C)$. \\

 Dans ce chapitre, on va déterminer la décomposition de
 $\pi_{\lambda}\vert_{B}$ et $\pi_{\lambda}\vert_{NA}$ en combinant
 plusieurs travaux et méthodes différents. Puis on va confirmer la
 conjecture de Duflo pour $G=SU(2,1)$ en interprétant la décomposition
 dans le cadre de la méthode des orbites. Enfin, on va en retirer des
 conséquences sur des systèmes différentiels qui semblent
 difficilement s'obtenir directement. \\

Maintenant on énonce le plan de ce chapitre:

1) Description des représentations (irréductibles et unitaires) de
$B_1$ et de $B$ associées aux orbites (coadjointes) fortement
régulières et admissibles.

2) Interprétation des travaux de Rossi-Vergne et confirmation de la conjecture de Duflo pour les séries discrètes holomorphes et anti-holomorphes de $G$.

3) Construction des séries discrètes de $G$ par des formes harmoniques (sens généralisé de $L^2$-Cohomologie).

4) Application de 3) à la décomposition de $\pi_{\lambda}\vert_{B}$ (et de $\pi_{\lambda}\vert_{B_1}$).

5) Etude du comportement asymptotique des solutions des systèmes
différentiels liés à 4).

6) Interprétation et application des travaux de Fabec et de ceux de Kraljevic.

7)  Confirmation des assertions (i) et (ii) de la conjecture de Duflo pour $G=SU(2,1)$.

8) Conséquences sur les systèmes différentiels dans 5).\\

\subsection{Description des représentations (irréductibles et unitaires) de $B_1$ et de $B$  associées aux orbites (coadjointes) fortement régulières et admissibles}\label{description-représentations}

Dans le chapitre précédent, on a déjà vu qu'il y a deux orbites
coadjointes fortement régulières et admissibles dans
$\mathfrak{b}_1^*$: $\Omega^-$ et $\Omega^+$, notons $f_{\pm}=\pm
E_2^* \in \Omega^{\pm}\subset \mathfrak{b}_1^*$. Le stabilisateur de
$f_{\pm}$ dans $B_{1}$ étant trivial,  il y a un seul
élément $\tau_{\pm}$ dans $X(f_{\pm})$, donc une seule représentation
unitaire irréductible associée à $\Omega^{\pm}$ , on la note
$\mathrm{T}_{\pm}$. Puisque $\Omega^-$ et $\Omega^+$ sont ouvertes,
$\mathrm{T}_{-}$ et $\mathrm{T}_{+}$ sont des séries discrètes de
$B_1$ (en fait il n'y a que deux séries discrètes pour $B_1$).

Maintenant en suivant la construction de Duflo, on va donner une
réalisation concrète pour $\mathrm{T}_{-}$ et $\mathrm{T}_{+}$.

Posons $n_{\pm}=f_{\pm}\vert_{\mathfrak{n}}$. Alors, $\mathfrak{n}$
est le radical unipotent de $\mathfrak{b}_{1}$ et le stabilisateur de
$n_{\pm}$ dans $B_{1}$ est le centre $Z$ de $N$. On en déduit
que $$\mathrm{T}_{\pm}=\text{Ind} \! \! \! \! \! \begin{array}{c} B_1
  \\ \uparrow \\ N \end{array}\! \! \! \! \! \mathrm{T}_{n_{\pm}}$$ où
$\mathrm{T}_{n_{\pm}}$ est la représentation irréductible unitaire de
$N$ associée à $n_{\pm}$ par la méthode des orbites de Kirillov.

On vérifie que
$\mathfrak{l}_{\pm}=\C\left(E_1\pm iE_1' \right)\oplus \C
E_2\subset\mathfrak{n}_{\C}$ est une polarisation positive de
$f_{\pm}$, c'est aussi une polarisation positive pour
$n_{\pm}=f_{\pm}\vert_{\mathfrak{n}}$.

 Nous réalisons alors $\mathrm{T}_{n_{\pm}}$ comme une induite
 holomorphe : Soit $\varphi$ une fonction lisse sur $N$ qui
 v\'{e}rifie $$X\ast\varphi= -n_{\pm}(X)\varphi \ \ \text{ pour tout}
 \ X \in \mathfrak{l}_{\pm}\ \ \ \ \ \ \ \ \ \ \ \ \ (\Delta)$$ où,
 $X\ast\varphi$ est l'action de $X$ sur $\varphi$ en tant
 qu'op\'{e}rateur diff\'{e}rentiel invariant à gauche.(ici, on
 prolonge naturellement l'action \`{a} $\mathfrak{n}_{\C}$, plus
 pr\'{e}cis\'{e}ent, pour $X, Y\in \mathfrak{n}$, $X\ast
 \varphi(x)=\left.\frac{d}{dt}\varphi(x\exp(tX)\right\vert_{t=0}$, et
 $(X+iY)\ast \varphi=X\ast \varphi+iY\ast \varphi$).

Comme $E_2 \in \mathfrak{l}_{\pm} $, on peut donc d\'{e}duire
facilement que pour tout $X \in \R E_2$, $\varphi(x\exp(X))=e^{-i
  n_{\pm}(X)}\varphi(x)$. Donc, $\vert\varphi \vert^2$ est bien
d\'{e}finie sur $N/\exp(\R E_2)$. Notons
$\mathbb{H}(n_{\pm},\mathfrak{l}_{\pm},\ N)$ le compl\'{e}t\'{e}
hilbertien de l'espace pr\'{e}hilbertien des fonctions lisses sur
$B_1$ qui v\'{e}rifient $(\Delta)$ avec la norme $$\Vert \varphi
\Vert^2=\int_{N/\exp(\R E_2)}\vert \varphi \vert^2dx<+\infty,$$ où
$dx$ est la mesure invariante à gauche sur $N/\exp(\R E_2)$. On sait
qu'elle existe, puisque $\exp(\R E_2)$ est distingué dans $N$. Il est
clair que $N$ op\`{e}re dans
$\mathbb{H}(n_{\pm},\mathfrak{l}_{\pm},\ N)$ par translations à
gauche, et cette repr\'{e}sentation que l'on note
$\rho(\ f_{\pm},\ \mathfrak{l}_{\pm},\ N )$ est bien une
r\'{e}alisation de $T_{n_{\pm}}$. Similairement, on peut définir une
représentation $\rho(\ f_{\pm},\ \mathfrak{l}_{\pm},\ B_1 )$ de $B_1$
dans $\mathbb{H}(n_{\pm},\mathfrak{l}_{\pm},\ B_1)$ (il suffit de
remplacer $N$ par $B_1$). Puisque $B_1$ est le produit semidirect de
$A$ et $N$, on déduit que

$$\rho(\ f_{\pm},\ \mathfrak{l}_{\pm},\ B_1 )\cong\text{Ind} \! \! \!
\! \! \begin{array}{c} B_1 \\ \uparrow \\ N \end{array}\! \! \! \! \!
\mathrm{T}_{n_{\pm}}=\mathrm{T}_{\pm}.$$

Maintenant soit $V=\R E_1\oplus \R E_1'$. Alors
$\ell_{\pm}=\mathfrak{l}_{\pm} \cap V_{\C}$ est un lagrangien complexe
totalement positif de $V$ pour $B_{n_{\pm}}$, et on a un isomorphisme
(d'espaces vectoriels réels): $$V\cong \ell_{\pm}$$

$$ x \longmapsto v_x \ \ \ \ \ \ \ \text{où}\ v_x \ \text{est l'unique
  élément de } \ell_{\pm} \ \text{tel que}\ x=v_x+\overline{v_x}.$$
D'où la structure complexe et hermitienne sur $V$ : $\langle x,
x\rangle=-i B_{n_{\pm}}(v_x, \overline{v_x}), \ x\in V $.

Soient $Sp(V)$ le groupe symplectique de $V$ pour $B_{n_{\pm}}$, et
$Mp(V)$ le groupe métaplectique qui agissent dans $N$ par
automorphismes. Alors il existe une unique représentation
$\mathrm{S}_{n_{\pm}}$ (la représentation métaplectique) de $Mp(V)$
dans $\mathbb{H}(n_{\pm},\mathfrak{l}_{\pm},\ N)$ qui vérifie

$$S_{n_{\pm}}(x.y)=S_{n_{\pm}}(x)\mathrm{T}_{n_{\pm}}(y)S_{n_{\pm}}(x)^{-1},\ \ x\in
Mp(V),\ y\in N .$$

Soit $U(\ell_{\pm})=\{ x\in Sp(V) \vert x.\ell_{\pm}=\ell_{\pm}
\}\cong U(V)$ le groupe unitaire de $V$ (pour la structure hermitienne
définie plus haut), et $\widetilde{U}(\ell_{\pm})$ son image
réciproque dans $Mp(V)$. Alors on a

$$\widetilde{U}(\ell_{\pm})\cong \{(u,z)\in U(\ell_{\pm})\times C^{\times}\}\vert z^2= \text{det}_{\ell_{\pm}}u \}.$$

Soit $\chi_{\pm}$ le caractère de $\widetilde{U}(\ell_{\pm})$ tel que
$\chi_{\pm}(u,z)=z$. Alors on définit une représentation
$S_{n_{\pm}}'$ de $U(\ell_{\pm})$ dans
$\mathbb{H}(n_{\pm},\mathfrak{l}_{\pm},\ N)$
par $$S_{n_{\pm}}'(x)\varphi(z)=\varphi(x^{-1}z).$$ Il vient alors

$$S_{n_{\pm}}(x)=\chi_{\pm}(x)S_{n_{\pm}}'(x),\ \ x\in \widetilde{U}(\ell_{\pm}).$$

Cela dit l'action adjointe de $M$ dans $\mathfrak{n}$ stabilise $V$ et est triviale sur $\R E_2$. Elle permet d'identifier $M$ à $U(V)$. La représentation $\mathrm{T}_{n_{\pm}}$ s'étend donc à $MN$ comme la représentation $S_{n_{\pm}}' \mathrm{T}_{n_{\pm}}$.

Maintenant soit $\widetilde{\mathrm{T}_{\pm}}=\text{Ind} \! \! \! \! \! \begin{array}{c} MAN \\ \uparrow \\ MN \end{array}\! \! \! \! \! S_{n_{\pm}}'\mathrm{T}_{n_{\pm}}$. Alors il est clair que $\widetilde{\mathrm{T}_{\pm}}\vert_{B_1}\cong \text{Ind} \! \! \! \! \! \begin{array}{c} B_1 \\ \uparrow \\ N \end{array}\! \! \! \! \! \mathrm{T}_{n_{\pm}}= \mathrm{T}_{\pm}$ de sorte que $\widetilde{\mathrm{T}_{\pm}}$ prolonge $T_{\pm}$ à $B=MB_1$. D'autre part, on vérifie facilement que $$(\widetilde{\mathrm{T}_{\pm}}(x)\varphi)(y)=\varphi (x^{-1}yx),\ \  \varphi\in\mathbb{H}(n_{\pm},\mathfrak{l}_{\pm},\ B_1) \ \ x\in M,\ y\in B_1.$$

On a vu que les $B$-orbites fortement r\'{e}guli\`{e}res et admissibles dans $\mathfrak{b}^*$ sont les $\Omega_{r,-}$ et $\Omega_{r,+}$, avec  $r+\frac{1}{2} \in \Z/3$. \\

 Maintenant, on va appliquer la th\'{e}orie de Duflo au sous-groupe de Borel $B=MAN$ d'alg\`{e}bre de Lie $\mathfrak{b}=\mathfrak{w}\oplus \mathfrak{a}\oplus \mathfrak{n}$.

On a vu que pour $f_{r,\pm}=rW^* \pm E_2^*\in \Omega_{r,\pm}$, $B(f_{r}^{\pm})=\exp\R W=M=U(V)$ est un tore de dimension 1. Donc  $B(f_{r,\pm})^{\mathfrak{b}}$ s'identifie à $\widetilde{U}(\ell_{\pm})$ qui est un tore connexe . Ainsi on en d\'{e}duit que si $f_{r,\pm}$ est admissible (i.e. $r+\frac{1}{2} \in \Z/3$), $X(f)$ contient un seul \'{e}l\'{e}ment que l'on note $\tau_{r,\pm}$, et $T_{f_{r,\pm},\tau_{r,\pm}}^{B}$ est la seule représentation unitaire irréductible associée à $\Omega_{r,\pm}$. Pour simplifier, pour $m\in \Z$, on  note \textbf{$\mathrm{T}_{m,-} := T_{f_{r,-},\tau_{r,-}}^{B}$, la représentation unitaire irréductible associée à $\Omega_{r,-}$, avec $r=\frac{m}{3}- \frac{1}{2}$, et $\mathrm{T}_{m,+} \triangleq T_{f_{r,-},\tau_{r,+}}^{B}$, mais avec $r=\frac{m}{3}+ \frac{1}{2}$}. D'autre part, puisque $B(f_{r}^{\pm})=\exp\R W$ est compact, selon la théorie de Duflo, l'ensemble \textbf{$ \{ \mathrm{T}_{m,\pm}: \ m\in \Z \}$ est exactement l'ensemble des séries discrètes de $B$}. On voit donc que \textbf{la formule de Plancherel pour $B$ ou $B_1$ ne fait apparaître que des séries discrètes}.

Dans la suite de cette section, on va donner une description concrète pour toutes les $\mathrm{T}_{m,\pm}$. Cette description est inspirée de Rossi-Vergne ([29]), et on va voir qu'elle est très utile dans la section \ref{décomposition} pour faire intervenir la formule de Plancherel. \\

Pour $m\in \Z$, on note $\sigma_{m}$ le caractère de $M$ tel que
$\sigma_{m}(\exp{tW})=e^{i(\frac{mt}{3})},\ \ t\in \R.$

On a déjà défini $\widetilde{\mathrm{T}_{\pm}}$. D'autre part puisque
$B_1$ est distingué dans $B$, il est \'{e}vident que l'on peut
prolonger trivialement $\sigma_{m}$ en une repr\'{e}sentation unitaire
irr\'{e}ductible de $B$ (c'est-\`{a}-dire que l'action de $B_1$ est
triviale), on la note encore $\sigma_{m}$. Il est facile de voir que
les représentations $\sigma_{m}\otimes \widetilde{\mathrm{T}_{\pm}}$
sont irréductibles, $\forall m\in \Z$, on va démontrer qu'elles ne
sont rien d'autres que $\mathrm{T}_{m,\pm}$.

\begin{theo}
Pour tout $m\in \Z$, $\mathrm{T}_{m,\pm} \cong \sigma_{m}\otimes
  \widetilde{\mathrm{T}_{\pm}}$.
\end{theo}

\noindent\begin{demo} Pour $m\in \Z$, posons $r_{m,
  \pm}=\frac{m}{3}\pm\frac{1}{2}$, et $f_{m, \pm}=r_{m, \pm}W^* \pm
E_2^*\in \Omega_{r,\pm}$. D'abord, il est clair que le plus grand
id\'{e}al nilpotent de $\mathfrak{b}$ est $\mathfrak{n}$, et
$n_{\pm}=f_{m, \pm}\vert_{\mathfrak{n}}$, où $n_{\pm}\in
{\mathfrak{n}}^*$ est ce que l'on a défini pour $B_1$. Posons aussi
$f_{m, \pm}^{1}=f_{m, \pm}\vert_{\mathfrak{b}_1}$ qui est en fait
indépendant de $m$. On a vu que $\mathfrak{l}_{\pm}=\C\left(E_1\pm
iE_1' \right)\oplus \C E_2\subset\mathfrak{n}_{\C}$ est une
polarisation positive pour $n_{m, \pm}$ et $f_{m, \pm}^{1}$. De plus
$\mathfrak{l}_{\pm}'=\mathfrak{m}\oplus \mathfrak{l}_{\pm}$ en est une
pour $f_{m, \pm}$. D'autre part par des calculs directs, on a $B(f_{m,
  \pm})=B(f_{m, \pm}^{1})=M$ et $B(n_{m, \pm})=M\exp\R E_2$. Donc
$B(n_{m, \pm})N=MN$, et on en déduit que les revêtements
métaplectiques $B(f_{m, \pm})^{\mathfrak{b}}$ et $B(f_{m,
  \pm}^{1})^{\mathfrak{b}_1}$ et l'image réciproque de $M$ dans
$B(n_{m, \pm})^{\mathfrak{n}}$ sont canoniquement isomorphes
à $$\widetilde{M}_{\pm}=\left\{\left.(x,z)\in \exp \R W\times
\C^{\times} \right\vert x=\exp tW,\ z^2=\text{det}(\exp
tW\vert_{\mathfrak{l}_{\pm}}) \right\}.$$ Comme pour $B_1$, on définit
le caractère $\chi_{\pm}$ pour $\widetilde{M}_{\pm}$ tel que
$\chi_{\pm}(\exp tW,z)=z$. On vérifie directement que $\text{det}(\exp
tW\vert_{\mathfrak{l}_{\pm}})=e^{\mp it}$ de sorte que
$\tau_{\pm}=\sigma_{m}\chi_{\pm}^{-1}$ est le seul élément dans
$X(f_{m, \pm})$, où $\sigma_{m}$ est le caractère unitaire de $M$ que
l'on a défini plus haut.

Alors, d'après les résultats de la section \ref{méthode-orbites}, on
a $$\mathrm{T}_{m,\pm}=\text{Ind}_{MN }^B(\tau_{\pm}\otimes
S_{n_{\pm}}T_{n_{\pm}})$$

$$=\text{Ind}_{MN }^B(\sigma_{m}\otimes S_{n_{\pm}}'\mathrm{T}_{n_{\pm}})$$

$$=\sigma_{m}\otimes \text{Ind}_{MN }^B S_{n_{\pm}}'\mathrm{T}_{n_{\pm}}$$

$$=\sigma_{m}\otimes \widetilde{T_{\pm}},$$ d'où le résultat.
\end{demo}
\begin{coro}
Soit $m\in\mathbb{Z}$. Alors le caractère central de la représentation
$T_{m,\pm}$ est la restriction du caractère $\sigma_{m}$ au centre
$Z_{G}$ de $B$.
\end{coro}
\begin{demo}
L'action de $Z_{G}$ par automorphismes dans $N$ étant triviale, le
caractère central de la représentation $\widetilde{T_{\pm}}$ est
trivial. D'où le corollaire.
\end{demo}

\subsection{Interprétation des travaux de Rossi-Vergne et confirmation de la conjecture de Duflo pour les séries discrètes holomorphes et anti-holomorphes de $G$ }\label{Rossi-Vergne}

Dans cette section, on garde les notations de la section
\ref{description-représentations}. On suppose que $\pi_{\lambda}$ est
une série discrète holomorphe (pour les séries discrètes
anti-holomorphes, on a exactement les mêmes résultats), i.e. $\lambda$
correspond à $\Delta_1^+=\{\alpha_{12},\ \alpha_{32},\ \alpha_{31}\}$,
ceci équivaut à $\lambda(H_{12}) \in \N^+$ et $\lambda(H_{31}) \in
\N^+$. Comme on a déjà expliqué, $\pi_\lambda$ est associée à
$f_0=-i\lambda\in \mathfrak{t}^*\subset \mathfrak{g}^*$ (au sens de
Duflo). Puisque $H=iH_{12}$ et $Z=i(-H_{12}-2H_{31})$, on a
$f_0(H)=-i\lambda(i H_{12})=\lambda(H_{12})$ et $f_0(Z)=-i\lambda
(i(-H_{12}-2H_{31}))=-\lambda(H_{12}+2H_{31})$ (donc surtout
$f_0(Z)+f_0(H)<0$). Donc on en déduit facilement que $\pi_{\lambda}$
est une série discrète holomorphe équivaut à dire que $f_0(H)\in \N^+$
et $f_0(Z)+f_0(H)$ est un entier pair strictement négatif. Dans la
suite de cette section on va exprimer les résultats en fonction de
$f_0(H)$ et $f_0(Z)$.\\

Il est clair que le param\`{e}tre de Blattner de $\pi_{\lambda}$ est
$\Lambda=\lambda+\delta_{G}-2\delta_K=\lambda+\alpha_{31}.$

Notons $U_\Lambda$ la repr\'{e}sentation unitaire irr\'{e}ductible de
$K$, dont le plus haut poids est $\Lambda$ (par rapport \`{a}
$\Delta_{1,K}^+=\{\alpha_{12}\}$). En appliquant la formule de Weyl,
la dimension de $U_\Lambda$ $$d_\Lambda=\left.\prod_{\alpha \in
  \Delta_{1,K}^+}\langle
\Lambda+\delta_K,\ \alpha\rangle\right/\prod_{\alpha \in
  \Delta_{1,K}^+}\langle
\delta_K,\ \alpha\rangle$$ $$=\left.\langle\Lambda+\delta_K,\ \alpha_{12}\rangle\right/\langle
\delta_K,\ \alpha_{12}\rangle$$ $$=\left.\langle\Lambda+\delta_K,\ H_{12}\rangle\right/\langle
\delta_K,\ H_{12}\rangle$$ $$=(\left.f_0(H) -1+1)\right/ 1=f_0(H).$$

Dans [29], Rossi et Vergne on donné une description de la restriction
d'une série discrète holomorphe d'un groupe de Lie simple à un
sous-groupe de Borel. Nous allons appliquer leur résultat au cas qui
nous intéresse et ensuite l'interpréter dans le cadre de la méthode
des orbites.

Rappelons que  les caractères de $M$ sont les  $\sigma_{m}$, $m\in\Z$
définis par
 $$\sigma_{m}(\exp{tW})=e^{i(\frac{mt}{3})},\ \ t\in \R.$$
On peut donc écrire
$$U_\Lambda\vert_{M}\cong\overset{f_0(H)-1}{\underset{m=0}{\bigoplus}}\sigma_{\Lambda_{m}},\ \text{où}
\ \Lambda_{m}\in \Z.$$ Alors d'après [29], on a :

\begin{theo}  La restriction de la série discrète holomorphe $\pi_{\lambda}$ \`{a} $B=MAN$ peut se d\'{e}composer de la mani\`{e}re suivante: $$\left.{\pi_\lambda}\right\vert_{B}\cong\overset{f_0(H)-1}{\underset{m=0}{\bigoplus}}\sigma_{\Lambda_{m}}\otimes\widetilde{\mathrm{T}_{-}}.$$

Donc on déduit que

$$\left.{\pi_\lambda}\right\vert_{B_1}\cong f_0(H).\mathrm{T}_{-}$$

\end{theo}

Rappelons qu'ici $\mathrm{T}_{-}$ est la série discrète de $B_{1}$
associée à l'orbite $\Omega^-$ et $\widetilde{\mathrm{T}_{-}}$ estson
prolongement à $B$ que l'on a défini dans la section
\ref{description-représentations}.

Puisque dans la section \ref{description-représentations}, on a montré
que $\mathrm{T}_{m,-} \cong \sigma_{m}\otimes
\widetilde{\mathrm{T}_{-}}$ (rappelons que $\mathrm{T}_{m,-}$ est la
série discrète de $B$ associée à $\Omega_{r,-}$, avec $r=\frac{m}{3}-
\frac{1}{2}$), on obtient que
$\left.{\pi_\lambda}\right\vert_{B}\cong\overset{f_0(H)-1}{\underset{m=0}{\bigoplus}}\mathrm{T}_{\Lambda_{m},-}$.

Dans la suite, on va calculer explicitement les $\Lambda_{m}$.

 Comme $M=\exp(\R W)$ est un tore de dimension 1, il suffit de
 calculer les valeurs propres de $dU_\Lambda(W)$.  Or on a,
 $W=\frac{1}{2}H-\frac{1}{6}Z$. Comme $Z$ est central dans
 $\mathfrak{k}$, $dU_\Lambda(Z)$ est l'opérateur scalaire
 $\Lambda(Z)Id$. Par ailleurs,
 $\Lambda(Z)=(\lambda+\alpha_{31})(Z)=i(f_{0}(Z)-3)$. D'autre part, la
 restriction de $U_{\Lambda}$ à $K'$ est la représentation
 irréductible de dimension $f_{0}(H)$ de $K'=SU(2)$. On sait alors que
 les valeurs propres de $dU_\Lambda(H)$ sont les $f_{0}(H)-1-2m$, avec
 $m$ entier et $0\leq m\leq f_{0}(H)-1$. On voit donc que les valeurs
 propres de $dU_\Lambda(W)$ sont
 $\frac{i(f_0(H)-1-2m)}{2}-\frac{i(f_0(Z)-3)}{6}=i(\frac{3f_0(H)-f_0(Z)}{6}-m)$.
 On en déduit donc que $\Lambda_{m}=\frac{3f_0(H)-f_0(Z)}{2}-3m$. Donc
 on peut résumer :

\begin{prop} Soit $\pi_{\lambda}$ une série discrète holomorphe de $G=SU(2,1)$ avec $\lambda$ son paramètre de Harish-Chandra qui vérifie que $\lambda(H_{12})= \N^+$  et $\lambda(H_{31}) \in \N^+$. Soit $f_0=-i\lambda\in \mathfrak{t}^*\subset \mathfrak{g}^*$ la forme linéaire associée (au sens de Duflo). Alors on a

$$\left.{\pi_\lambda}\right\vert_{B}\cong\overset{f_0(H)-1}{\underset{m=0}{\bigoplus}}\mathrm{T}_{(\frac{3f_0(H)-f_0(Z)}{2}-3m),-}$$

 et
 $$\left.{\pi_\lambda}\right\vert_{B_1}\cong f_0(H).\mathrm{T}_{-}$$
  \end{prop}

  Pour les séries discrètes anti-holomorphes, on a un résultat complètement analogue.

Maintenant, on va démontrer la conjecture de Duflo pour $\pi_\lambda$ holomorphe ou anti-holomorphe.

Supposons toujours que $\pi_\lambda$ est holomorphe, et soit
$f_0=-i\lambda$ la forme linéaire associée (au sens de Duflo). Donc
l'orbite coadjointe associée $\mathcal{O}_{\lambda}$ de $\pi_\lambda$
est $G.f_0$. Rappelons que l'on a noté les projections $\text{p}:
\mathcal{O}_{\lambda}\longrightarrow
\text{p}(\mathcal{O}_{\lambda})\subset \mathfrak{b}^* $ et
$\text{p}_1: \mathcal{O}_{\lambda}\longrightarrow
\text{p}_1(\mathcal{O}_{\lambda})\subset \mathfrak{b}_1^* $. On a vu
que dans notre cas $f_0(Z)+f_0(H)<0$. Donc d'après la formule de la
section \ref{propreté-p1} qui donne $\text{p}_1(K.f_0)$, on a
$\text{p}_1(\mathcal{O}_{\lambda})=\Omega^-$. Donc d'après la
proposition précédente et le théorème 5.3, \textbf{les assertions (i)
  et (ii) de la conjecture de Duflo sont confirmées pour
  $(G, \pi_\lambda,  B_1)$}.

Soit maintenant $k\in K$, on a vu dans la section \ref{propreté-p} que $\text{p}(k.f_0)\in \mathfrak{b}_{fr}^*$ si et seulement si $\langle k.f_0,E_2\rangle\neq 0$. Or  $\langle k.f_0,E_2\rangle=\langle k.f_0,H\rangle+\langle f_0,Z\rangle$ avec $\vert \langle k.f_0,H\rangle\vert\leq \vert\langle f_0,H\rangle\vert$, d'après ce qui précède, $\langle k.f_0,E_2\rangle< 0$. Donc on a $\text{p}(k.f_0)\subset \mathfrak{b}_{fr}^*$. Donc $\text{p}(G.f_0)\subset \mathfrak{b}_{fr}^*$. Plus précisément, d'après la formule $(**)$ de la section \ref{propreté-p} et ce qui la précède, on déduit que

$$\text{p}(\mathcal{O}_{\pi_\lambda})=\bigcup_{k\in K} B.(r(k)W^*-E_2^*)=\bigcup_{k\in K}\Omega_{r(k),-},$$

\noindent o\`{u} $r(k)=\frac{1}{3}\langle f_0,Z\rangle+\frac{\langle f_0,H\rangle^2-\langle f_0,Z\rangle^2}{2(\langle k.f_0,H\rangle+\langle f_0,Z\rangle)}$.\\

On a déjà vu que $\left\{ \langle k.f_0,H\rangle\vert k\in K \right\}=[-\vert\langle f_0,H\rangle\vert, \vert\langle f_0,H\rangle\vert]$, on en déduit donc facilement  que $\{r(k)\vert k\in K \}=[\frac{-(3f_0(H)+f_0(Z))}{6},\frac{3f_0(H)-f_0(Z)}{6}]$. Donc on a $$\text{p}(\mathcal{O}_{\pi_\lambda})=\bigcup_{\frac{-(3f_0(H)+f_0(Z))}{6} \leq r \leq \frac{3f_0(H)-f_0(Z)}{6}} \Omega_{r,-}.$$

Remarquons que  $-(3f_0(H)+f_0(Z))$ et $3f_0(H)-f_0(Z)$ sont des entiers pairs, on en déduit qu'il y a $3f_0(H)$ $B$-orbites coadjointes fortement régulières et admissibles dans $\mathfrak{b}^*$ qui sont contenues dans  $\text{p}(\mathcal{O}_{\pi_\lambda})$: Elles sont $$\Omega_{(\frac{3f_0(H)-f_0(Z)}{6}+ \frac{1-m}{3}-\frac{1}{2}),-},$$ où $0\leq m \leq 3f_0(H)-1$ et $m\in \Z$. On l'\'{e}crit comme une proposition. Voici\\

\begin{prop} Soit $\pi_{\lambda}$ une série discrète holomorphe de $G=SU(2,1)$ avec $\lambda$ son paramètre de Harish-Chandra qui vérifie que $\lambda(H_{12})\in \N^+$  et $\lambda(H_{31})\in \N^+$ et $f_0=-i\lambda\in \mathfrak{t}^*\subset \mathfrak{g}^*$ la forme linéaire associée (au sens de Duflo). Alors il y a $3f_0(H)$ $B$-orbites coadjointes fortement régulières et admissibles de $\mathfrak{b}^*$ qui sont contenues dans  $\text{p}(\mathcal{O}_{\pi_\lambda})$: Elles sont

$$\Omega_{(\frac{3f_0(H)-f_0(Z)}{6}+ \frac{1-m}{3}-\frac{1}{2}),-}$$

où $0\leq m \leq 3f_0(H)-1$ et $m\in \Z$. \end{prop}

D'après la section \ref{description-représentations}, pour $m\in \Z$, $\mathrm{T}_{m,-}$ est associée à l'orbite $\Omega_{(\frac{m}{3}-\frac{1}{2}),-}$, on peut donc déduire des propositions 6.4 et 6.5 et du théorème 5.5 que \textbf{les assertions
(i) et (ii) de la conjecture de Duflo sont confirmées pour  $(G, \pi_{\lambda},  B)$}.

Pour $\pi_\lambda$ anti-holomorphe, on peut confirmer de la même manière les assertions
(i) et (ii) de la conjecture de Duflo pour $(G, \pi_\lambda,  B_1)$ et $(G, \pi_\lambda,  B)$.

On peut donc faire un résumé:

\textbf{Pour $\pi_\lambda$ holomorphe ou anti-holomorphe, les assertions
(i) et (ii) de la conjecture de Duflo sont confirmées pour $(G, \pi_\lambda,  B_1)$ et $(G, \pi_\lambda,  B)$}. \\

\noindent \textbf{Remarque}. Il y a $3f_0(H)$ $B$-orbites coadjointes
fortement régulières et admissibles de $\mathfrak{b}^*$ qui sont
contenues dans $\text{p}(O_{\pi_\lambda})$ et chaque telle orbite
correspond à une (seule) représentation irréductible et unitaire de
$B$. cependant, il n'y a que $f_0(H)$ telles orbites qui figurent dans
la décomposition $\left.{\pi_\lambda}\right\vert_{B}$. On va voir un
phénomène analogue pour $\pi_\lambda$ ni holomorphe ni
anti-holomorphe.

\subsection{Construction des séries discrètes de $G$ par des formes harmoniques (sens généralisé de $L^2$-Cohomologie)}\label{L2C}

Dans cette section, on va rappeler des éléments sur la construction des séries discrètes d'un groupe réductif par la méthode des formes harmoniques. Dans la section suivante, on va appliquer cette méthode à la décomposition de $\pi_{\lambda}\vert_{B}$ (et de $\pi_{\lambda}\vert_{B_1}$), pour $G=SU(2,1)$.

Dans toute la suite de cette section, on ne suppose plus que $G=SU(2,1)$.\\

En 1982, Hersant a construit les "formes harmoniques et les modules de
cohomologie relative des alg\`{e}bres de Lie" g\'{e}n\'{e}ralisant
certaines constructions de Schmid ([30,31]), Narasimhan-Okamoto ([23])
et Penney. Dans la suite, on va d'abord faire un petit rappel de la
construction de Hersant, puis expliquer surtout comment elle
g\'{e}n\'{e}ralise la construction de Narasimhan-Okamoto. Dans la
section suivante, on va travailler dans le cadre de Hersant.\\

\textbf{Rappel de la th\'{e}orie de Hersant [11]} \\

Soient $G$ un groupe de Lie r\'{e}el, $K$ et $Z$ deux sous-groupes
ferm\'{e}s de $G$ qui v\'{e}rifient: (1) $Z$ est inclus dans le centre
de G. (2) $Z$ est inclus dans $K$, et $K/Z$ est compact.  Soient
$\mathfrak{g},\mathfrak{k},\mathfrak{z}$ les alg\`{e}bres de Lie
correspondantes, et
$\mathfrak{g}_{\C},\mathfrak{k}_{\C},\mathfrak{z}_{\C}$ leurs
complexifiées respectives.

Supposons qu'il existe une sous-alg\`{e}bre complexe $\mathfrak{e}$ de
$\mathfrak{g}_{\C}$ qui est stable par l'action adjointe de $K$, et
telle que $\mathfrak{g}_{\C}=\mathfrak{e}+\overline{\mathfrak{e}}$ et
$\mathfrak{e}\bigcap \overline{\mathfrak{e}}=\mathfrak{k}_{\C}$ (ici
"$^{-}$" signifie la conjugaison dans $\mathfrak{g}_{\C}$ par rapport
\`{a} la forme r\'{e}elle $\mathfrak{g}$)

Soit $V$ un $(\mathfrak{e},K)$-module qui est $K$-unitaire :
C'est-\`{a}-dire que $V$ est l'espace d'une repr\'{e}sentation
unitaire $\tau$ de $K$ et d'une repr\'{e}sentation not\'{e}e $\tau$
aussi, de $\mathfrak{e}$ vérifiant $$\forall k\in K, \forall
x\in \mathfrak{e}, \tau(\text{Ad}k.x)=\tau(k)\tau(x)\tau(k^{-1}),
\text{et} \ \forall x\in\mathfrak{k},
\tau(x)=\left.\frac{d}{dt}\right|_{t=0}\tau(\exp tx).$$ Supposons de
plus que $Z$ agit scalairement dans $V$ suivant le caract\`{e}re
$\chi$.

Maintenant soit $\pi$ une représentation unitaire de $G$ dans $\H$. On
note aussi $\pi$ la repr\'{e}sentation de $\mathfrak{g}_{\C}$ (et de
son alg\`{e}bre enveloppante universelle) dans $\H^\infty$ (ici
$\H^\infty$ est le sous-espace des vecteurs $C^{\infty}$ de $\pi$ dans
$\H$).

Consid\'{e}rons $\Lambda(\mathfrak{e}^*)\otimes
V^*\otimes\H^\infty=\oplus_{p}\Lambda^{p}(\mathfrak{e}^*)\otimes
V^*\otimes\H^\infty$ le $\mathfrak{e}$-complexe standard à valeurs
dans $V^*\otimes\H^\infty$ (pour $\tau^*\otimes \pi$). Notons
$\partial$ le cobord standard de Hoschild-Serre pour
$\Lambda(\mathfrak{e})\otimes V^*\otimes\H^\infty$, et
$\Theta=\text{Ad}^{*}\otimes\tau^*\otimes\pi$ la
repr\'{e}sentation de $K$ dans $\Lambda(\mathfrak{e}^*)\otimes
V^*\otimes\H^\infty$ d\'{e}finie par:
$(\Theta(k)\omega)(X_1,....,X_p)=(\tau^*(k)\otimes
\pi(k)).\omega(\text{Ad}k^{-1}X_1,...,\text{Ad}k^{-1}X_p)$ où $k\in K$
et $\omega\in \Lambda^{p}(\mathfrak{e}^*)\otimes V^*\otimes\H^\infty$
est une $p$-forme.  On peut v\'{e}rifier facilement que les
op\'{e}rateurs $\Theta(k)$ commutent au cobord $\partial$.

Consid\'{e}rons le sous-espace de $\Lambda(\mathfrak{e}^*)\otimes
V^*\otimes\H^\infty$ constitué des formes $\omega$ telles que
$i(X)\omega=0$ pour  tout $X\in\mathfrak{k}^{\C}$, o\`{u}
$(i(X)\omega)(Y_1,....,Y_{p-1})=\omega(X,Y_1,....,Y_{p-1})$ pour
$\omega$ une $p$-forme. Cet espace s'identifie naturellement
à $$\Lambda((\mathfrak{e}/\mathfrak{k}_{\C})^*)\otimes
V^*\otimes\H^\infty$$ qui est stable par la repr\'{e}sentation
$\Theta$. Notons $$[\Lambda((\mathfrak{e}/\mathfrak{k}_{\C})^*)\otimes
  V^*\otimes\H^\infty]^K=\left\{\omega\in
\Lambda((\mathfrak{e}/\mathfrak{k}^{\C})^*)\otimes
V^*\otimes\H^\infty| \ \Theta(k)\omega=\omega, \forall k\in
K\right\}.$$ On peut v\'{e}rifier que
$[\Lambda((\mathfrak{e}/\mathfrak{k}_{\C})^*)\otimes
  V^*\otimes\H^\infty]^K$ est stable par l'op\'{e}rateur $\partial$,
donc d\'{e}finit par restriction de $\partial$, un sous-complexe que
l'on note $$[\Lambda((\mathfrak{e}/\mathfrak{k}_{\C})^*)\otimes
  V^*\otimes\H^\infty]^K\stackrel{\partial_K}{\longrightarrow}[\Lambda((\mathfrak{e}/\mathfrak{k}_{\C})^*)\otimes
  V^*\otimes\H^\infty]^K.$$

Puisque $K/Z$ est supposé compact, $\text{Ad}K$ est un groupe
compact. Donc sur $\Lambda(\mathfrak{e}^*)$ (ainsi que sur
$\Lambda((\mathfrak{e}/\mathfrak{k}_{\C})^*)$ ), il y a un produit
scalaire $\text{Ad}K-$invariant. Donc l'espace
$\Lambda(\mathfrak{e}^*)\otimes V^*\otimes\H$ a une structure
naturelle de produit tensoriel hilbertien, pour laquelle la
repr\'{e}sentation $\Theta$ de $K$ est unitaire.

Hersant a d\'{e}montr\'{e} les assertions qui suivent. Le sous-espace
que l'on note $[\Lambda((\mathfrak{e}/\mathfrak{k}_{\C})^*)\otimes
  V^*\otimes\H]^K$ des vecteurs invariants par $\Theta(K)$ est
l'adh\'{e}rence de l'espace
$[\Lambda((\mathfrak{e}/\mathfrak{k}_{\C})^*)\otimes
  V^*\otimes\H^\infty]^K$ dans $\Lambda(\mathfrak{h}^*)\otimes
V^*\otimes\H$. L'op\'{e}rateur
$[\Lambda((\mathfrak{e}/\mathfrak{k}_{\C})^*)\otimes
  V^*\otimes\H^\infty]^K\stackrel{\partial_K}{\longrightarrow}[\Lambda((\mathfrak{e}/\mathfrak{k}_{\C})^*)\otimes
  V^*\otimes\H^\infty]^K$ a un adjoint formel $\delta_K$ sur le
m\^{e}me espace, et si on pose
$\square_{K}=\partial_{K}\delta_{K}+\delta_{K}\partial_{K}$,
alors $$[\text{Im}(\partial_K+\delta_K)]^\bot=\ker(\partial_K)\cap\ker(\delta_K)=\ker\square_{K},$$
\noindent o\`{u} $[\text{Im}(\partial_K+\delta_K)]^\bot$ est
l'orthogonal de $\text{Im}(\partial_K+\delta_K)$ dans
$[\Lambda((\mathfrak{e}/\mathfrak{k}_{\C})^*)\otimes
  V^*\otimes\H]^K$. En particulier,
 $\ker(\partial_K)\cap\ker(\delta_K)$ est fermé dans
$[\Lambda((\mathfrak{e}/\mathfrak{k}_{\C})^*)\otimes
  V^*\otimes\H]^K$.

Pour simplifier, on appelle le sous-complexe construit ci-dessus le
\textbf{"complexe de Hersant" et on le note $C^*(\mathfrak{e},V,K,
  G,\pi )$}.  Puisque dans toute la suite, on prend toujours "$Z$"
trivial pour tous les groupes concernés, on n'indique pas "$Z$" dans
la notation du complexe. De plus dans toute la suite, si le "$K$"
concerné est aussi trivial, on remplacera $C^*(\mathfrak{e},V,K, G,\pi
)$ par $C^*(\mathfrak{e},V, G,\pi )$.\\

Maintenant soit $R$ la repr\'{e}sentation de $G$ dans
$\H=L^2(G/Z,\chi)$ par translations \`{a} droite, o\`{u}
$L^2(G/Z,\chi)$ est l'espace des fonctions $f$ de carr\'{e}
int\'{e}grable modulo $Z$ qui v\'{e}rifient $f(zg)=f(gz)=\chi(z)f(g),
\forall g\in G, \forall z\in Z$, on note aussi $R$ la
repr\'{e}sentation de $\mathfrak{g}_{\C}$ (et de son alg\`{e}bre
enveloppante universelle ) dans $\H^\infty$. Ici $\H^\infty$ est le
sous-espace des vecteurs lisses dans $\H$ et $L^2(G/Z,\chi)$ est
relatif à la mesure de Haar invariante à gauche. Donc la translation à
droite de $G$ dans $\H$ a le terme \og $\Delta^{1/2}$ \fg devant, si
$G$ n'est pas unimodulaire, o\`{u} $\Delta$ est la fonction module de
$G$.

Hersant a également d\'{e}montr\'{e} que le sous-espace
$\ker(\partial_K)\cap\ker(\delta_K)$ du complexe
$C^*(\mathfrak{e},V,K, G, R )$ est  invariant par la
repr\'{e}sentation unitaire $1\otimes 1\otimes l$ de $G$ dans
$\Lambda(\mathfrak{e}^*)\otimes V^*\otimes\H$ (o\`{u} $1\otimes
1\otimes l$ signifie que $G$ agit trivialement dans
$\Lambda(\mathfrak{e}^*)\otimes V^*$ et par translations
r\'{e}guli\`{e}re \`{a} gauche dans $\H$). On note la
sous-repr\'{e}sentation associ\'{e}e $\pi(\mathfrak{e},\chi)$. Comme
la représentation $1\otimes
1\otimes l$ laisse clairement invariant le sous-espace des $q$-formes,
on a $\pi(\mathfrak{e},\chi)=\bigoplus
\pi^q(\mathfrak{e},\chi)$, o\`{u} $\pi^q(\mathfrak{e},\chi)$ est la
sous-repr\'{e}sentation associ\'{e}e de $G$ dans
$(\ker(\partial_K)\cap\ker(\delta_K))\cap\Lambda^{q}(\mathfrak{e}^*)\otimes
V^*\otimes\H$.\\

Maintenant, soit $G$ un groupe de Lie semi-simple d'alg\`{e}bre de Lie
$\mathfrak{g}$, $K$ un sous-groupe compact maximal de $G$
d'alg\`{e}bre de Lie $\mathfrak{k}$. Soit
$\mathfrak{g}=\mathfrak{k}\oplus \mathfrak{p}$ la décomposition de
Cartan associ\'{e}e et soit
$\mathfrak{g}_{\C},\mathfrak{k}_{\C},\mathfrak{p}_{\C}$ les
complexifi\'{e}s respectifs de
$\mathfrak{g},\mathfrak{k},\mathfrak{p}$. On fixe une décomposition
d'Iwasawa $G=KAN$. Supposons que sur l'espace homog\`{e}ne $G/K$, il
existe une structure complexe $G$-invariante, alors il existe deux
sous-alg\`{e}bres complexes $\mathfrak{p}_+,\mathfrak{p}_-$ de
$\mathfrak{g}^{\C}$ telles
que $$\mathfrak{p}_+\oplus\mathfrak{p}_-=\mathfrak{p}_{\C},\ \mathfrak{p}_+=\overline{\mathfrak{p}_-},\ [\mathfrak{k}_{\C},\mathfrak{p}_+]\subseteq
\mathfrak{p}_+.$$ Donc $[\mathfrak{k}_{\C},\mathfrak{p}_-]\subseteq
\mathfrak{p}_-$. De plus l'espace des vecteurs tangents
anti-holomorphes est canoniquement identifi\'{e} \`{a}
$\mathfrak{p}_+$. Dans ce cas, $\mathfrak{k}$ contient une
sous-alg\`{e}bre de Cartan $\mathfrak{t}$ de $\mathfrak{g}$ (donc
l'ensemble des s\'{e}ries discr\`{e}tes de $G$ n'est pas vide). Notons
$\Delta$ l'ensemble des racines de $\mathfrak{g}_{\C}$ relatif à
$\mathfrak{t}_{\C}$, et $\Delta_K$ (resp. $\Delta_n$ )l'ensemble des
racines compactes (resp. non compactes) pour $\Delta$.  Alors il
existe un et un seul sous-ensemble $\Delta^{+}_n \subset\Delta_n$ tel
que $\mathfrak{p}_+=\bigoplus_{\alpha\in
  \Delta^+_n}\mathfrak{g}^{\alpha}$, o\`{u} $\mathfrak{g}^{\alpha}$
est l'espace radiciel associ\'{e} \`{a} $\alpha$, et il existe un
ensemble de racines positives $\Delta^{+}$ de $\Delta$, tel que $
\Delta^{+}_n\subset \Delta^{+}$. C'est-\`{a}-dire que $ \Delta^{+}_n$
est l'ensemble de racines positives non-compactes par rapport à
$\Delta^{+}$, on note $\Delta^{+}_K$ l'ensemble de racines positives
associé (donc $\Delta^{+}= \Delta^{+}_n\cup \Delta^{+}_K$) et
$\mathfrak{F}$ l'ensemble des formes alg\'{e}briquement
int\'{e}grables sur $\mathfrak{t}_{\C}$, i.e. l'ensemble des formes
lin\'{e}aires complexes $\lambda$ sur $\mathfrak{t}_{\C}$ telles ques
$\lambda(H_\alpha)$ soit entier pour tout $\alpha\in \Delta$, ici
$H_\alpha$ est la coracine de $\alpha$. Pour simplifier, on suppose
que le complexifié $G_{\C}$ de $G$ existe et est simplement connexe,
de sorte que toute forme alg\'{e}briquement int\'{e}grable est
analytiquement int\'{e}grable. Posons

$$\mathfrak{F}'=\left\{\lambda\in \mathfrak{F}\vert(\lambda+\rho)(H_\alpha)\neq 0\ \ \forall\alpha \in \Delta\right\}$$

$$\mathfrak{F}_0'=\left\{\lambda\in \mathfrak{F}'\vert(\lambda+\rho)(H_\alpha)> 0\ \ \forall\alpha \in \Delta^{+}_K\right\}$$

\noindent o\`{u} $\rho=\frac{1}{2}\sum_{\alpha\in \Delta^{+}}\alpha$
est la demi-somme des racines positives associ\'{e}es \`{a}
$\Delta^{+}$.

Soient $\Lambda\in \mathfrak{F}_0'$ et $\tau_{\Lambda}$ la
repr\'{e}sentation unitaire irr\'{e}ductible de $K$ avec le plus haut
poids $\Lambda$ par rapport \`{a} $\Delta_K$, et $V_{\Lambda}$
l'espace hilbertien associ\'{e}. Notons $\tau_{\Lambda}^*$ la
repr\'{e}sentation contragr\'{e}diente de $\tau_{\Lambda}$, et
$V_{\Lambda}^*$ l'espace hilbertien associ\'{e}. Soient
$Q_{\Lambda}=\left\{\alpha\in
\Delta^{+}_n\vert(\Lambda+\rho)(H_\alpha)>0\right\}$ et
$q_{\Lambda}=\text{card}(Q_{\Lambda})$.\\

Maintenant, on va énoncer un théorème de Narasimhan-Okamoto dans le
cadre de la théorie de Hersant. On reprend les notations dans "Rappel
de la th\'{e}orie de Hersant" et celles ci-dessus. Supposons que $Z$
est trivial ( $\chi$ aussi), $K$ est le sous-groupe compact maximal de
$G$, $\mathfrak{e}=\mathfrak{p}_{+}\oplus \mathfrak{k}_{\C}$
($(\mathfrak{e}/\mathfrak{k}_{\C})^*\cong \mathfrak{p}_{-}$), et
$V=V_{\Lambda}$, $\tau=\tau_{\Lambda}$, l'action de $\mathfrak{p}_{+}$
dans $V_{\Lambda}$ étant triviale. On note la représentation
$\pi^{q}(\mathfrak{e},1)\triangleq \pi^{q}(\mathfrak{e})$. Alors  le
 th\'{e}or\`{e}me suivant est dû à Narasimhan et Okamoto ([23]) :

\begin{theo}
Si $q\neq q_{\Lambda}$, alors $\pi^{q}(\mathfrak{e})$ est réduite à
zéro, et $\pi^{q_{\Lambda}}(\mathfrak{e})$ est une s\'{e}rie
discr\`{e}te de $G$. Plus pr\'{e}cis\'{e}ment,
$\pi^{q_{\Lambda}}(\mathfrak{e})$ est la repr\'{e}sentation
contragr\'{e}diente de la s\'{e}rie discr\`{e}te de $G$ dont le
param\`{e}tre de Harish-Chandra est $\Lambda+\rho$, donc le
param\`{e}tre de Harish-Chandra de $\pi^{q_{\Lambda}}(\mathfrak{e})$
est $-(\Lambda+\rho)$.
\end{theo}

Maintenant, on veut réaliser $\pi_{\lambda}$ comme une
$\pi^{q_{\Lambda}}(\mathfrak{e})$. Sans perte de généralité, on peut
supposer $\Delta^{+}_{\lambda}\cap \Delta_{K}=\Delta^{+}_{K}$, où
$\Delta^{+}_{\lambda}$ est l'ensemble de racines positives par rapport
à $\lambda$. Soient $W_{K}$ le groupe de Weyl pour $\Delta_{K}$ et
$w_{K}$ le seul élément dans $W_{K}$ tel que
$w_{K}(\Delta^{+}_{K})=-\Delta^{+}_{K}$. Alors il est clair que
$\Lambda=-w_{K}\lambda-\rho\in\mathfrak{F}_0'$ et
$-(\Lambda+\rho)=w_{K}.\lambda$, donc
$\pi^{q_{\Lambda}}(\mathfrak{e})=\pi_{w_{K}.\lambda}$. Or
$\pi_{w_{K}.\lambda}\cong\pi_{\lambda}$. Donc on peut réaliser
$\pi_{\lambda}$ par $\pi^{q_{\Lambda}}(\mathfrak{e})$ avec
$\Lambda=-w_{K}\lambda-\rho$.\\

\subsection{ Application de \ref{L2C} à la décomposition de $\pi_{\lambda}\vert_{B}$ (et de $\pi_{\lambda}\vert_{B_1}$)}\label{décomposition}

Maintenant on revient sur $G=SU(2,1)$ et on reprend toutes les
notations des chapitres précédents qui concernent $G$. On se donne
$\lambda\in i\mathfrak{t}^{*}$ correspondant à une série discète
$\pi_{\lambda}$ de $G$. On suppose, comme il est loisible que
$\Delta^{+}_{\lambda}$ est l'un des $\Delta^{+}_{j}$, $j=1,2,3$.

Dans notre situation, on peut prendre $\mathfrak{p}_{+}$ et
$\mathfrak{p}_{-}$ de la section précédente comme

 $$ \mathfrak{p}_{+}=\left\{\left.\left(
\begin{array}{cc}
0 & 0 \\
 Y  & 0 \\
\end{array}
\right)  \right\vert Y \in M_{1,2}(\C) \right\}$$

$$ \mathfrak{p}_{-}=\left\{\left.\left(
\begin{array}{cc}
0 & X \\
 0  & 0 \\
\end{array}
\right)  \right\vert X \in M_{2,1}(\C) \right\}$$

\noindent o\`{u} $\mathfrak{p}_\C=\left\{\left.\left(
\begin{array}{cc}
0 & X \\
 Y & 0 \\
\end{array}
\right) \right\vert X \in M_{2,1}(\C),\ Y \in M_{1,2}(\C)\right\}$ est
le complexifié de $\mathfrak{p}$.  Soit
$\mathfrak{e}=\mathfrak{p}_+\oplus\mathfrak{k}_{\C}$. On vérifie
facilement que l'ensemble de racines positives $\Delta^{+}=\Delta_1^+$
défini dans le chapitre 3 est tel que $\Delta_{K}^{+}=\alpha_{12}$ et
$\mathfrak{p}_{+}=\oplus_{\alpha\in\Delta_{n}^{+}}\mathfrak{g}^{\alpha}$. Par
suite, $\rho=\alpha_{32}$ et la symétrie $s_{12}$ par rapport à la
racine $\alpha_{12}$ est l'élément dans $W_K$ tel que
$s_{12}(\Delta_{K}^{+})=-\Delta_{K}^{+}$. Donc si
$\Delta_{\lambda}^{+}=\Delta_j^+$ ($j\in\{1,2,3\}$), d'après ce que
l'on a vu, $\pi_{\lambda}=\pi^{q_{\Lambda}}(\mathfrak{e})$, avec
$\Lambda=-s_{12}\lambda-\alpha_{32}$. Or
$\alpha_{32}=s_{12}(\alpha_{31})$, donc
$\Lambda=-s_{12}(\lambda+\alpha_{31})$ et on en déduit que
$V^*_{\Lambda}\cong V_{-\Lambda}\cong V_{\lambda+\alpha_{31}}$.

Cependant on a aussi
$\mathfrak{p}_+\oplus\mathfrak{k}_{\C}=\mathfrak{h}_{+}\oplus\mathfrak{k}_{\C}$,
où $\mathfrak{h}_{+}=\C(E_1'+iE_1)\oplus \C (S-iE_2/2)$. On en
déduit que $\mathfrak{e}_{\mathfrak{b}_1}=\mathfrak{h}_{+}$
(resp. $\mathfrak{e}_{\mathfrak{b}}=\mathfrak{h}_{+}\oplus
\mathfrak{m}_{\C}$) qui est considérée comme une sous-algèbre de
$(\mathfrak{b}_1)_{\C}$ (resp. $\mathfrak{b}_{\C}$ ), avec "$K, Z$"
triviaux, "$V=V_{\Lambda}$" dans lequel l'action de
$\mathfrak{e}_{\mathfrak{b}_1}$ provient par restriction de celle de
$\mathfrak{e}$ (resp. avec "$K=M$", "$Z$" trivial, "$V=V_{\Lambda}$"
dans lequel l'action de $\mathfrak{e}_{\mathfrak{b}}$ provient par
restriction de celle de $\mathfrak{e}$) vérifient les "hypothèses de
Hersant". Comme $G=KB_{1}=KB$, on voit alors que le complexe de
Hersant $C^*(\mathfrak{e},V_{\Lambda}, K, G )$ peut aussi
s'interpréter comme $C^*(\mathfrak{e}_{\mathfrak{b}_1},V_{\Lambda},
B_1 )$ ou $C^*(\mathfrak{e}_{\mathfrak{b}},V_{\Lambda}, M, B )$, et que
les espaces de cohomologies sont les mêmes

 $$H^{q_{\Lambda}}(\mathfrak{e},V_{\Lambda}, G)= H^{q_{\Lambda}}(\mathfrak{e}_{\mathfrak{b}},V_{\Lambda}, B)=H^{q_{\Lambda}}(\mathfrak{e}_{\mathfrak{b}_1},V_{\Lambda}, B_1),$$

\noindent chacun d'eux fournissant la restriction de la représentation
$\pi_{\lambda}=\pi^{q_{\Lambda}}(\mathfrak{e})$ au sous-groupe
correspondant. Donc $\pi_{\lambda}\vert_{B_1}$
(resp. $\pi_{\lambda}\vert_{B}$) est l'action de $B_1$ par
translations à gauche dans
$H^{q_{\Lambda}}(\mathfrak{e}_{\mathfrak{b}_1},V_{\Lambda}, B_1)$
(resp. $H^{q_{\Lambda}}(\mathfrak{e}_{\mathfrak{b}},V_{\Lambda}, B)$).

Notons $\xi$ la représentation de
$\mathfrak{p}_{+}\oplus\mathfrak{k}_{\C}$ dans $V_\Lambda^* \cong
V_{\lambda+\alpha_{31}}$ d\'{e}finie par
$\xi\vert_{\mathfrak{p}_{+}}=0$ et
$\xi\vert_{\mathfrak{k}_{\C}}=\tau^*_\Lambda=\tau_{\lambda+\alpha_{31}}$.
On  note $\xi\vert_{\mathfrak{h}_{+}}$ encore par $\xi$. On désigne
l'espace dans lequel $\mathrm{T}_{-}$ (resp. $\mathrm{T}_{+}$) agit
par $\H_{-}$ (resp. $\H_{+}$), et $\H_{-}^\infty$
(resp. $\H_{+}^\infty$) le sous-espace des vecteurs $C^{\infty}$ de
$\mathrm{T}_{-}$ (resp. $\mathrm{T}_{+}$). Pour $q_{\Lambda}=0,1,2$,
définissons

$$(\delta_\pm)_{q_{\Lambda}}:\Lambda^{q_{\Lambda}}(\mathfrak{h}_{+})^*\otimes V_{\lambda+\alpha_{31}}\otimes \H_{\pm}^\infty\longrightarrow
\Lambda^{q_{\Lambda}+1}(\mathfrak{h}_{+})^*\otimes V_{\lambda+\alpha_{31}}\otimes \H_{\pm}^\infty,$$

\noindent le cobord standard de Hochschild-Serre pour le
$\mathfrak{h}_{+}$-module $ V_{\lambda+\alpha_{31}}\otimes
\H_{\pm}^\infty$, où l'action de $\mathfrak{h}_{+}$ dans
$\H_{\pm}^\infty$ est induite par $\mathrm{T}_{\pm}$, mais celle de
$\mathfrak{h}_{+}$ dans $ V_{\lambda+\alpha_{31}}$ est $\xi+(\text{tr
  ad}_{\mathfrak{b}_{1}}).\text{id}/2$ (non pas simplement $\xi$). Soit
${(\delta_\pm)_{q_{\Lambda}}}^*$  l'adjoint formel de
$(\delta_\pm)_{q_{\Lambda}}$ tel qu'il est explicité dans
  [11]. Alors il est clair que le complexe ci-dessus n'est rien autre
  que le complexe de Hersant
  $C^*(\mathfrak{e}_{\mathfrak{b}_1},V_{\lambda+\alpha_{31}},
  B_1,\mathrm{T}_{\pm})$ mais pour lequel l'action de
  $\mathfrak{e}_{\mathfrak{b}_1}=\mathfrak{h}_{+}$ est $\xi+(\text{tr
    ad}_{\mathfrak{b}_{1}}).\text{id}/2$. Notons alors
  $\mathfrak{H}_{\pm}^{q_{\Lambda}}=\ker(\delta_\pm)_{q_{\Lambda}}\cap\ker
  {(\delta_\pm)_{q_{\Lambda}}}^*$ et
  $\mathfrak{H}_{\pm}=\bigoplus_{q_{\Lambda}=0}^{2}\mathfrak{H}_{\pm}^{q_{\Lambda}}$.

Maintenant, considérons le $(\mathfrak{e}_{\mathfrak{b}}, M)$-module
$V_{\lambda+\alpha_{31}}$, où l'action de $M$ dans
$V_{\lambda+\alpha_{31}}$ est induite par celle de $K$, mais l'action
de $\mathfrak{e}_{\mathfrak{b}}$ dans $V_{\lambda+\alpha_{31}}$ est
$\xi+(\text{tr ad}_{\mathfrak{b}}).\text{id}/2$ (rappelons que on a
déjà défini le $\mathfrak{e}$-module $\xi$ plus haut, ici on désigne
$\xi\vert_{\mathfrak{e}_{\mathfrak{b}}}$ encore par $\xi$). Donc à
partir de ce module, on peut construire le complexe de Hersant
$C^*(\mathfrak{e}_{\mathfrak{b}},V_{\lambda+\alpha_{31}}, M,
B,\mathrm{T}_{m,\pm})$, où $\mathrm{T}_{m,\pm}$ est la représentation
unitaire irréductible de $B$ que l'on a déjà définie. Puisque
$\mathrm{T}_{m,\pm}=\sigma_m\otimes\widetilde{\mathrm{T}_{\pm}}$ qui
prolonge $\mathrm{T}_{\pm}$, et
$(\mathfrak{e}_{\mathfrak{b}})/\mathfrak{m}_{\C}\cong
\mathfrak{h}_{+}$, on déduit que le
"$\ker(\partial_K)\cap\ker(\delta_K)$" (pour la notation, voir le
"Rappel de la th\'{e}orie de Hersant") pour
$C^*(\mathfrak{e}_{\mathfrak{b}},V_{\lambda+\alpha_{31}}, M,
B,\mathrm{T}_{m,\pm})$ n'est rien d'autre que
$(\mathfrak{H}_{\pm})_{\sigma_{-m}}$, où pour $m\in \Z$,
$(\mathfrak{H}_{\pm})_{\sigma_{m}}$ désigne le sous-espace de
$\mathfrak{H}_{\pm}$ constitué des vecteurs de poids $\sigma_{m}$ pour
la représentation
$\text{Ad}^{*}\otimes\tau_{\lambda+\alpha_{31}}\otimes
\widetilde{\mathrm{T}_{\pm}}$ de $M$. On note
$(\mathfrak{H}_{\pm})^{q_{\Lambda}}_{\sigma_{m}}=(\mathfrak{H}_{\pm})_{\sigma_{m}}\cap
\mathfrak{H}_{\pm}^{q_{\Lambda}}$.

Rappelons que puisque $\Delta^{+}_{\lambda}=\Delta_j^+$
($j\in\{1,2,3\}$), on a
$\pi_{\lambda}=\pi^{q_{\Lambda}}(\mathfrak{e})$, avec
$\Lambda=-s_{12}\lambda-\alpha_{32}$. Maintenant on \'{e}nonce le
premier th\'{e}or\`{e}me de cette section.

\begin{theo}
$$
  \text{(\textrm{i})}\ \ \pi_{\lambda}\vert_{B_1}=\pi^{q_{\Lambda}}(\mathfrak{e})\vert_{B_1}\cong
  \dim(\mathfrak{H}_{-}^{q_{\Lambda}})\mathrm{T}_{+}\oplus
  \dim(\mathfrak{H}_{+}^{q_{\Lambda}})\mathrm{T}_{-}$$

$$\text{(\textrm{ii})}
  \ \ \pi_{\lambda}\vert_{B}=\pi^{q_{\Lambda}}(\mathfrak{e})\vert_{B}\cong
  \sum_{m\in\Z}\left[\dim((\mathfrak{H}_{-})^{q_{\Lambda}}_{\sigma_{m}})\mathrm{T}_{m,+}\oplus\dim((\mathfrak{H}_{+})^{q_{\Lambda}}_{\sigma_{m}})\mathrm{T}_{m,-}\right]$$

\end{theo}

\noindent\begin{demo} On va d'abord montrer (\textrm{i}): Puisque
$H^{q_{\Lambda}}(\mathfrak{e}_{\mathfrak{b}_1},V_{\Lambda}, B_1)$ est
$B_1$-invariant pour $1\otimes1\otimes l$, on peut le
d\'{e}sint\'{e}grer par la d\'{e}composition de Plancherel, autrement
dit
 $$H^{q_{\Lambda}}(\mathfrak{e}_{\mathfrak{b}_1},V_{\Lambda}, B_1)\cong \sum_{\pi\in\widehat{B_1}_d}\H_\pi\otimes\text{Hom}_{B_1}(\H_\pi,H^{q_{\Lambda}}(\mathfrak{e}_{\mathfrak{b}_1},V_{\Lambda}, B_1)).$$

Où $\widehat{B_1}_d$ d\'{e}signe l'ensemble des s\'{e}ries
discr\`{e}tes de $B_1$, $\H_\pi$ est l'espace dans lequel $\pi$ agit
et
$\text{Hom}_{B_1}(\H_\pi,H^{q_{\Lambda}}(\mathfrak{e}_{\mathfrak{b}_1},V_{\Lambda},
B_1))$ est l'espace des op\'{e}rateurs d'intrelacement entre $\H_\pi$
et $H^{q_{\Lambda}}(\mathfrak{e}_{\mathfrak{b}_1},V_{\Lambda},
B_1)$. Dans notre cas, $B_1$ ne poss\`{e}de que 2 s\'{e}ries
discr\`{e}tes: $\mathrm{T}_{-}$ et $\mathrm{T}_{+}$ dont l'une est la
contragr\'{e}diente de l'autre.  En effet, dune part $B_1$ a
uniquement deux orbites coadjointes ouvertes, $\Omega_{+}$ et
$\Omega_{-}$, et d'autre part, on vérifie facilement que la
désintégration de la restriction de $T_{+}$ à $N$ ne fait intervenir
que les représentations de caractère central $\exp tE_{2}\mapsto e^{iat}$ avec
$a>0$. Pour d\'{e}crire concr\`{e}tement l'espace
$\text{Hom}_{B_1}(\H_\pi,H^{q_{\Lambda}}(\mathfrak{e}_{\mathfrak{b}_1},V_{\Lambda},
B_1))$, on peut suivre la m\'{e}thode de Schmid pour les groupes de
Lie semi-simples ( [30] ), sauf que dans notre cas, il faut faire un
ajustement concernant la fonction module, puisque le groupe $B_1$
n'est pas unimodulaire. En fait sous l'identification de Plancherel
 $$ L^2({B_1})\cong \sum_{\pi\in\widehat{B_1}_d}\H_\pi\otimes \H_{\pi^{\vee}}, \ \text{où}\  \pi^{\vee} \ \text{est la représentation contragrédiente de}\ \pi , $$

 \noindent la différentielle de l'action de $B_1$ dans
 $\H_{\pi^{\vee}}$ ne correspond pas à l'action régulière à droite $R$
 (dans $L^2({B_1})$), mais correspond à $R-(\text{tr
   ad}_{\mathfrak{b}_{1}}).\text{id}/2$, et le reste est expliqué dans
 la section 5 de [11].

 On peut démontrer (\textrm{ii}) similairement, et il suffit de
 remarquer que la représentation contragrédiente de
 $\mathrm{T}_{m,\pm}$ est $\mathrm{T}_{-m,\mp}$.
\end{demo}\\

 \noindent \textbf{Remarque}. Il est clair que
 $\dim(\mathfrak{H}_{\pm}^{q_{\Lambda}})=\sum_{m\in\Z}\dim((\mathfrak{H}_{\pm})^{q_{\Lambda}}_{\sigma_{m}})$\\

Maintenant, on commence par traiter les séries discrètes ni
holomorphes ni anti-holomorphes, sauf indication contraire, on garde
toutes les notations concernant $G=SU(2,1)$ dans les sections et
chapitres précédents. Supposons donc que $\pi_{\lambda}$ est une
s\'{e}rie discr\`{e}te ni holomorphe ni anti-holomorphe, c'est-à-dire
que $\lambda$ correspond à $\Delta_3^+$: $\lambda(H_{12})\in \N^+$ et
$\lambda(H_{13})\in \N^+$ avec $\lambda(H_{12}) > \lambda(H_{13})$, on
peut vérifier facilement que dans ce cas, $q_{\Lambda}=1$. D'autre
part, puisque $f_0(H)=\lambda(H_{12})$ et
$f_0(Z)=2\lambda(H_{13})-\lambda(H_{12})$, on déduit facilement que
$\lambda$ correspond à $\Delta_3^+$ si et seulement si $f_0(H)\in
\N^+$, $f_0(H)+f_0(Z)\in 2\N^+$ et $|f_0(H)|>|f_0(Z)|$. Dans la suite
on va exprimer les résultats concernés en fonction de $f_0(H)$ et
$f_0(Z)$. \\

 On va d'abord construire le système différentiel ordinaire (d'ordre
 1) $D^{\pm}_m$ sur $]0,+\infty[$ pour chaque $m\in \N$: Si
    $m\geqslant f_0(H)$, $D^{\pm}_m$ a $2f_{0}(H)$ fonctions inconnues, et
    si $0\leqslant m\leqslant f_0(H)-1$, $D^{\pm}_m$ a $2m+1$
    fonctions inconnues. Plus pr\'{e}cis\'{e}ment: si $m\geq f_0(H)$,
    le syst\`{e}me $D^{\pm}_m$ est le suivant:
$$(y^{\pm}_m)'(t)=(t^{-1}A^{\pm}_{m}+t^{-2}B^{\pm}_{m}+t^{-3}C^{\pm}_{m})y^{\pm}_m(t)\ \ \ \ (t\in]0,+\infty[) \ \   (D^{\pm}_m,\ \ \ m\geqslant f_0(H))$$
\noindent o\`{u}
$$y^{\pm}_m(t)=\left(\begin{array}{c} a^{\pm}_{m,1}(t)\\ \vdots\\a^{\pm}_{m,n}(t)\\ \vdots\\ a^{\pm}_{m,2f_0(H)}(t)\end{array}\right)
$$
$$A^{\pm}_{m}=\left(\begin{array}{ccccccc}r^{\pm}_{m}& & & & & &\\ &A^{\pm}_{m,0}& & & & &\\ & & \ddots & & & &\\ & & & A^{\pm}& & & \\ & & & & \ddots& & \\ & & & & & A_{m,(f_0(H)-2)}& \\ & & & & & & s^{\pm}_{m}\end{array}\right)$$

avec $r^{-}_{m}=s^{+}_{m}=\frac{f_0(Z)-f_0(H)}{2}$, $s^{-}_{m}=r^{+}_{m}=-(\frac{f_0(Z)+f_0(H)}{2})$ et
$$A^{-}_{m,n}=\left(\begin{array}{cc}-(n+1+\frac{f_0(Z)-f_0(H)}{2})& \sqrt{2}i(n+1)\\ -\sqrt{2}i(f_0(H)-n-1)& (n+1+\frac{f_0(Z)-f_0(H)}{2})\end{array}\right),\ \ 0\leq n\leq f_0(H)-2$$

$$A^{+}_{m,n}=\left(\begin{array}{cc}(\frac{f_0(Z)+f_0(H)}{2}-n-1)& -\sqrt{2}i(n+1)\\ \sqrt{2}i(f_0(H)-n-1)& (n+1-\frac{f_0(Z)+f_0(H)}{2})\end{array}\right),\ \ 0\leq n\leq f_0(H)-2$$

$$B^{\pm}_{m}=\left(\begin{array}{ccccc}B^{\pm}_{m,0}& & & &\\& \ddots& & &\\& &B^{\pm}_{m,n}& &\\& & &\ddots&\\& & & & B^{\pm}_{m,f_{0}(H)-1}\end{array}\right) $$
\noindent avec
$$B^{\pm}_{m,n}=\left(\begin{array}{cc}0&\pm 1\\ \pm 2(m-n)&0\end{array}\right),\ \ 0\leq n\leq f_0(H)-1$$

\noindent et

$$C^{\pm}_m=\left(\begin{array}{ccccc}C^{\pm}_{m,0}& & & &\\& \ddots& & &\\& &C^{\pm}_{m,n}& &\\& & &\ddots&\\& & & & C^{\pm}_{m,f_{0}(H)-1}\end{array}\right)$$ avec $$C^{\pm}_{m,n}=\left(\begin{array}{cc}-1&0\\ 0&1\end{array}\right),\ \ 0\leq n\leq f_0(H)-1$$

Et si $0\leqslant m\leqslant f_0(H)-1$, le syst\`{e}me $D^{\pm}_m$ est le suivant: $$(y^{\pm}_m)'(t)=(t^{-1}A^{\pm}_{m}+t^{-2}B^{\pm}_{m}+t^{-3}C^{\pm}_{m})y^{\pm}_m(t)\ \ \ \ (t\in]0,+\infty[) \ \   (D^{\pm}_m,\ \ \ 0\leqslant m\leqslant f_0(H)-1)$$
\noindent o\`{u} $$y^{\pm}_m(t)=\left(\begin{array}{c} a^{\pm}_{m,1}(t)\\ \vdots\\a^{\pm}_{m,n}(t)\\ \vdots\\ a^{\pm}_{m,2m+1}(t)\end{array}\right)
$$
 $$A^{\pm}_{m}=\left(\begin{array}{cccccc}r^{\pm}_{m}& & & & & \\ &A^{\pm}_{m,0}& & & & \\ & & \ddots & & & \\ & & & A^{\pm}_{m,n}& &  \\ & & & & \ddots&  \\ & & & & & A^{\pm}_{m,m-1} \\\end{array}\right)$$ avec $r^{-}_{m}=\frac{f_0(Z)-f_0(H)}{2}$, $r^{+}_{m}=-(\frac{f_0(Z)+f_0(H)}{2})$ et
$$A^{-}_{m,n}=\left(\begin{array}{cc}-(n+1+\frac{f_0(Z)-f_0(H)}{2})& \sqrt{2}i(n+1)\\ -\sqrt{2}i(f_0(H)-n-1)& (n+1+\frac{f_0(Z)-f_0(H)}{2})\end{array}\right),\ \ 0\leq n\leq m-1$$

$$A^{+}_{m,n}=\left(\begin{array}{cc}(\frac{f_0(Z)+f_0(H)}{2}-n-1)& -\sqrt{2}i(n+1)\\ \sqrt{2}i(f_0(H)-n-1)& (n+1-\frac{f_0(Z)+f_0(H)}{2})\end{array}\right),\ \ 0\leq n\leq m-1$$

$$B^{\pm}_{m}=\left(\begin{array}{cccccc}B^{\pm}_{m,0}& & & & &\\& \ddots& & & &\\& &B^{\pm}_{m,n}& & &\\& & &\ddots& &\\& & & & B^{\pm}_{m,m-1}&\\& & & & & 0\\\end{array}\right) $$
 \noindent avec
 $$B^{\pm}_{m,n}=\left(\begin{array}{cc}0&\pm 1\\ \pm 2(m-n)&0\end{array}\right),\ \ 0\leq n\leq m-1$$
\noindent et

 $$C^{\pm}_m=\left(\begin{array}{cccccc}C^{\pm}_{m,0}& & & & &\\& \ddots& & & &\\& &C^{\pm}_{m,n}& & &\\& & &\ddots& &\\& & & & C^{\pm}_{m,m-1} &\\ & & & & &-1\\\end{array}\right)$$

 \noindent avec

 $$C^{\pm}_{m,n}=\left(\begin{array}{cc}-1&0\\ 0&1\end{array}\right),\ \ 0\leq n\leq m-1.$$

 Définissons alors

 $$(D^{\pm}_{m})^{\infty}=\left\{y_{m}\ \text{solution de} \ D^{\pm}_{m};\ \ \int^{+\infty}_{0}\frac{\Vert y_{m}(t) \Vert^2}{t}dt<+\infty\right\}.$$

 Maintenant on énonce le deuxième théorème de cette section:

 \begin{theo} Pour $\pi_{\lambda}$ ni holomorphe ni anti-holomorphe avec $f_0(H)\in \N^+$, $f_0(H)+f_0(Z)\in 2\N^+$ et $|f_0(H)|>|f_0(Z)|$, on a que

 $$\pi_{\lambda}\vert_{B}\cong\sum_{m=0}^{+\infty}\dim((D^{-}_{m})^{\infty}).\mathrm{T}_{[3m-\frac{(3f_0(H)+f_0(Z))}{2}],+}\bigoplus$$

 $$\sum_{m=0}^{+\infty}\dim((D^{+}_{m})^{\infty}).\mathrm{T}_{[-(3m+ \frac{f_0(Z)-3f_0(H)}{2})],-} $$

 \noindent et  $$\pi_{\lambda}\vert_{B_1}\cong \left(\sum_{m=0}^{+\infty}\dim((D^{-}_{m})^{\infty})\right)\mathrm{T}_+ \bigoplus \left(\sum_{m=0}^{+\infty}\dim((D^{+}_{m})^{\infty})\right)\mathrm{T}_- $$
 \end{theo}

\noindent \textbf{Remarque}. Puisque $\dim((D^{\pm}_{m})^{\infty})\leqslant 2N_1$, on déduit surtout que $\pi_{\lambda}$ est $B$-admissible.\\

\noindent \begin{demo}
D'après le théorème précédent et la remarque qui le suit, il s'agit de montrer:

Pour $m\in\N$ et $l=\pm(3m-\frac{(3f_0(H)\pm f_0(Z))}{2})$

$$\dim((\mathfrak{H}_{\mp})^{1}_{\sigma_{l}})=\dim((D^{\mp}_{m})^{\infty})$$
\noindent et pour les autres $l\in \Z$,

$$\dim((\mathfrak{H}_{\pm})^{1}_{\sigma_{l}})=0.$$\\

On va d'abord traiter
$(\mathfrak{H}_{-})^{1}=\ker(\delta_{-})_{1}\cap\ker
{(\delta_{-})_{1}^*}$.  Notons $H_1\triangleq S,Z_1\triangleq -E_2/2,
X_1\triangleq E_1'/\sqrt{2}$ et $Y_1\triangleq E_1/\sqrt{2}$, alors
dans $\mathfrak{b}_{1}=\mathfrak{a}\oplus\mathfrak{n}$, $[X_1,Y_1]=Z_1$,
$[H_1,X_1]=X_1$, $[H_1,Y_1]=Y_1$ et $[H_1,Z_1]=2Z_1$. Notons
$J_0=H_1+iZ_1$, et $J_1=X_1+iY_1$, alors $\C J_0\oplus\C
J_1=\mathfrak{h}_{+}$. Donc si on note $\varphi_0$, $\varphi_1$ la
base duale de $J_0$ et $J_1$, tous les \'{e}l\'{e}ments dans
$\Lambda^1(\mathfrak{h}_{+})^*\otimes V_{\lambda+\alpha_{31}}\otimes
\H_{\pi_1}^\infty$ s'\'{e}crivent sous la forme $\varphi_0\otimes
\upsilon_0+\varphi_1\otimes \upsilon_1$ avec $\upsilon_0,
\upsilon_1\in V_{\lambda+\alpha_{31}}\otimes \H_{\pi_1}^\infty$, et on
peut v\'{e}rifier directement que les \'{e}quations
 $$(\delta_{-})_{1}( \varphi_0\otimes \upsilon_0+\varphi_1\otimes
\upsilon_1 )=0 \ \ \text{dans} \ \Lambda^2(\mathfrak{h}_{+})^*\otimes
V_{\lambda+\alpha_{31}}\otimes \H_{-}^\infty$$

$$(\delta_{-})_{1}^*((\varphi_0\otimes \upsilon_0+\varphi_1\otimes
\upsilon_1)=0 \ \ \text{dans} \ V_{\lambda+\alpha_{31}}\otimes
\H_{-}^\infty$$ sont \'{e}quivalentes à

$$J_0.\upsilon_1=J_1.\upsilon_0+\upsilon_1$$ et

$$J_0^*.\upsilon_0+J_1^*.\upsilon_1=0.$$

Dans les 2 derni\`{e}res \'{e}quations, $J_0, J_1$ sont en tant
qu'op\'{e}rateurs dans $ V_{\lambda+\alpha_{31}}\otimes
\H_{-}^\infty$(par rapport \`{a} la représentation $( \ \xi+(\text{tr
  ad}\ ).\text{id})\otimes \mathrm{T}_-$) et $J_0^*,J_1^*$ sont leurs
adjoints formels respectifs.  Donc $\ker(\delta_{-})_{1}\cap\ker
{(\delta_{-})_{1}^*}$ est l'ensemble des \'{e}l\'{e}ments
$\varphi_0\otimes \upsilon_0+\varphi_1\otimes \upsilon_1$ qui v\'{e}rifient
les 2 derni\`{e}res \'{e}quations. Dans la suite, on va d\'{e}terminer
explicitement les actions de $J_0, J_1$ et $J_0^*, J_1^*$ dans
$V_{\lambda+\alpha_{31}}\otimes \H_{-}^\infty$. Pour ceci, il suffit
de d\'{e}terminer celles dans $V_{\lambda+\alpha_{31}}$ et dans
$\H_{-}^\infty$ repectivement. On va d'abord d\'{e}terminer leurs
actions dans $\H_{-}^\infty$.\\

On a déjà vu que $\mathrm{T}_-\cong \text{Ind}_{N}^{AN}\rho( n_{-},
\mathfrak{l}_{-}, N )$ (voir la section
\ref{description-représentations}), où
$\rho(\ n_{-},\ \mathfrak{l}_{-},\ N )$ est une induite holomorphe.

Comme $n_{-}=-E_2^*$, on a $n_{-}(X_1)=n_{-}(Y_1)=0$ et $n_{-}(Z_1)>0$
, sans perte de g\'{e}n\'{e}ralit\'{e}, on peut supposer que
$n_{-}=Z_1^*$ (par rapport \`{a} la base duale $X_1^*,\ Y_1^*,
\ Z_1^*$ dans $\mathfrak{n}^*$ ), car
$\mathrm{T}_-\cong\text{Ind}_{N}^{AN}
\rho(n_{-},\mathfrak{l}_{-},N)\cong \text{Ind}_{N}^{AN} \rho(r
n_{-},\mathfrak{l}_{-},N)$, $\forall r>0$. Soit $F_{1}$ un élément de
l'espace de la représentation $\rho(n_{-},\mathfrak{l}_{-},N)$, et
$g=\exp(xX_1+yY_1+zZ_1)\in N$. Il est clair que
$F_1(g)=e^{-iz}F_1(\exp(xX_1+yY_1))$ (rappelons que par d\'{e}finition
$\forall X\in \mathfrak{l}_{-}$, on a $X*F_1+if_0(X)F_1=0$). On peut
donc identifier $F_1$ à une fonction $\C\longrightarrow\C :
x+iy\longmapsto F_1(\exp(xX_1+yY_1))$. Pour simplifier, on la note
encore $F_1$, et sous cette identification, par des calculs directs,
on peut obtenir que $X_1*F_1=\frac{\partial F_1}{\partial
  x}+\frac{i}{2}yF_1$, $Y_1*F_1=\frac{\partial F_1}{\partial
  y}-\frac{i}{2}xF_1$, et $Z_1*F_1=-iF_1$. On en d\'{e}duit
que $$(\frac{\partial}{\partial x}+i\frac{\partial}{\partial
  y})F_1+\frac{1}{2}(x+iy)F_1=0\ \ \ \ \ (\Delta).$$ On peut
v\'{e}rifier que $F_0(x+iy)=e^{-\frac{1}{4}(x^2+y^2)}$ est une
solution particuli\`{e}re de $(\Delta)$. Puisque $F_0$ ne s'annule en
aucun point, $F_1$ s'\'{e}crit sous la forme
$F_1=F_2.F_{0} $ avec $F_2$ holomorphe (car $(\frac{\partial}{\partial
  x}+i\frac{\partial}{\partial y})F_2=0$). Donc $\rho( n_{-},
\mathfrak{l}_{-}, N )$ peut se r\'{e}aliser dans l'espace
hilbertien $$\H_{N,Z_1^*}=\left\{F;\ \ F \text{ est holomorphe
  et} \int_{\C}{\vert F(\omega)\vert}^2
e^{-\frac{{\vert\omega\vert}^2 }{2}}d\omega < \infty\right\},$$ où
$\omega=x+iy$.

Notons $\H_{AN,Z_1^*}$ l'espace hilbertien o\`{u} agit
$\text{Ind}_{N}^{AN} \rho( n_{-}, \mathfrak{l}_{-}, N) \cong
\mathrm{T}_-$ et soit $f\in \H_{AN,Z_1^*} $. Vu que $AN/N=A$, $f$ est
une fonction de $AN$ dans $\H_{N,Z_1^*}$ qui v\'{e}rifie $$f(gn)=\rho(
n_{-}, \mathfrak{l}_{-}, N)(n^{-1})f(g), \text{et} \int_{A}\Vert
f(g)\Vert^2 dg < \infty .$$

Il est clair que le groupe $R_+^*$ (pour la multiplication) est
isomorphe \`{a} $A=\exp\R H_1$ par $t\longmapsto \exp(\log t)H_1$, et
la mesure de Haar pour $R_+^*$ est $\frac{dt}{t}$ (ici $dt$ est la
mesure de Lebesgue). Donc on peut identifier $\H_{AN,Z_1^*}$ \`{a}
l'espace hilbertien $$\left\{f: R_+^*\times \C\longrightarrow \C| \ f
\ \text{est mesurable et pour presque tout} \ t\in R_+^*,
\omega\longmapsto f(t,\omega) \right.$$ $$\left.\text{est holomorphe
  et} \int_{R_+^*\times \C}\vert f(t,\omega)\vert^2
e^{-\frac{{\vert\omega\vert}^2
  }{2}}\frac{dt}{t}d\omega<\infty\right\}.$$

On note cet espace encore $\H_{AN,Z_1^*}$ et $\pi_{AN}$ la
repr\'{e}sentation associée (qui est \'{e}quivalente à
$\mathrm{T}^-$). D'autre part, on voit facilement que  $f\in
\H_{AN,Z_1^*} $ correspond à une application $\phi\R^{*}_{+}\mapsto \H(
n_{-}, \mathfrak{l}_{-}, N)$ (voir la section
\ref{description-représentations}),
avec $f(t,\omega)=\phi(t)(g).e^{\frac{1}{4}(x^2+y^2)}.e^{iz}$ o\`{u}
$\omega=x+iy$ et
$g=\exp(xX_1+yY_1+zZ_1)$. Puisque $$\exp(-uH_1).\exp(\log
tH_1).\exp(xX_1+yY_1+zZ_1)$$ $$=\exp(\log(e^{-u}t)H_1).\exp(xX_1+yY_1+zZ_1),$$
on d\'{e}duit facilement que pour $f\in
\H_{AN,Z_1^*}$, $$(\pi_{AN}(\exp(uH_1))f(t,\omega)=f(e^{-u}t,\omega).$$
Donc si de plus $f\in \H_{AN,Z_1^*}^\infty$, on a
$(d\pi_{AN}(H_1)f)(t,\omega)=-t\frac{\partial f}{\partial
  t}(t,\omega)$.

De m\^{e}me, comme $$\exp(-uZ_1).\exp(\log
tH_1)\exp(xX_1+yY_1+zZ_1)$$ $$=\exp(\log t
H_1)\exp(xX_1+yY_1+(z-ut^{-2})Z_1),$$ on obtient que
$(\pi_{AN}(\exp(uZ_1))f(t,\omega)=e^{iut^{-2}}f(t,\omega)$ . Donc
$(d\pi_{AN}(Z_1)f)(t,\omega)=it^{-2}f(t,\omega)$ pour $f\in
\H_{AN,Z_1^*}^\infty$.

De la m\^{e}me fa\c{c}on, on peut obtenir que $$(\pi_{AN}(\exp(uX_1+vY_1))f)(t,\omega)=e^{\left(\frac{1}{2}t^{-1}\omega(u-iv)-\frac{t^{-2}}{4}(u^2+v^2)\right)}f(t,\omega-t^{-1}(u+iv)).$$
Donc $$(d\pi_{AN}(X_1)f)(t,\omega)=\frac{1}{2}t^{-1}\omega f(t,\omega)-t^{-1}\frac{\partial f}{\partial \omega}(t,\omega)$$ et $$(d\pi_{AN}(Y_1)f)(t,\omega)=-\frac{i}{2}t^{-1}\omega f(t,\omega)-it^{-1}\frac{\partial f}{\partial \omega}(t,\omega).$$
Donc $$(d\pi_{AN}(J_0)f)(t,\omega)=-t^{-2}f(t,\omega)-t\frac{\partial f}{\partial t}(t,\omega)$$
$$(d\pi_{AN}(J_1)f)(t,\omega)=t^{-1}\omega f(t,\omega).$$

Avant de d\'{e}terminer $J_0^*$ et $J_1^*$, on va d'abord d\'{e}crire
l'action du groupe compact $M$ dans $\H_{AN,Z_1^*}$. On a vu que
$(\nu(m)\varphi)(x)=\varphi(m^{-1}xm)$ pour $\varphi\in
\H(h_0,\mathfrak{b}_1^-,AN)$ et $m\in M$. On a aussi vu $M=\exp{\R
  W}$, avec $\exp{\text{ad}(-wW)}.X_1$ $=(\cos w)X_1+(\sin w)Y_1$ et
$\exp{\text{ad}(-wW)}$ $.Y_1=-(\sin w)X_1+(\cos w)Y_1$, où $$W=\left(
\begin{array}{ccc}
i/3& 0&0 \\
 0 & -2i/3&0\\
0&0&i/3
\end{array}
\right) $$ est d\'{e}fini dans le chapitre 3.

Donc si $m=\exp{wW}$ et $g=\exp(hH_1).\exp(xX_1+yY_1+zZ_1))$, alors
$m^{-1}gm$ $=\exp(hH_1).\exp((x\cos w-y\sin w)X_1+(x\sin w +y\cos
w)Y_1+zZ_1)$. 
On en déduit que pour $f\in \H_{AN,Z_1^*}$, on a
$(m.f)(t,\omega)=f(t,\omega e^{i(w)})$.  Donc les \'{e}l\'{e}ments dans
$\H_{AN,Z_1^*}$ qui sont de la forme $f(t)\omega^n$ sont des vecteurs
poids de $M$, ($m.(f(t)\omega^n)=(f(t)\omega^n).e^{inw}$), et le
sous-espace engendr\'{e} par ce genre d'\'{e}l\'{e}ments est dense
dans $\H_{AN,Z_1^*}$. En fait, on peut v\'{e}rifier directement que
$\langle f(t)\omega^n,g(t)\omega^m\rangle=0$ si $n\neq m$ où
$\langle,\rangle$ est le produit scalaire de $\H_{AN,Z_1^*}$, et tout
\'{e}l\'{e}ment $f(t,\omega)\in \H_{AN,Z_1^*} $ s'\'{e}crit
$f(t,\omega)=\sum_{n=0}^{+\infty} f_n(t)\omega^n$.

Puisque $\int_{\C}{\vert \omega\vert}^{2n}e^{-\frac{{\vert \omega \vert}^2}{2}}d\omega=2\pi.2^n.n!$, on en d\'{e}duit que $$\H_{AN,Z_1^*}=\left\{\sum_{n=0}^{+\infty} f_n(t)\omega^n;\ \ \sum_{n=0}^{+\infty}(2^n.n!.\int_{R_+^*}\frac{{\vert f_n(t) \vert}^2}{t}dt) <\infty\right\}.$$

Maintenant on va calculer $J_0^*$ et $J_1^*$. Comme $\pi_{AN}$ est unitaire, pour tout $X\in \mathfrak{h}\oplus \mathfrak{n}$, on a $d\pi_{AN}(X)^*=-d\pi_{AN}(X)$. Donc on en d\'{e}duit imm\'{e}diatement le calcul de $J_0^*$ et $J_1^*$: $$J_0^*=d\pi_{AN}(H_1+iZ_1)^*=d\pi_{AN}(-H_1+iZ_1)=-t^{-2}+t\frac{\partial }{\partial t}$$  $$J_1^*=d\pi_{AN}(X_1+iY_1)^*=d\pi_{AN}(-X_1+iY_1)=2t^{-1}\frac{\partial}{\partial \omega}.$$

Maintenant, on va calculer l'action de $J_0$, $J_1$, $J_0^*$ et $J_1^*$ dans $V_{\lambda+\alpha_{31}}$. Les op\'{e}rateurs sont induits par la repr\'{e}sentation $\xi+(\text{tr ad}).\text{id}$ de $\mathfrak{h}_{+}=\C J_0+\C J_1$ dans $V_{\lambda+\alpha_{31}}$. D'abord, comme dans la section \ref{Rossi-Vergne}, on peut d\'{e}terminer que $\dim(V_{\lambda+\alpha_{31}})=f_0(H)$ (rappelons que $f_0(H)=-i\lambda(i H_{12})=\lambda(H_{12})$). D'autre part puisque $\mathfrak{k}\cong \mathfrak{su}(2)\oplus z(\mathfrak{k})$ o\`{u} $z(\mathfrak{k})$ est le centre de $\mathfrak{k}$, et les repr\'{e}sentations irr\'{e}ductibles de $\mathfrak{su}(2)$ sont d\'{e}termin\'{e}es par la dimension, et $z(\mathfrak{k})$ agit scalairement. Il est bien connu que l'on peut r\'{e}aliser $V_{\lambda+\alpha_{31}}$ dans $\C[\omega_1,\omega_2]$ l'espace des polyn\^{o}mes complexes homog\`{e}nes (2 variables) de degr\'{e} $f_0(H)-1$, et le produit scalaire ${\langle,\rangle}_1$ est d\'{e}fini par ${\langle P,Q\rangle}_1=\int_{\C^2}P\overline{Q}e^{-\frac{(\vert \omega_1\vert^2+\vert \omega_2\vert^2)}{2}}d\omega_1d\omega_2$. Notons $v_n=\omega_1^{n}.\omega_2^{(f_0(H)-1-n)}$ ($n=0,1..,f_0(H)-1$). Alors $\tau_{\lambda+\alpha_{31}}(H_{12})v_n=(f_0(H)-1-2n)v_n$ et $\tau_{\lambda+\alpha_{31}}(E_{12})v_n=-nv_{n-1}$ (par convention, $v_{-1}=0$). Rappelons que  $$H_{12}=\left(
\begin{array}{ccc}
1& 0&0 \\
 0 & -1&0\\
0&0&0
\end{array}
\right)$$
et $$E_{12}=\left(
\begin{array}{ccc}
0& 1&0 \\
 0& 0&0\\
0&0&0
\end{array}
\right)$$

Notons $$Z_0'=\left(
\begin{array}{ccc}
 1/2&0&0\\
 0&1/2&0\\
 0&0&-1
 \end{array}
\right) $$
alors $\C Z_0'=z(\mathfrak{k})^{\C}$. Donc $Z_0'$ agit scalairement dans $\C[\omega_1,\omega_2]$.  Puisque $Z_0'=H_{13}-H_{12}/2$ et $(\lambda+\alpha_{31})(H_{13}-H_{12}/2)=(\lambda+\alpha_{31})(H_{13}-H_{12}/2)=\frac{f_0(Z)}{2}-3/2$, \textbf{on note cette quantit\'{e} $\alpha$}, $\tau_{\lambda+\alpha_{31}}(Z_0')=\alpha.\text{id}$.

 Puisque $$J_0=H_1+iZ_1=H_1-iE_2/2=H_{13}+2E_{31}$$ avec $$H_{13} \in \mathfrak{k}^{\C},\ \ E_{31}\in \mathfrak{p}_{+},$$ on a $$\xi(J_0)=\tau_{\lambda+\alpha_{31}}(H_{13}).$$

Or, $H_{13}=\frac{H_{12}}{2}+Z_0'$, donc $\xi(J_0)=\tau_{\lambda+\alpha_{31}}(\frac{H_{12}}{2}+Z_0')$. D'autre part, par des calculs directs, on voit que $\text{tr ad}=4H_1^*$, o\`{u} $H_1^*$ est l'\'{e}l\'{e}ment par rapport à la base duale  $H_1^*,\  X_1^*,\ Y_1^*,\ Z_1^ *$ dans $\mathfrak{b}_1^*$, et $H_1^*(J_0)=1$. Donc on en d\'{e}duit que $$J_0(v_n)=((\xi+2H_1^*\text{id})(J_0))(v_n)=(\frac{f_0(H)-1-2n}{2}+\alpha+2)v_n.$$

De m\^{e}me, $$J_1=E_1'/\sqrt{2}+iE_1/\sqrt{2}=-\sqrt{2}i E_{12}-\sqrt{2}i E_{32}$$ avec $E_{12}\in \mathfrak{k}_{\C}$ et $E_{32}\in\mathfrak{p}_{+}$.
Donc $\xi(J_1)=\tau_{\lambda+\alpha_{31}}(-\sqrt{2}iE_0)$, or $\text{tr ad}(J_1)=0$, et on a $$J_1(v_n)=-\sqrt{2}i(\tau_{\lambda+\alpha_{31}}(E_0))(v_n)=\sqrt{2}inv_{n-1}.$$
Par des calculs directs, on obtient que ${\langle v_n, v_m\rangle}_1=0$ si $n\neq m$, et ${\langle v_n, v_n\rangle}_1= 4\pi^2.2^{f_0(H)-1}.n!.(f_0(H)-1-n)!$. Donc on en d\'{e}duit que $J_0^*=J_0$ ( car $\frac{f_0(H)-1-2n}{2}-\alpha+2\in \R$ ), et $J_1^*(v_{n-1})=-\sqrt{2}i(f_0(H)-n)v_n$.

Maintenant, on peut calculer $\ker(\delta_{-})_{1}\cap\ker{(\delta_{-})_{1}}^*$. Mais avant d'effectuer les calculs explicitement, on remarque que les vecteurs $v_n$ constituent une base de vecteurs propres pour le groupe compact $M$ dans $V_{\lambda+\alpha_{31}}$. Puisque $W=iH_{12}/2-iZ_0'/3$, on obtient que  $$\exp(wW)v_n=e^{iw(\frac{f_0{H}-1-2n}{2}-\frac{\alpha}{3})}v_n=e^{iw(\frac{(3f_0(H)-f_0(Z))}{6}-n)}v_n.$$

Dans la suite, on va calculer $\ker(\delta_{-})_{1}\cap\ker{(\delta_{-})_{1}}^*$. On a vu que $\ker(\delta_{-})_{1}\cap\ker{(\delta_{-})_{1}}^*$ est constitu\'{e} des \'{e}l\'{e}ments $\varphi_{0}\otimes \mu_0+\varphi_{1}\otimes \mu_1$ avec $\mu_0, \ \mu_1\in V_{\lambda+\alpha_{31}}\otimes \H_{-}^\infty$ qui v\'{e}rifient $$J_0.\mu_1=J_1.\mu_0+\mu_1$$ et $$J_0^*.\mu_0+J_1^*.\mu_1=0.$$

On peut \'{e}crire $\mu_0=\sum_{n=0}^{f_0(H)-1}v_n\otimes\phi_n $ et $\mu_1=\sum_{n=0}^{f_0(H)-1}v_n\otimes\psi_n$ avec $\phi_n, \psi_n\in\H_{-}^\infty$.

Donc, on a $J_0.\mu_1=\sum_{n=0}^{f_0(H)-1}v_n\otimes J_0(\psi_n)+\sum_{n=0}^{f_0(H)-1}J_0(v_n)\otimes\psi_n $. D'apr\`{e}s ce qui pr\'{e}c\`{e}de, on obtient que $$J_0.\mu_1=\sum_{n=0}^{f_0(H)-1}v_n\otimes(-t\frac{\partial \psi_n}{\partial t}-t^{-2}\psi_n+\frac{f_0(H)-2n+2\alpha +3}{2}\psi_n).$$

De m\^{e}me  $$J_1.\mu_0=\sum_{n=0}^{f_0(H)-1}v_n\otimes t^{-1}\omega.\phi_n +\sqrt{2}i\sum_{n=0}^{f_0(H)-1}v_{n-1}\otimes n\phi_n ,$$ donc $$J_1.\mu_0+\mu_1=\sum_{n=0}^{f_0(H)-1}v_n\otimes(t^{-1}\omega.\phi_n+\psi_n+\sqrt{2}i(n+1)\phi_{n+1}).$$

Ici par convention $\phi_{-1}=\phi_{f_0(H)}=\psi_{f_0(H)}=\psi_{-1}=0$, et cette convention s'applique dans toute la suite de la section. Donc on en d\'{e}duit que l'\'{e}quation $$J_0.\mu_1=J_1.\mu_0+\mu_1$$ est \'{e}quivalente \`{a} $$-t\frac{\partial \psi_n}{\partial t}-t^{-2}\psi_n+\frac{f_0(H)-2n+2\alpha +3}{2}\psi_n=t^{-1}\omega.\phi_n+\psi_n+\sqrt{2}i(n+1)\phi_{n+1}$$ pour $0\leq n \leq f_{0}(H)-1$, c'est-à-dire $$-t\frac{\partial \psi_n}{\partial t}+(-t^{-2}+\frac{f_0(H)-2n-2\alpha +1}{2})\psi_n=t^{-1}\omega.\phi_n+\sqrt{2}i(n+1)\phi_{n+1} (*).$$
De la m\^{e}me mani\`{e}re, on peut obtenir que  $$J_0^*.\upsilon_0=\sum_{n=0}^{f_0(H)-1} v_n\otimes(t\frac{\partial \phi_n}{\partial t}-t^{-2}\phi_n+\frac{f_{0}(H)-2n+2\alpha +3}{2}\phi_n)$$ et $$J_1^*.\mu_1=\sum_{n=0}^{f_0(H)-1}v_n\otimes2t^{-1}\frac{\partial \psi_n}{\partial \omega} -\sqrt{2}i\sum_{n=0}^{f_0(H)-1}v_{n+1}\otimes(f_0(H)-1-n)\psi_n .$$
Donc $$J_0^*.\mu_0+J_1^*.\upsilon_1=0$$ est \'{e}quivalente \`{a} $$2t^{-1}\frac{\partial \psi_n}{\partial \omega}-\sqrt{2}i(f_0(H)-n)\psi_{n-1}+t\frac{\partial \phi_n}{\partial t}+(-t^{-2}+\frac{f_0(H)-2n+2\alpha +3}{2})\phi_n=0\ (**)$$ pour $0\leq n \leq f_0(H)-1$.

On \'{e}crit $\phi_n=\sum_{m=0}^{+\infty}a_{m,n}(t)\omega^m$ et
$\psi_n=\sum_{m=0}^{+\infty}b_{m,n}(t)\omega^m$, $0\leq n \leq
f_0(H)-1$. Donc $$-t\frac{\partial \psi_n}{\partial
  t}=\sum_{m=0}^{+\infty}-tb_{m,n}'(t)\omega^m$$ $$(-t^{-2}+\frac{f_{0}(H)-2n+2\alpha
  +1}{2})\psi_n=\sum_{m=0}^{+\infty}(-t^{-2}+\frac{f_0(H)-2n+2\alpha
  +1}{2})b_{m,n}(t)\omega^m,$$

$$t^{-1}\omega.\phi_n=\sum_{m=0}^{+\infty}t^{-1}a_{m,n}(t)\omega^{m+1}$$ et $$\sqrt{2}i(n+1)\phi_{n+1}=\sum_{m=0}^{+\infty}\sqrt{2}i(n+1)a_{m,n+1}(t)\omega^m.$$ Donc de $(*)$, on d\'{e}duit que $$-tb_{m,n}'(t)+(-t^{-2}+\frac{f_0(H)-2n+2\alpha +1}{2})b_{m,n}(t)-t^{-1}a_{m-1,n}(t)-$$ $$\sqrt{2}i(n+1)a_{m,n+1}(t)=0\ \ \ \ \ \ \ \ \ \ \ \  (\triangleright)$$ pour $m\geq 0,\ 0\leq n \leq f_0(H)-1$. Ici par convention $a_{-1,n}=b_{-1,n}=a_{m,f_0(H)}=b_{m,f_0(H)}=0$ et cette convention s'applique dans toute la suite.

De m\^{e}me, on peut obtenir de $(**)$ que $$ta_{m,n}'(t)+(-t^{-2}+\frac{f_0(H)-2n+2\alpha +3}{2})a_{m,n}(t)-$$ $$\sqrt{2}i(f_0(H)-n)b_{m,n-1}(t)+2t^{-1}(m+1)b_{m+1,n}(t)=0\ \ \ \ \ \  (\triangleright\triangleright).$$

Puisque la variable $t$ est strictement positive, de $(\triangleright)$ et $(\triangleright\triangleright)$, on d\'{e}duit que pour $m$ fix\'{e}, on a $$b'_{m-n,f_0(H)-1-n}(t)=-t^{-2}a_{m-n-1,f_0(H)-1-n}(t)+$$ $$(-t^{-3}+\frac{2n-f_0(H)+2\alpha+3}{2}t^{-1})b_{m-n,f_0(H)-1-n}(t)-t^{-1}\sqrt{2}i(f_0(H)-n)a_{m-n,f_0(H)-n}(t)$$ et $$a'_{m-n-1,f_0(H)-1-n}(t)=t^{-1}\sqrt{2}i(1+n)b_{m-n-1,f_0(H)-2-n}(t)+$$ $$(t^{-3}-\frac{2n-f_0(H)+2\alpha+5}{2}t^{-1})a_{m-n-1,f_0(H)-1-n}(t)-2t^{-2}(m-n)b_{m-n,f_0(H)-1-n}(t).$$

Or $a_{m,f_0(H)}=a_{m,-1}=b_{m,-1}=b_{m,f_0(H)}=0$. Donc pour chaque
$m\geq 0$, on obtient un syst\`{e}me diff\'{e}rentiel $D^{-}_m$
d'ordre 1 comme expliqué ci-après. Si $m\geqslant f_0(H)$, $D^{-}_m$ a
$2f_0(H)$ fonctions inconnues, et si $0\leqslant m\leqslant f_0(H)-1$,
$D^{-}_m$ a $2m+1$ fonctions inconnues. Plus pr\'{e}cis\'{e}ment: si
$m\geq f_0(H)$ on obtient le syst\`{e}me d'ordre 1 de $2f_0(H)$
fonctions inconnues suivante:
$$(y^{-}_m)'(t)=(t^{-1}A^{-}_{m}+t^{-2}B^{-}_{m}+t^{-3}C^{-}_{m})y^{-}_m(t)\ \ \ \ (t\in]0,+\infty[)
    \ \ (D^{-}_m,\ \ \ m\geqslant f_0(H))$$ o\`{u}
$$y^{-}_m(t)=\left(\begin{array}{c}
      b_{m,f_0(H)-1}(t)\\a_{m-1,f_0(H)-1}(t)\\ \vdots\\b_{m-n,f_0(H)-1-n}(t)\\a_{m-n-1,f_0(H)-1-n}(t)\\ \vdots\\ b_{m-(f_0(H)-1),0}(t)\\ a_{m-f_0(H),0}(t)\end{array}\right)$$
    et $A^{-}_{m}$, $B^{-}_{m}$ et $C^{-}_{m}$ sont les matrices que
    l'on a définies auparavant.

Si $0\leqslant m\leqslant N_1-1$, on obtient le  syst\`{e}me d'ordre 1 de 2m+1 fonctions inconnues suivante : $$(y^{-}_m)'(t)=(t^{-1}A^{-}_{m}+t^{-2}B^{-}_{m}+t^{-3}C^{-}_{m})y^{-}_m(t)\ \ \ \ (t\in]0,+\infty[) \ \   (D^{-}_m,\ \ 0\leqslant m\leqslant f_0(H)-1 )$$ o\`{u} $$y^{-}_m(t)=\left(\begin{array}{c} b_{m,f_0(H)-1}(t)\\a_{m-1,f_0(H)-1}(t)\\ \vdots\\b_{m-n,f_0(H)-1-n}(t)\\a_{m-n-1,f_0(H)-1-n}(t)\\ \vdots\\ b_{1,f_0(H)-m}(t)\\ a_{0,f_0(H)-m}(t)\\ b_{0,f_0(H)-m-1}(t)\end{array}\right)$$ et $A^{-}_{m}$, $B^{-}_{m}$ et $C^{-}_{m}$ sont les matrices que l'on a définies auparavant (pour cette raison-là, on appelle les systèmes encore $D^{-}_m$).

 On voit facilement que les syst\`{e}mes ($\triangleright$) et ($\triangleright \triangleright$) sont \'{e}quivalents à $D^{-}_m$ ($m\geqslant0$). Le point important est que les syst\`{e}mes $D^{-}_m$ sont deux \`{a} deux ind\'{e}pendants  et toute $a_{m,n}\in C^{\infty}$ ( la m\^{e}me chose pour $b_{m,n}$ ) est impliqu\'{e}e dans un et un seul syst\`{e}me $D^{-}_m$.

D'autre part, remarquons que $[W,J_0]=0$ et $[W,J_1]=iJ_1$. On en
déduit successivement que $\exp(wW).J_0=J_0$,
$\exp(wW).J_1=e^{iw}.J_1$,  $\exp(wW).\varphi_0=\varphi_0$ et
$\exp(wW).\varphi_1=e^{-iw}\varphi_1$. D'apr\`{e}s ce qui
pr\'{e}c\`{e}de, il en résulte que

$$\exp(wW).(\varphi_0\otimes v_{f_0(H)-n-1}\otimes a_{m-n-1,f_0(H)-n-1}(t)\omega^{m-n-1} )=$$



$$=e^{iw(m-\frac{(3f_0(H)+f_0(Z))}{6})}\varphi_0\otimes
v_{f_0(H)-n-1}\otimes a_{m-n-1,N_1-n-1}(t)\omega^{m-n-1} $$ et
$$\exp(wW).(J_1^*\otimes v_{f_0(H)-n-1}\otimes b_{m-n,f_0(H)-n-1}(t)\omega^{m-n}  )=$$



$$=e^{iw(m-\frac{(3f_0(H)+f_0(Z))}{6})}\varphi_1\otimes v_{f_0(H)-n-1}\otimes b_{m-n,f_0(H)-n-1}(t)\omega^{m-n}. $$

Pour $m\in N$, désignons par $F^{-}_{m}$ l'espace formé des vacteurs
\begin{equation*}
\begin{split}
\sum_{n=0}^{N_1-1}( \varphi_0\otimes v_{f_0(H)-n-1}\otimes
a_{m-n-1,f_0(H)-n-1}(t)\omega^{m-n-1} \\
+\varphi_1\otimes v_{f_0(H)-n-1}\otimes b_{m-n,f_0(H)-n-1}(t)\omega^{m-n})
\end{split}
\end{equation*}
o\`{u}, $a_{m-n-1,f_0(H)-n-1}(t)$, $\ b_{m-n,f_0(H)-n-1}(t)$ sont
solutions de $D^{-}_m$ avec $m\geqslant 0$. D\'{e}finissons
\'{e}galement les espaces $$(F^{-}_m)^{\infty}=\left\{g_{m} \in
F^{-}_m;\ \ \int_{\R_+^*}\frac{\Vert g_{m}(t)
  \Vert^2}{t}dt<+\infty\right\}.$$ Rappelons que l'on a déjà défini

$$(D^{-}_m)^{\infty}=\left\{y_{m} \ \text{solution de} \ D^{-}_m;\ \ \int_{\R_+^*}\frac{\Vert y_{m}(t) \Vert^2}{t}dt<+\infty\right\}.$$

Alors il est clair que $\text{dim}(F^{-}_m)^{\infty}=\text{dim}(D^{-}_m)^{\infty}$. D'autre part, d'apr\`{e}s ce qui pr\'{e}c\`{e}de, on déduit facilement que pour $m\in\N$ et $l=(3m-\frac{(3f_0(H)+ f_0(Z))}{2})$

$$\dim((\mathfrak{H}_{-})^{1}_{\sigma_{l}})=\dim((F^{-}_{m})^{\infty})$$
\noindent et pour les autres $l\in \Z$,

$$\dim((\mathfrak{H}_{-})^{1}_{\sigma_{l}})=0.$$ \\

Il reste donc à traiter $(\mathfrak{H}_{+})^{1}=\ker(\delta_{+})_{1}\cap\ker {(\delta_{+})_{1}^*}$. Comme pour  $(\mathfrak{H}_{-})^{1}$, on peut obtenir que $\ker(\delta_{+})_{1}\cap\ker {(\delta_{+})_{1}^*}$ est constitu\'{e} des \'{e}l\'{e}ments $\varphi_{0} \otimes\mu_0+\varphi_1\otimes\mu_1$ avec $\mu_0,\mu_1\in V_{\lambda+\alpha_{31}}\otimes \mathbb{H}_{+}^\infty $ qui v\'{e}rifient $$J_0.\mu_1=J_1.\mu_0+\mu_1$$ et $$J_0^*.\mu_0+J_1^*.\mu_1=0.$$ Ici, comme dans $(\mathfrak{H}_{-})^{1}$, $J_0=H_1+iZ_1$ et $J_1=X_1+iY_1$, et $H_1,Z_1,X_1,Y_1$ sont comme ceux que l'on a d\'{e}finis auparavant. $J_0^*,J_1^*$ sont leurs adjoints formels respectivement. D'autre part, puisque $\mathrm{T}_{+}\cong \mathrm{T}_{-}^{\vee}$, en identifiant $\mathbb{H}_{+}$ \`{a} $\mathbb{H}_{-}$ (donc $\mathbb{H}_{+}^\infty$ \`{a} $\mathbb{H}_{-}^\infty$), les op\'{e}rateurs $J_0=H_1+iZ_1$ et $J_1=X_1+iY_1$ dans $\mathbb{H}_{+}^\infty$ sont \'{e}quivalents \`{a} $H_1-iZ_1$ et $X_1-iY_1$ respectivement dans  $\mathbb{H}_{-}^\infty$ (et $J_0^*,J_1^*$ aux $-H_1-iZ_1$ et $-X_1-iY_1$ respectivement). Dans la suite on va garder ces identifications, et travailler dans $V_{\lambda+\alpha_{31}}\otimes\mathbb{H}_{-}^\infty $.

On écrit $\mu_0=\sum_{n=0}^{f_0(H)-1}v_n\otimes\phi_n $ et $\mu_1=\sum_{n=0}^{f_0(H)-1}v_n\otimes\psi_n $ avec $\phi_n, \psi_n\in\H_{-}^\infty$, où les $v_n\in V_{\lambda+\alpha_{31}}$ sont ceux que l'on a d\'{e}finis auparavant.

 On obtient que $$J_0.\mu_1=\sum_{n=0}^{f_0(H)-1}v_n\otimes(-t\frac{\partial \psi_n}{\partial t}+t^{-2}\psi_n+\frac{f_0(H)-2n+2\alpha +3}{2}\psi_n) $$ et $$J_1.\mu_0=\sum_{n=0}^{f_0(H)-1}v_n\otimes-2t^{-1}\frac{\partial \phi_n}{\partial \omega} +\sqrt{2}i\sum_{n=0}^{f_0(H)-1}v_{n-1}\otimes n\phi_n .$$

Donc $$J_1.\mu_0+\mu_1=\sum_{n=0}^{f_0(H)-1}v_n\otimes(-2t^{-1}\frac{\partial \phi_n}{\partial \omega}+\psi_n+\sqrt{2}i(n+1)\phi_{n+1}),$$ ici par convention $\phi_{-1}=\phi_{f_0(H)}=\psi_{f_0(H)}=\psi_{-1}=0$ bien entendu. Donc on en d\'{e}duit que l'\'{e}quation $$J_0.\mu_1=J_1.\mu_0+\mu_1$$ est \'{e}quivalente \`{a} $$-t\frac{\partial \psi_n}{\partial t}+(t^{-2}+\frac{f_0(H)-2n+2\alpha +1}{2}\psi_n=-2t^{-1}\frac{\partial \phi_n}{\partial \omega}+\sqrt{2}i(n+1)\phi_{n+1}\ \ (1')$$ pour $0\leq n \leq f_0(H)-1$.

De m\^{e}me, on a $$J_0^*.\mu_0=\sum_{n=0}^{f_0(H)-1}v_n\otimes(t\frac{\partial \phi_n}{\partial t}+t^{-2}f_n+\frac{f_0(H)-2n+2\alpha +3}{2}\phi_n) $$ et $$J_1^*.\mu_1=\sum_{n=0}^{f_0(H)-1}v_n\otimes-t^{-1}\omega .\psi_n  -\sqrt{2}i\sum_{n=0}^{f_0(H)-1}v_{n+1}\otimes(f_0(H)-1-n)\psi_n .$$

Donc $$J_0^*.\mu_0+J_1^*.\mu_1=0$$ est \'{e}quivalente \`{a} $$-t^{-1}\omega .\psi_n-\sqrt{2}i(f_0(H)-n)\psi_{n-1}+t\frac{\partial \phi_n}{\partial t}+(t^{-2}+\frac{f_0(H)-2n+2\alpha +3}{2})\phi_n=0\ \ \ (2') $$ pour $0\leq n \leq f_0(H)-1$.

On \'{e}crit $\phi_n=\sum_{m=0}^{+\infty}a_{m,n}(t)\omega^m$ et $\psi_n=\sum_{m=0}^{+\infty}b_{m,n}(t)\omega^m$, $0\leq n \leq f_0(H)-1$.
Alors, de la même manière que dans l'étude de $(\text{ker}\delta_1)\bigcap(\text{ker}\delta_1^*)$,  on obtient de $(1')$ et $(2')$ que $$a_{m,n}'(t)=\sqrt{2}it^{-1}(N_1-n)b_{m,n-1}-$$ $$(t^{-3}+\frac{N_1-2n+2\alpha+3}{2})a_{m,n}(t)+t^{-2}b_{m-1,n}(t)\ \ (\triangleleft)$$ et $$b_{m,n}'(t)=2(m+1)t^{-2}a_{m+1,n}(t)+$$ $$(t^{-3}+\frac{N_1-2n+2\alpha+1}{2}t^{-1})b_{m,n}(t)-\sqrt{2}i(n+1)t^{-1}a_{m,n+1}(t) \ \ (\triangleleft\triangleleft).$$

Donc, comme dans l'étude de $\ker(\delta_{-})_{1}\cap\ker
{(\delta_{-})_{1}^*}$, pour chaque $m\geq 0$, on obtient un
syst\`{e}me diff\'{e}rentiel $D^{+}_m$ d'ordre 1. Si $m\geqslant N_1$,
$D^{+}_m$ a $2N_1$ fonctions inconnues, et si $0\leqslant m\leqslant
f_0(H)-1$, $D^{+}_m$ a $2m+1$ fonctions inconnues. Plus
pr\'{e}cis\'{e}ment: si $m\geqslant f_0(H)$, $D^{+}_m$ est le système
suivant :

$$(y^{+}_m)'(t)=(t^{-1}A^{+}_{m}+t^{-2}B^{+}_{m}+t^{-3}C^{+}_{m})y^{+}_m(t)\ \ \ \ (t\in]0,+\infty[) \ \   (D^{+}_m,\ \ \ m\geqslant f_0(H))$$ o\`{u} $$y^{+}_m(t)=\left(\begin{array}{c} a_{m,0}(t)\\b_{m-1,0}(t)\\ \vdots\\a_{m-n,n}(t)\\b_{m-n-1,n}(t)\\ \vdots\\ a_{m-(f_0(H)-1),f_0(H)-1}(t)\\ b_{m-f_0(H),f_0(H)-1}(t)\end{array}\right)$$ et les $A^{+}_{m}, B^{+}_{m}, C^{+}_{m}$ sont les matrices que l'on a définies auparavant.

 Si $0\leqslant m\leqslant f_0(H)-1$, $D^{+}_m$ est le système suivant : $$(y^{+}_m)'(t)=(t^{-1}A^{+}_{m}+t^{-2}B^{+}_{m}+t^{-3}C^{+}_{m})y^{+}_m(t)\ \ \ \ (t\in]0,+\infty[) \ \   (D^{+}_m,\ \ 0\leqslant m\leqslant f_0(H)-1 )$$ o\`{u} $$y^{+}_m(t)=\left(\begin{array}{c} a_{m,0}(t)\\b_{m-1,0}(t)\\ \vdots\\a_{m-n,n}(t)\\b_{m-n-1,n}(t)\\ \vdots\\ a_{1,m-1}(t)\\ b_{0,m-1}(t)\\ a_{0,m}\end{array}\right)$$ et les $A^{+}_{m}, B^{+}_{m}, C^{+}_{m}$ sont les matrices que l'on a définies auparavant (également pour cette raison-là, on appelle les systèmes $D^{+}_m$).

 Comme dans $\ker(\delta_{-})_{1}\cap\ker {(\delta_{-})_{1}^*}$, on
 voit facilement que les syst\`{e}mes ($\triangleleft$)et
 ($\triangleleft \triangleleft$) sont \'{e}quivalents à $D^{+}_m$
 ($m\geqslant0$).  De plus les syst\`{e}mes $D^{+}_m$ sont deux \`{a}
 deux ind\'{e}pendants et toute $a_{m,n}\in C^{\infty}$ (la m\^{e}me
 chose pour $b_{m,n}$) et est impliqu\'{e}e dans un et un seul
 syst\`{e}me $D^{+}_m$.

 D'autre part, il faut remarquer que puisque $\mathrm{T}_{+}\cong \mathrm{T}_{-}^{\vee}$, on a  $$\exp(wW).a_{m,n}(t)\omega^m=e^{-iwm}a_{m,n}(t)\omega^m.$$ 

Donc on en d\'{e}duit que  $$\exp(wW).(\varphi_0\otimes v_{n}\otimes a_{m-n,n}(t)\omega^{m-n} )=$$

$$=e^{-iw(m+\frac{(f_0(Z)-3f_0(H))}{6})}\varphi_0\otimes v_{n}\otimes
a_{m-n,n}(t)\omega^{m-n} ,$$
et

$$\exp(wW).(\varphi_1\otimes v_{n}\otimes b_{m-n-1,n}(t)\omega^{m-n-1}  )=$$

$$=e^{-iw(m+\frac{(f_0(Z)-3f_0(H))}{6})}\varphi_1\otimes v_{n}\otimes b_{m-n-1,n}(t)\omega^{m-n-1}. $$
D\'{e}finissons les espaces
$$F^{+}_m=\left\{\sum_{n=0}^{f_0(H)-1}(\varphi_0\otimes v_{n}\otimes a_{m-n,n}(t)\omega^{m-n} +\varphi_1\otimes v_{n}\otimes b_{m-n-1,n}(t)\omega^{m-n-1} )\right\}$$ o\`{u}, $a_{m-n,n}(t)$, $b_{m-n-1,n}(t)$  sont solutions de $D^{+}_m$ avec $m\geqslant 0$. D\'{e}finissons \'{e}galement les espaces $$(F^{+}_m)^{\infty}=\left\{g_{m} \in E_m;\ \ \int_{\R_+^*}\frac{\Vert g_{m}(t) \Vert^2}{t}dt<+\infty\right\}$$ et rappelons que l'on a déjà défini

$$(D^{+}_m)^{\infty}=\left\{y_{m}\ \text{solution de} \ D^{+}_m;\ \ \int_{\R_+^*}\frac{\Vert y_{m}(t) \Vert^2}{t}dt<+\infty\right\}.$$

 Alors il est clair que $\dim((F^{+}_m)^{\infty})=\dim((D^{+}_m)^{\infty})$. D'autre part, d'apr\`{e}s ce qui pr\'{e}c\`{e}de, on déduit facilement que pour $m\in\N$ et $l=-(3m-\frac{(3f_0(H)- f_0(Z))}{2})$

$$\dim((\mathfrak{H}_{+})^{1}_{\sigma_{l}})=\dim((F^{+}_{m})^{\infty})$$
\noindent et pour les autres $l\in \Z$,

$$\dim((\mathfrak{H}_{+})^{1}_{\sigma_{l}})=0.$$ \\

La démonstration est donc bien achevée.
\end{demo}\\

Maintenant on va retraiter le cas s\'{e}ries discr\`{e}tes holomorphes à
l'aide du théorème 6.7. On verra que non seulement on peut retrouver
les résultats de Rossi-Vergne pour $G=SU(2,1)$ (théorème 6.3), mais
cette m\'{e}thode est particuli\`{e}rement simple pour traiter les
s\'{e}ries discr\`{e}tes holomorphes (et celles anti-holomorphes) de
$G$.

Supposons maintenant que $\pi_{\lambda}$ est une s\'{e}rie disc\`{e}te
holomorphe avec le param\`{e}tre de Harish-Chandra $\lambda$ qui
vérifie que $\lambda(H_{12})\in \N^+$ et $\lambda(H_{31})\in \N^+$. On
peut facilement vérifier que dans ce cas $q_\Lambda=0$. On a déjà vu
dans la section \ref{Rossi-Vergne} que dire que $\pi_{\lambda}$ est
une série discrète holomorphe équivaut à dire que $f_0(H)\in \N^+$ et
$f_0(Z)+f_0(H)$ est un entier pair strictement négatif. Dans la suite,
on va exprimer les résultats en fonction de $f_0(H)$ et $f_0(Z)$.\\

 Il est bien clair que $\ker{(\delta_{\pm})_{0}}^*=\Lambda^0(\mathfrak{h}_{+})^*\otimes V_{\lambda+\alpha_{31}}\otimes \H_{\pm}^\infty$. Donc $\mathfrak{H}_{\pm}^{0}=\ker{(\delta_{\pm})_{0}}\cap\ker{(\delta_{\pm})_{0}}^*=\ker{(\delta_{\pm})_{0}}$. On peut facilement voir que $$\ker{(\delta_{\pm})_{0}}=\left\{\mu\in V_{\lambda+\alpha_{31}}\otimes \H_{\pm}^\infty;\ \ J_0{\mu}=0\ \ \text{et}\ \ J_1{\mu}=0\right\},$$ ici les op\'{e}rateurs $J_0$ et $J_1$ sont ceux que l'on a d\'{e}finis auparavant.

D'abord, on va calculer $\ker{(\delta_{-})_{0}}$. \'{E}crivons $\mu=\sum_{n=0}^{N_1-1}v_n\otimes \phi_n $  avec $\phi_n\in\H_{-}^\infty$, et les $v_n\in V_{\lambda+\alpha_{31}}$ sont ceux que l'on a d\'{e}finis auparavant. Alors \`{a} l'aide des informations que l'on a eues pour les s\'{e}ries discr\`{e}tes ni holomorphes ni anti-holomorphes, on obtient que les condition $J_0{\mu}=0$ et $J_1{\mu}=0$ sont \'{e}quivalentes \`{a} $$-t\frac{\partial \phi_n}{\partial t}-t^{-2}\phi_n+\frac{f_0(H)-2n+2\alpha+3}{2}\phi_n=0,$$ et $$t^{-1}\omega \phi_n+\sqrt{2}i(n+1)\phi_{n+1}=0$$ respectivement. Ici, comme pour les s\'{e}ries discr\`{e}tes ni holomorphes ni anti-holomorphes $\alpha=\frac{f_0(Z)}{2}-3/2$. Donc en \'{e}crivant $\phi_n=\sum_{m=0}^{+\infty}$ $a_{m,n}(t)\omega^m$, on obtient de la premi\`{e}re \'{e}quation que $a_{m,n}(t)\in \mathbb{C}t^{\frac{f_0(H)+f_0(Z)}{2}-n}e^{\frac{t^{-2}}{2}}$. Or

$$\int_{0}^{+\infty}\frac{(\mid t^{\frac{f_0(H)+f_0(Z)}{2}-n}e^{\frac{t^{-2}}{2}}\mid)^2}{t}dt=+\infty.$$

Donc $\mathfrak{H}_{-}^{0}=\ker(\delta_{-})_{0}={0}$.

Maintenant, il reste \`{a} calculer $\mathfrak{H}_{+}^{0}=\ker(\delta_{+})_{0}$. Comme dans $\mathfrak{H}_{-}^{0}$, en identifiant $\H_{+}$ \`{a} $\H_{-}$, on obtient que $$\ker(\delta_{+})_{0}=\left\{\mu=\sum_{n=0}^{f_0(H)-1}v_n\otimes \phi_n \right\},$$ o\`{u} $\phi_n\in\H_{\pi_1}^\infty$ v\'{e}rifient
$$-t\frac{\partial \phi_n}{\partial t}+t^{-2}\phi_n+\frac{f_0(H)-2n+2\alpha+3}{2}\phi_n=0$$ et

$$-2t^{-1}\frac{\partial \phi_n}{\partial \omega}+\sqrt{2}(n+1)i\phi_{n+1}=0.$$

 Donc, en \'{e}crivant
 $\phi_n=\sum_{m=0}^{+\infty}a_{m,n}(t)\omega^m$, on obtient de la
 premi\`{e}re \'{e}quation que $a_{m,n}(t)\in
 \mathbb{C}t^{\frac{f_0(H)+f_0(Z)}{2}-n}e^{-\frac{t^{-2}}{2}}$, et de
 la seconde \'{e}quation, on d\'{e}duit que
 $\phi_{f_0(H)-1}=a_{0,f_0(H)-1}$ et
 $a_{m,n}(t)=\sqrt{2}(n+1)it\frac{a_{m-1,n+1}(t)}{2m}$. On peut
 v\'{e}rifier directement que si $a_{m-1,n+1}(t)\in
 \mathbb{C}t^{\frac{f_0(H)+f_0(Z)}{2}-(n+1)}e^{-\frac{t^{-2}}{2}}$,
 alors $\sqrt{2}(n+1)it\frac{a_{m-1,n+1}(t)}{2m}\in
 \mathbb{C}t^{\frac{f_0(H)+f_0(Z)}{2}-n}e^{-\frac{t^{-2}}{2}}$. D'autre
 part, puisque $\frac{f_0(H)+f_0(Z)}{2}$ est en entier strictement
 négatif, on a $$\int_{0}^{+\infty}\frac{(\mid
   t^{\frac{f_0(H)+f_0(Z)}{2}-n}e^{-\frac{t^{-2}}{2}}\mid)^2}{t}dt<+\infty.$$

 On en d\'{e}duit donc $\mathfrak{H}_{+}^{0}=\ker(\delta_{+})_{0}\cong\bigoplus_{n=0}^{f_0(H)}\mathbb{C}v_n=V_{\lambda+\alpha_{31}}$. Donc d'apr\`{e}s le théorème 6.7, on obtient exactement le r\'{e}sultat du théorème 6.3.

\subsection{ Étude asymptotique des solutions des systèmes différentiels $D^{\pm}_{m}$}

Dans cette section, on va étudier le comportement asymptotique des solutions des systèmes $D^{\pm}_{m}$, on va d'abord traiter $D^{-}_{m}$.

 En faisant un changement de variable $z=1/t$, on peut transformer le syst\`{e}me $D^{-}_{m}$ dans la section précédente en $$x_m'(z)=\{z^{-1}(-A^{-}_{m})+(-B^{-}_{m})+z(-C^{-}_{m})\}x_m(z)\ \ \ \ ((D^{-}_{m})') \ \ (\text{où}\ \ x_m(z)=y_m(1/z))$$

 Il est clair que $$\int_{0}^{+\infty}\frac {\parallel x_m(z)\parallel^2}{z}dz= \int_{0}^{+\infty}\frac {\parallel y_m(t)\parallel^2}{t}dt.$$

Puisque au voisinage de $0$, le syst\`{e}me $(D^{-}_{m})'$ a une singularité de première esp\`{e}ce, il est bien connu qu'au voisinage de $0$, il existe une solution fondamentale $X_m(z)=U_m(z)z^{\Delta^{-}_m}z^{N_m}$, o\`{u} $U_m(z)$ est une transformation analytique en $0$ ( c-est-\`{a}-dire $U_m(z)$ est analytique au voisinage de $0$ et $U_m(0)$ est une matrice inversible ), $\Delta^{-}_m$ est une matrice constante diagonale qui est \'{e}quivalente \`{a} la partie semi-simple de $-A^{-}_{m}$, et $N_m$ est une matrice constante strictement triangulaire inf\'{e}rieure. D'autre part, par des calculs directs, on voit que $-A^{-1}_{m}$ est diagonalisable et pour $m\geq f_0(H)$ (resp.$0\leq m< f_0(H)$) , elle a $f_0(H)+1$ (resp. $m+1$) valeurs propres positives et $f_0(H)-1$ ( resp. $m$ ) n\'{e}gatives. En fait, par exemple, pour $m\geq f_0(H)$, les valeurs propre de $-A^{-}_{m}$ sont $\frac{f_0(H)-f_0(Z)}{2}, \frac{f_0(H)+f_0(Z)}{2}$ qui sont positives et
$$\pm \sqrt{(n+1+\frac{f_0(Z)-f_0(H)}{2})^2+2(n+1)(f_0(H)-n-1)}=$$

$$\pm\sqrt{-[(n+1)-\frac{f_0(H)+f_0(Z)}{2}]^2+[(\frac{f_0(H)-f_0(Z)}{2})^2+(\frac{f_0(H)+f_0(Z)}{2})^2]},$$
pour $0\leq n\leq f_0(H)-2$. Donc on en d\'{e}duit que le sous-espace
que l'on note $(E^{-}_{m})^{0}$ des solutions $x_m(z)$ de
$(D^{-}_{m})'$ telles que $\int_{0}^{+\infty}\frac {\parallel
  x_m(z)\parallel^2}{z}dz$ converge au voisinage de $0$ est de
dimension $f_0(H)+1$, si $m\geq f_0(H)$, et de dimension $m+1$, si
$0\leq m< f_0(H)$. (Pour plus de d\'{e}tails, voir les résultats de la
section 9.2 de l'appendice qui nous ont été communiqu\'{e}s par
C. Sabbah).

Maintenant on \'{e}tudie le comportement asymptotique des solutions au
voisinage de l'infini. Apr\`{e}s une permutation des coordonnées du
système $(D^{-}_{m})'$ ($m\geqslant f_0(H)$), on r\'{e}\'{e}crit
$$(D^{-}_{m})':
z\widetilde{x_m}'(z)=z^2[\widetilde{C^{-}_{m}}+z^{-1}\widetilde{B^{-}_{m}}+z^{-2}\widetilde{A^{-}_{m}}]\widetilde{x_m}(z)$$
avec $$\widetilde{x_m}(z)=\left(\begin{array}{c}
  b_{m,f_0(H)-1}(z)\\b_{m-1,f_0(H)-2}(z)\\ \vdots\\ b_{m-(f_0(H)-1),0}(t)\\a_{m-1,f_0(H)-1}(z)\\a_{m-2,f_0(H)-2}(z)\\ \vdots\\ a_{m-f_0(H),0}(z)\end{array}\right)$$
$$\widetilde{C_m}=\left(\begin{array}{cc}I_{f_0(H)}&0\\0&-I_{f_0(H)}\end{array}\right),$$ o\`{u}

$I_{f_0(H)}$ est la matrice carr\'{e}e identit\'{e} de taille $f_0(H)$, $$ \widetilde{B^{-}_{m}}=\left(\begin{array}{cc}0&\widetilde{B^{-}_{m}}^{12}\\\widetilde{B^{-}_{m}}^{21}&0\end{array}\right)$$ o\`{u} $\widetilde{B^{-}_{m}}^{12}$ et $\widetilde{B^{-}_{m}}^{21}$ sont des matrices carr\'{e}es de taille $f_0(H)$ et on ne pr\'{e}cise pas leurs formes, puisque elles n'affectent pas notre argument suivant.  $$\widetilde{A^{-}_{m}}=\left(\begin{array}{cc}\widetilde{A^{-}_{m}}^{11}&\widetilde{A^{-}_{m}}^{12}\\\widetilde{A^{-}_{m}}^{21}&\widetilde{A^{-}_{m}}^{22}\end{array}\right)$$

\noindent où les $\widetilde{A^{-}_{m}}^{ij}$($i,j\in\{1,2\}$) sont des matrices carr\'{e}es de taille $f_0(H)$. Pour la m\^{e}me raison, on ne pr\'{e}cise pas leurs formes exactes. Il est connu qu'il existe au voisinage de l'infini, une unique transformation analytique formelle $\widetilde{T_m}(z)$ de la forme $$\widetilde{T_m}(z)=\left(\begin{array}{cc}I_{f_0(H)}&\widetilde{T_m}^{12}(z)\\\widetilde{T_m}^{21}(z)&I_{f_0(H)}\end{array}\right)$$ qui transforme $(D^{-}_{m})'$ en $\widetilde{(D^{-}_{m})'}$ (c'est-\`{a}-dire qu'au voisinage de l'infini, toute solution $\widetilde{x_m}(z)=\widetilde{T_m}(z)\widetilde{y_m}(z)$, o\`{u} $\widetilde{y_m}(z)$ est une solution de $\widetilde{(D^{-}_{m})'}$). Ici le nouveau syst\`{e}me $\widetilde{(D^{-}_{m})'}$ est

$$\widetilde{D_m'}:\ \ z\widetilde{y_m}'(z)=z^2\widetilde{H_m}(z)\widetilde{y_m}(z),$$ avec

$$\widetilde{H_m}(z)=\left(\begin{array}{cc}\widetilde{H_m}^{11}(z)&0\\0&\widetilde{H_m}^{22}(z)\end{array}\right)$$
(donc $\widetilde{D_m'}$ se divise en deux syst\`{e}mes ind\'{e}pendants). Les $\widetilde{T_m}^{ij}(z)$,  $\widetilde{H_m}^{ii}(z)$ ($i,j\in\{1,2\}$) sont analytiques au voisinage de l'infini, et  leurs coefficients peuvent se calculer par récurrence concrètement (pour la proc\'{e}dure de calculs, voir les pages 42 et 43 de l'ouvrage de W.Balser qui s'intitule \og Formal power series and linear systems of meromorphic ordinary differential equations\fg). Ce qui est important dans notre situation, c'est que $\widetilde{H_m}^{11}(z)$ et $\widetilde{H_m}^{22}(z)$ sont de la forme $$\widetilde{H_m}^{11}(z)=[I_{f_0(H)}+z^{-2}\widetilde{G_m}^{11}(z)]$$ et $$\widetilde{H_m}^{22}(z)=[-I_{f_0(H)}+z^{-2}\widetilde{G_m}^{22}(z)],$$ o\`{u} les fonctions matricielles $\widetilde{G_m}^{11}(z)$ et $\widetilde{G_m}^{22}(z)$ sont analytiques au voisinage de l'infini. Donc on peut d\'{e}duire que les solutions pour les deux petits syst\`{e}mes $$z\widetilde{y_m}'^{jj}(z)=z^2\widetilde{H_m}^{jj}(z)\widetilde{y_m}^{jj}(z)\ j=1,2$$ sont de la forme $$\widetilde{y_m}^{11}(z)=\exp{(\frac{x^2}{2})}\widetilde{x_m}^{11}(z)$$ et $$\widetilde{y_m}^{22}(z)=\exp{(-\frac{x^2}{2})}\widetilde{x_m}^{22}(z)$$ o\`{u} $\widetilde{x_m}^{11}(1/z)$ et $\widetilde{x_m}^{22}(1/z)$ sont les solutions de deux syst\`{e}mes de la premi\`{e}re esp\`{e}ce au voisinage de 0 respectivement. D'autre part, il est clair qu'au voisinage de l'infini, les solutions de $(D^{-}_{m})'$ sont de la forme $$y_m(z)=\widetilde{T_m}(z)\left(\begin{array}{c}\widetilde{y_m}^{11}(z)\\\widetilde{y_m}^{22}(z)\end{array}\right).$$

Donc d'apr\`{e}s la forme des solutions fondamentales sur les
syst\`{e}mes de la premi\`{e}re esp\`{e}ce (voir ce que l'on a fait
pour étudier le comportement asymptotique en 0 des solutions de
$(D^{-}_{m})'$. Pour plus de d\'{e}tails, voir la section 9.2 de
l'appendice), on d\'{e}duit que le sous-espace que l'on note
$(E^{-}_{m})^{\infty}$ des solutions $y_m(z)$ de
$\widetilde{(D^{-}_{m})'}$ telles que $\int_{0}^{+\infty}\frac
{\parallel y_m(z)\parallel^2}{z}dz$ converge au voisinage de $+\infty$
sont de la
forme $$\widetilde{T_m}(z)\left(\begin{array}{c}0\\\widetilde{y_m}^{22}(z)\end{array}\right).$$
Donc $(E^{-}_{m})^{\infty}$ est de dimension $f_0(H)$. D'autre part,
il est clair que
$\text{dim}(D^{-}_{m})^{\infty}=\text{dim}((E^{-}_{m})^{0}\cap
(E^{-}_{m})^{\infty})$. Or vu que la dimension de l'espace des
solutions de $(D^{-}_{m})'$ est de $2f_0(H)$, on d\'{e}duit que $1\leq
\dim((D^{-}_{m})^{\infty})\leq f_0(H)$. De même, pour $0\leq m<
f_0(H)$, on peut obtenir que
$\text{dim}((E^{-}_{m})^{\infty})=m$. Puisque la dimension de l'espace
des solutions de $(D^{-}_{m})'$ est de $2m+1$, $
\text{dim}(D^{-}_{m})^{\infty}\leqslant m $. Donc surtout
$\text{dim}(D^{-}_{0})^{\infty}=0$.\\

De la même manière, on peut étudier le comportement asymptotique des
solutions de $D^{+}_{m}$. On peut obtenir que pour $m\geqslant
f_0(H)$, on a que $1\leq \dim((D^{+}_{m})^{\infty})\leq f_0(H)$, et
pour $0\leq m< f_0(H)$, on a que $\dim((D^{+}_{m})^{\infty})\leqslant
m $. Donc surtout $\dim((D^{+}_{0})^{\infty})=0$.

 D'après ce qui précède, on déduit que $\sum_{m=0}^{+\infty}\dim((D^{\pm}_{m})^{\infty})=+\infty$. Selon le théorème 6.8, on a donc $$\pi_{\lambda}\vert_{B_1}\cong \left( +\infty \right)\mathrm{T}_+ \bigoplus \left( +\infty \right)\mathrm{T}_- .$$

Donc d'après le théorème 5.4, on déduit que les assertions
(i) et (ii) de la conjecture de Duflo sont confirmées pour $(G, \pi_{\lambda}, B_1)$, avec $\pi_{\lambda}$ ni holomorphe ni anti-holomorphe. D'autre part dans la section 6.2, on a déjà confirmé les assertions
(i) et (ii) de la conjecture de Duflo pour  $(G, \pi_{\lambda}, B_1)$, avec $\pi_{\lambda}$ holomorphe ou anti-holomorphe. Donc \textbf{les assertions
(i) et (ii) de la conjecture de Duflo sont confirmées pour $(G, \pi_{\lambda}, B_1)$ pour toutes les séries discrètes}.

\subsection{ Interprétation et application des travaux de Fabec et de
  ceux de Kraljevic}

Dans cette section, on va réétudier la décomposition de
$\pi_{\lambda}\vert_{B}$ pour $\pi_{\lambda}$ ni holomorphe ni
anti-holomorphe à l'aide des travaux de Fabec et ceux de Kraljevic.
En combinant les résultats que l'on a obtenus dans les sections 6.4 et
6.5, on arrive à décomposer explicitement $\pi_{\lambda}\vert_{B}$.\\

Soit $B=MAN$ le sous-groupe de Borel de $G=SU(2,1)$ que l'on a défini
auparavant. Soit $R^{p,\alpha}$ la série principale de $G=SU(2,1)$ par
rapport à $B$, o\`{u} $p\in \widehat{M}$ (c'est-\`{a}-dire un
caract\`{e}re unitaire de $M$) et $\alpha\in
\mathfrak{a}_{\C}^*$. Rappelons que $R^{p,\alpha}$ est la
représentation induite par la représentation de $B$ de dimension $1$ :
$man\mapsto \exp((\alpha+\rho)(\log a)p(m))$, où $m\in M$, $a\in A$,
$n\in N$ et
$\rho=\frac{1}{2}\tr(\text{ad}\vert_{\mathfrak{n}})$. Comme dans [20], on identifie $\widehat{M}$ au r\'{e}seau des formes
lin\'{e}aires $f$ dans $\mathfrak{m}^*$ avec $f(W)\in \Z/3$, et
$\widehat{K}$ au r\'{e}seau des formes lin\'{e}aires $f$ dans
$\mathfrak{t}^*$ avec $f(W)\in \Z/3$ et $f(H_{12})\in \N$. On
identifie également $\alpha$ à $\alpha(S)\in \C$, où
$S\in\mathfrak{a}$ est ce que l'on a défini dans le chapitre $3$.

Kraljevic ([19] et [20]) a d\'{e}termin\'{e} explicitement tous les
sous-quotients irr\'{e}ductibles unitarisables pour $R^{p,\alpha}$. En
particulier, il a donn\'{e} une condition necessaire et suffisante
pour qu'un sous-quotient irr\'{e}ductible unitarisable soit une
s\'{e}rie discr\`{e}te et il a aussi montr\'{e} que $R^{p,\alpha}$ a
les m\^{e}me sous-quotients irr\'{e}ductibles (unitarisables ou pas)
que $R^{p,-\alpha}$.

Fabec ([9]) a construit  pour toute $R^{p,\alpha}$ deux sous-espaces
ferm\'{e}s invariants $N_2\subset N_{1}$. Notamment il a d\'{e}montr\'{e}
que \\

\begin{theo}
Si $\alpha=3p+2k$, avec $3p+k\geqslant 1$ et $k\geqslant 1$ (ici $k\in
\Z$), alors $R^{p,\alpha}$ est unitarisable sur $N_2$. De plus en
notant ce sous-quotient unitaire encore $R^{p,\alpha}$, on a
$$R^{p,\alpha}\vert_{B}\cong \left.\{\sum_{n\geqslant 3p+k}^{+\infty}\tau^{p-n,+}\bigoplus \sum_{n\geqslant k}^{+\infty}\tau^{p+n,-}\}\right.$$

o\`{u} $$\tau^{m,+}\cong\sigma_{3m}\otimes \overline{\mathrm{T}^-}\cong \mathrm{T}_{3m,-}$$ et

$$\tau^{m,-}\cong\sigma_{3m}\otimes\overline{\mathrm{T}^+}\cong \mathrm{T}_{3m,+}$$

 $\forall m\in \Z/3$
\end{theo}

\noindent \textbf{Remarque}. Dans l'article de Fabec ([9]), il y a une erreur sur un indice de somme dans le th\'{e}or\`{e}me concern\'{e} (voir le th\'{e}or\`{e}me 6.3 de [9]), et le th\'{e}or\`{e}me pr\'{e}c\'{e}dent est la version corrig\'{e}e.\\

Or d'apr\`{e}s Kraljevic (voir le th\'{e}or\`{e}me 6, la proposition
2 et (iv) de la proposition 3 dans [20]), lorsque $3p+k>1$ et $k>1$
($k\in \Z$ ), on peut montrer que $R^{p,\alpha}$ a un seul
sous-quotient unitarisable irr\'{e}ductible (et il correspond \`{a}
"$\pi_{1,1}(r)$" avec "$j(p,\alpha)=1, j(p,-\alpha)=2$" dans
[19]). Donc on en d\'{e}duit qu'il est exactement celui dans le
th\'{e}or\`{e}me pr\'{e}c\'{e}dent. D'autre part, en appliquant le
th\'{e}or\`{e}me 6 de [19], on peut montrer qu'il correspond \`{a} une
s\'{e}rie discr\`{e}te ni holomorphe ni anti-holomorphe. Plus
pr\'{e}cis\'{e}ment, le param\`{e}tre de Harish-Chandra associ\'{e}
$\lambda$ vérifie que $\lambda(H_{12})=\alpha$ ,
$\lambda(H_{13})=3p+k$ et $\lambda(H_{32})=k$. Donc en appliquant le
th\'{e}or\`{e}me pr\'{e}c\'{e}dent de Fabec, \textbf{on peut
  décomposer explicitement toutes ces s\'{e}ries discr\`{e}tes ni
  holomorphes ni anti-holomorphes dont le paramètre de Harish-Chandra
  $\lambda$ vérifie $\lambda(H_{13})>1$ et $\lambda(H_{32})>1$}. Du
coup il reste \`{a} décomposer les s\'{e}ries discr\`{e}tes ni
holomorphes ni anti-holomorphes $\pi_{\lambda}\vert_{B}$ avec
$\lambda(H_{13})=1$ ou $\lambda(H_{32})=1$.

Maintenant, on consid\`{e}re $R^{\frac{2k+1}{3},-1}$ avec $k\in \Z$ et
$k\neq 0$ et $k\neq -1$. D'apr\`{e}s ([19]), on sait que
$R^{\frac{2k+1}{3},-1}$ a les m\^{e}me sous-quotients
irr\'{e}ductibles que $R^{\frac{2k+1}{3},1}$. Or on peut d\'{e}duire
de ([19]) que $R^{\frac{2k+1}{3},1}$ a trois sous-quotients
irr\'{e}ductibles unitarisables. Plus pr\'{e}cis\'{e}ment, si $k\geq
1$, alors il s'agit de la s\'{e}rie discr\`{e}te holomorphe dont le
param\`{e}tre de Harish-Chandra associ\'{e} $\lambda$ v\'{e}rifie que
$\lambda(H_{12})=k$ ,$\lambda(H_{31})=1$, de la s\'{e}rie discr\`{e}te
ni holomorphe ni anti-holomorphe dont le param\`{e}tre de
Harish-Chandra associ\'{e} $\lambda_{k}$ v\'{e}rifie que
$\lambda_{k}(H_{12})=k+1$, $\lambda_{k}(H_{13})=1$ et d'une repr\'{e}sentation
unitaire irr\'{e}ductible qui n'est pas une s\'{e}rie discr\`{e}te. Et
si $k\leq -2$, alors il s'agit de la s\'{e}rie discr\`{e}te
anti-holomorphe dont le param\`{e}tre de Harish-Chandra associ\'{e}
$\lambda$ v\'{e}rifie que $\lambda(H_{12})=-1-k$, $\lambda(H_{23})=1$,
de la s\'{e}rie discr\`{e}te ni holomorphe ni anti-holomorphe dont le
param\`{e}tre de Harish-Chandra associ\'{e} $\lambda_{k}$ v\'{e}rifie que
$\lambda_{k}(H_{12})=-k$, $\lambda_{k}(H_{32})=1$ et d'une repr\'{e}sentation
unitaire irr\'{e}ductible qui n'est pas une s\'{e}rie
discr\`{e}te. Ainsi il est clair que chaque s\'{e}rie discr\`{e}te ni
holomorphe ni anti-holomorphe dont le param\`{e}tre de Harish-Chandra
associ\'{e} $\lambda$ v\'{e}rifie que $\lambda(H_{13})=1$ ou
$\lambda(H_{32})=1$ figure comme un sous-quotient irr\'{e}ductible
dans une (et une seule) $R^{\frac{2k+1}{3},1}$ (donc dans une et une
seule $R^{\frac{2k+1}{3},-1}$). Or d'apr\`{e}s le th\'{e}or\`{e}me 6.5
de [9], on sait que $R^{\frac{2k+1}{3},-1}$ est unitarisable sur
$N_2$, et, notant encore $R^{\frac{2k+1}{3},-1}$ le sous-quotient
$N_{2}$, on a

si $k\geq 1$, $$R^{\frac{2k+1}{3},-1}|_{B}\cong
\tau^{-\frac{k-1}{3},+}$$

et si $k\leq -2$,

$$R^{\frac{2k+1}{3},-1}\vert_{B}\cong \tau^{-\frac{k+2}{3},-}.$$

Donc on peut d\'{e}duire de la proposition 6.4, du th\'{e}or\`{e}me
6.8 et de la section 6.5 que le sous-quotient unitaire sur $N_2$ est
la repr\'{e}sentation unitaire irr\'{e}ductible qui n'est pas une
s\'{e}rie discr\`{e}te. D'autre part, d'apr\`{e}s le th\'{e}or\`{e}me
6.6 de [9], $R^{\frac{2k+1}{3},-1}$ est unitarisable sur $N_1/N_2$, et

si $k\geq 1$, on a $$R^{\frac{2k+1}{3},-1}\vert_{B}\cong
\left.\sum_{n\geqslant
  0}^{+\infty}\tau^{(\frac{2k+1}{3}+n),-}\bigoplus \sum_{n\geqslant
  0}^{k-1}\tau^{(\frac{2k+1}{3}-n),+} \bigoplus\sum_{n\geqslant
  k+1}^{+\infty}\tau^{(\frac{2k+1}{3}-n),+}\right. $$

si $k\leq -2$, on a $$R^{\frac{2k+1}{3},-1}\vert_{B}\cong \left.\sum_{n\geqslant 0}^{+\infty}\tau^{(\frac{2k+1}{3}-n),+}\bigoplus \sum_{n\geqslant 0}^{-k-2}\tau^{(\frac{2k+1}{3}+n),-} \bigoplus\sum_{n\geqslant -k}^{+\infty}\tau^{(\frac{2k+1}{3}+n),-}\right.$$

 Donc la proposition 6.4 nous permet de d\'{e}duire que le
 sous-quotient unitarisable de $R^{\frac{2k+1}{3},-1}$ sur $N_1/N_2$
 contient forc\'{e}ment la s\'{e}rie discr\`{e}te ni holomorphe ni
 anti-holomorphe de paramètre $\lambda_{k}$. Or d'apr\`{e}s le
 théorème 6.8 et la section 6.5, on sait que si $k\geq 1$, alors la
 s\'{e}rie discr\`{e}te ni holomorphe ni anti-holomorphe de paramètre
 $\lambda_{k}$ ne contient pas le terme "$\tau^{(\frac{2k+1}{3}),+}$"
 et si $k\leq -2$, elle ne contient pas le terme
 "$\tau^{(\frac{2k+1}{3}),-}$". Donc on en d\'{e}duit que le
 sous-quotient unitarisable $N_1/N_2$ de $R^{\frac{2k+1}{3},-1}$ est
 la somme de la s\'{e}rie discr\`{e}te ni holomorphe ni
 anti-holomorphe et de la s\'{e}rie discr\`{e}te holomorphe (ou
 anti-holomorphe). Donc en appliquant encore une fois la proposition
 6.3, on peut d\'{e}composer explicitement dans cette situation
 $\pi_{\lambda}\vert_{B}$.  On arrive donc à d\'{e}composer
 explicitement $\pi_{\lambda}\vert_{B}$, pour toutes les séries
 discr\`{e}te ni holomorphes ni anti-holomorphes (de $G$). On résume
 tout ceci comme un théorème:

\begin{theo}
Pour $\pi_{\lambda}$ ni holomorphe ni anti-holomorphe avec $\lambda(H_{12})=n_1\in \N^+$ et $\lambda(H_{13})=n_2\in \N^+$ ($n_1> n_2$), on a que $$\pi_{\lambda}\vert_{B}\cong
\sum_{m=0}^{+\infty}\mathrm{T}_{[3m+2n_1-n_2],+}\bigoplus  \sum_{m=0}^{+\infty}\mathrm{T}_{[-(3m+ n_1+n_2)],-} $$ c'est-à-dire

$$\pi_{\lambda}\vert_{B}\cong
\sum_{m=0}^{+\infty}\mathrm{T}_{[3m+\frac{3f_0(H)-f_0(Z)}{2}],+}\bigoplus  \sum_{m=0}^{+\infty}\mathrm{T}_{[-(3m+ \frac{3f_0(H)+f_0(Z)}{2})],-}.$$
Donc on en déduit aussi une nouvelle fois le résultat que l'on a obtenu dans la section 6.5 :

$$\pi_{\lambda}\vert_{B_1}\cong \left( +\infty \right)\mathrm{T}_+ \bigoplus \left( +\infty \right)\mathrm{T}_- .$$
\end {theo}

\subsection{ Confirmation des assertions
(i) et (ii) de la conjecture de Duflo pour $G=SU(2,1)$ }

Jusqu'à présent, on a confirmé les assertions (i) et (ii) de la
conjecture de Duflo pour $(G, \pi_{\lambda}, B_1)$ pour toutes les
séries discrètes $\pi_{\lambda}$, et pour $(G, \pi_{\lambda}, B)$,
pour $\pi_{\lambda}$ holomorphe ou anti-holomorphe. Maintenant, le
théorème 6.9 nous permet de les confirmer pour
$(G, \pi_{\lambda}, B)$, pour $\pi_{\lambda}$ ni holomorphe ni
anti-holomorphe.

Soit $\pi_{\lambda}$ ni holomorphe ni anti-holomorphe avec
$\lambda(H_{12})\in \N^+$ et $\lambda(H_{13})\in \N^+$
($\lambda(H_{12})>\lambda(H_{13})$). Comme dans le chapitre 5,
$\mathcal{O}_{\pi_{\lambda}}=G.f_0$ est son orbite coadjointe associée
dans $\mathfrak{g}^*$ avec $f_0=-\lambda$. Donc comme dans la section
6.2, on peut obtenir que $$\{\langle k.f_0,E_2\rangle\vert \ k\in
K\}=\{\langle k.f_0,H\rangle+\langle f_0,Z\rangle\vert \ k\in
K\}$$ $$=\{x_1.\langle f_0,H\rangle+\langle f_0,Z\rangle\vert -1\leq
x_1\leq 1\}=[f_0(Z)-f_0(H), f_0(H)+f_0(Z)].$$

Puisque $|f_0(H)|> |f_0(Z)|$, d'après la formule $(**)$ de la section
5.4 et ce qui la précède, on déduit que

$$\text{p}(G.f_0)\cap \mathfrak{b}_{fr}^*=(\bigsqcup_{r\leq\frac{-(3f_0(H)+f_0(Z))}{6}}\Omega_{r,-})\bigcup (\bigsqcup_{r\geq\frac{3f_0(H)-f_0(Z)}{6}}\Omega_{r,+}).$$

Donc on en déduit que l'ensemble des $B$-orbites coadjointes fortement
r\'{e}guli\`{e}res et admissibles de $\mathfrak{b}^*$ qui sont
contenues dans $\text{p}(\mathcal{O}_{\pi_{\lambda}})$ est

$$\{\Omega_{r_{N},-}: N\in \N\}\bigcup\{\Omega_{l_{N},+}: N\in \N\}$$

o\`{u} $r_{N}=-(\frac{3f_0(H)+f_0(Z)}{6}+\frac{N}{3})-\frac{1}{2}$, et
$l_{N}=\frac{3f_0(H)-f_0(Z)}{6}+\frac{N}{3}+\frac{1}{2}$. On
l'\'{e}crit comme une proposition. Voici\\

\begin{prop} Soit $\pi_{\lambda}$ ni holomorphe ni anti-holomorphe avec $\lambda(H_{12})\in \N^+$ et $\lambda(H_{13})\in \N^+$ ($\lambda(H_{12})>\lambda(H_{13})$), et $\mathcal{O}_{\pi_{\lambda}}$ son orbite coadjointe associée dans $\mathfrak{g}^*$. Alors l'ensemble des $B$-orbites
coadjointes fortement r\'{e}guli\`{e}res et admissibles de
$\mathfrak{b}^*$ qui sont contenues dans
$\text{p}(\mathcal{O}_{\pi_{\lambda}})$ est

$$\{\Omega_{r_{N},-}: N\in \N\}\bigcup\{\Omega_{l_{N},+}: N\in \N\}$$ o\`{u} $r_{N}=-(\frac{3f_0(H)+f_0(Z)}{6}+\frac{N}{3})-\frac{1}{2}$, et $l_{N}=\frac{3f_0(H)-f_0(Z)}{6}+\frac{N}{3}+\frac{1}{2}$. \end{prop}

Puisque pour $m\in \Z$, $\mathrm{T}_{m,\pm}$ correspond à $\Omega_{(\frac{m}{3}\pm\frac{1}{2}),+}$, les théorèmes 5.5 et 6.10 et la proposition 6.11 nous permettent de confirmer les assertions
(i) et (ii) de la conjecture de Duflo  pour $(G, \pi_{\lambda}, B)$, pour $\pi_{\lambda}$ ni holomorphe ni anti-holomorphe.\\

Donc les assertions
(i) et (ii) de la conjecture de Duflo sont confirmées pour $(G, \pi_{\lambda}, B_1)$ et $(G, \pi_{\lambda}, B)$ pour toutes les séries discrètes $\pi_{\lambda}$ de $G=SU(2,1)$.\\

\noindent \textbf{Remarque}. Comme on a vu pour les  séries discrètes holomorphes, pour les séries discrètes ni holomorphes ni holomorphes, il n'y a qu'également un tiers des $B$-orbites coadjointes fortement régulières et admissibles dans la projection $\text{p}(\mathcal{O}_{\pi_{\lambda}})$ qui interviennent dans la décomposition $\mathcal{O}_{\pi_{\lambda}}\vert_{B}$.

\subsection{ Conséquences sur $(D_{m}^{\pm})^{\infty}$ }

Dans le théorème 6.8, on a vu que pour $\pi_{\lambda}$ ni holomorphe ni anti-holomorphe avec $f_0(H)\in \N^+$, $f_0(H)+f_0(Z)\in 2\N^+$ et $|f_0(H)|>|f_0(Z)|$, on a que

 $$\pi_{\lambda}\vert_{B}\cong\sum_{m=0}^{+\infty}\dim((D^{-}_{m})^{\infty}).\mathrm{T}_{[3m-\frac{(3f_0(H)+f_0(Z))}{2}],+}\bigoplus$$

 $$\sum_{m=0}^{+\infty}\dim((D^{+}_{m})^{\infty}).\mathrm{T}_{[-(3m+ \frac{f_0(Z)-3f_0(H)}{2})],-}. $$
  Donc d'après le théorème 6.10, on déduit que

  \begin{prop}
  Pour les systèmes différentiels $D_{m}^{\pm}$ et les sous-espace $(D_{m}^{\pm})^{\infty}$ définis dans la section $5.4$, on a

  $$\text{dim}((D_{m}^{\pm})^{\infty})=0, \ \ \text{si }\ 0\leqslant m\leqslant f_0(H)-1$$

  et

  $$\text{dim}((D_{m}^{\pm})^{\infty})=1, \ \ \text{si }\ m\geqslant f_0(H).$$
  \end{prop}

  Ce résultat semblerait difficilement s'obtenir par des méthodes "directes".

\pagebreak

  \begin{center}\section{Comparaison entre la construction de Duflo et celle d'Auslander-Kostant pour la conjecture}\end{center}

  Dans cette section, on va interpréter la conjecture dans le cadre d'Auslander-Kostant pour les groupes résolubles algébriques.

  D'abord on fait quelques rappels concernant la théorie d'Auslander-Kostant\\

   \textbf{La théorie d'Auslander-Kostant} : \\

   Soit $G$ un groupe réel résoluble algébrique tel que $[G,G]\subset
   G^{\mathfrak{u}}$ d'algèbre de Lie $\mathfrak{g}$, ici
   $G^{\mathfrak{u}}$ est le radical unipotent de $G$. Soit $g\in
   \mathfrak{g}^*$ une forme linéaire. On dit que $g$ est
   \textbf{entière}, s'il existe un caractère unitaire $\chi_{g}$ de
   $G(g)_0$ de différentielle $ig\vert_{\mathfrak{g}(g)}$, où
   $G(g)\subset G$ est le stabilisateur de $g$ par rapport à l'action
   coadjointe de $G$ dans $\mathfrak{g}^*$, $\mathfrak{g}(g)$ est
   l'algèbre de Lie de $G(g)$ et $G(g)_0$ est la composante neutre de
   $G(g)$. On note $\mathfrak{g}^*_{\text{ent}}$ l'ensemble des formes
   linéaires entières. Si $g\in \mathfrak{g}^*_{\text{ent}}$, on note
   $Y(g)$ l'ensemble des représentations unitaires irréductibles de
   $G(g)$ dont la restriction à $G(g)_0$ est un multiple de $\chi_{g}$
   .

   Soient $\mathfrak{l}$ une polarisation en $g\in \mathfrak{g}^*$ et
   $\mathfrak{a}$ un idéal de $\mathfrak{g}$, on dit que
   $\mathfrak{l}$ est admissible pour $\mathfrak{a}$ si
   $\mathfrak{l}\cap \mathfrak{a}_{\C}$ est une polarisation en
   $a=g\vert_{\mathfrak{a}}$.\\

   \begin{theo}

   Soit $g\in \mathfrak{g}^*$, alors il existe une polarisation
   positive $\mathfrak{l}$ en $g$ vérifiant la condition de Pukanszky,
   qui soit $\mathfrak{g}^{\mathfrak{u}}$-admissible ( ici
   $\mathfrak{g}^{\mathfrak{u}}$ est le radical unipotent de
   $\mathfrak{g}$), telle que $\mathfrak{l}$ soit $G(g)$-invariante et
   que $l\cap \mathfrak{g}^{\mathfrak{u}}$ soit
   $G(g\vert_{\mathfrak{g}^{\mathfrak{u}}})$-invariante.

   \end{theo}

   Soit $g\in \mathfrak{g}^*$ et soit $\mathfrak{l}$ une polarisation positive en $g$ vérifiant les conditions du théorème précédent. On définit alors la forme linéaire $\theta\in\mathfrak{g}(g)^*$ par
   $$\theta(X)=\frac{1}{2}\text{Im}\ \text{tr}({\text{ad}X}\vert_{\mathfrak{l}})\ , \ X\in \mathfrak{g}(g) $$
Il s'avère que $\theta$ est indépendant du choix de $\mathfrak{l}$.

 Soit $\widetilde{g}\in \mathfrak{g}^*$ une forme linéaire telle que

   $(i)\ \ \widetilde{g}\vert_{\mathfrak{g}^{\mathfrak{u}}}=g\vert_{\mathfrak{g}^{\mathfrak{u}}}.$

   $(ii)\ \ \widetilde{g}\vert_{\mathfrak{g}(g)}=g\vert_{\mathfrak{g}(g)}-\theta.$

   Une telle forme linéaire existe, car $\mathfrak{g}$ étant
   résoluble, on a
   $(\mathfrak{g}(g))^{\mathfrak{u}}=\mathfrak{g}^{\mathfrak{u}}(g)$,
   tandis que $\theta\vert_{\mathfrak{g}^{\mathfrak{u}}(g)}=0$. Alors
   l'orbite $\widetilde{\mathcal{O}}=G.\widetilde{g}$ de
   $\widetilde{g}$ ne dépend que de l'orbite $\mathcal{O}=G.g$ de $g$
   et l'application $\mathcal{O}\longmapsto\widetilde{\mathcal{O}}$
   induit une bijection de $\mathfrak{g}^*_{\text{ent}}/G$ sur
   $\mathfrak{g}^*_{\text{ad}}/G$, où $\mathfrak{g}^*_{\text{ad}}$ est
   l'ensemble des formes linéaires admissibles (dans le sens de
   Duflo). Rappelons que l'on a déjà défini ce qu'est une forme
   linéaire admissible (au sens de Duflo) dans le chapitre 2.

   Soit $g\in \mathfrak{g}^*_{\text{ent}}$ et $\widetilde{g}\in
   \mathfrak{g}^*_{\text{ad}}$ comme ci-dessus, alors on a
   $G(g)=G(\widetilde{g})$. De plus il existe un caractère
   $\chi_{\theta}$ de $G(\widetilde{g})^{\mathfrak{g}}$ de
   différentielle $i\theta$. Alors l'application
   $\tau\longmapsto\widetilde{\tau}=\chi_{\theta}^{-1}\tau$ induit une
   bijection de $Y(g)$ sur $X(\widetilde{g})$, où $X(\widetilde{g})$
   est ce que l'on a définit dans le chapitre 2.

   Si $g\in\mathfrak{g}^*_{\text{ent}}$ et $\tau\in Y(g)$, Auslander-Kostant lui font correspondre une représentation $\pi_{g,\tau}\in\widehat{G}$.

   Si $g\in\mathfrak{g}^*_{\text{ad}}$ et $\tau\in X(g)$, Duflo lui associe une représentation $T_{g,\tau}\in \widehat{G}$ (ce que l'on a expliqué dans le chapitre 2).\\

   Alors on a le théorème suivant :

   \begin{theo}

   Pour tout $g\in \mathfrak{g}^*_{\text{ent}}$ et tout $\tau\in Y(g)$, on a

   $$\pi_{g,\tau}\cong T_{\widetilde{g},\widetilde{\tau}} $$
   \end{theo}

  Maintenant on revient sur $G=SU(2,1)$, et on garde les notations dans les chapitres précédents qui concernent $G$. On va appliquer la théorie d'Auslander-Kostant à $B=MAN$ et $B_1=AN$.

  D'abord pour $B_1$, on a vu qu'il y a deux $B_1$-orbites coadjointes fortement régulières (dans ${\mathfrak{b}_1}^*$) $\Omega^{\pm}$ qui sont également admissibles et ouvertes. Donc il est évident que $\Omega^{\pm}\subset {\mathfrak{b}_1}^*_{\text{ad}} $. Comme les stabilisateurs concernés sont triviaux, on déduit des paragraphes précédents que $\forall f\in \Omega^{\pm} $, on a $f\in {\mathfrak{b}_1}^*_{\text{ent}}$ et $\widetilde{f}=f$. De plus $Y(f)=X(f)=\{\tau\}$ et $\pi_{f,\tau}=T_{f,\tau}=T_{\pm}$.\\

Maintenant, on traite $B$, on a vu que les $B$-orbites  coadjointes fortement régulières (dans $\mathfrak{b}^*$) sont les $\Omega_{r,\pm}=B.(r W^*\pm E_2^*)$, et $\Omega_{r,\pm}$ est admissible $\Leftrightarrow r+\frac{1}{2}=\Z/3$. Soit $f_{m,\pm}=(\frac{m}{3}\pm\frac{1}{2})W^*\pm E_2^*\in b^*$. On peut vérifier directement que $\mathfrak{l}_{\pm}=\C W\oplus \C (E_1+\pm iE_1')\oplus \C E_2$ est une polarisation positive en $f_{m,\pm}$ qui vérifie le théorème précédent. Puisque $\text{tr}({\text{ad}W}\vert_{\mathfrak{l}_{\pm}})=\mp i$, où $W\in \mathfrak{m}$ est ce que l'on a défini dans le chapitre 3, on déduit que $\widetilde{g_{m,\pm}}=f_{m,\pm}$, avec $g_{m,\pm}=\frac{m}{3}W^*\pm E_2^*\in b^*$. Comme il y a un seul élément dans $X(f_{m,\pm})$, il y en a seul aussi dans $Y(g_{m,\pm})$, et on note $\pi_{g_{m,\pm}}$ la seule représentation unitaire irréductible associée à $g_{m,\pm}$ au sens d'Auslander-Kostant. Alors on a $\pi_{g_{m,\pm}}\cong\mathrm{T}_{m,\pm}$.

Maintenant, on part du paramètre de Blattner $\Lambda$. On a vu que si le paramètre de Harish-Chandra $\lambda$ correspond à une série discrète holomorphe, alors $\Lambda=\lambda+\alpha_{31}$. On peut aussi vérifier facilement que si $\lambda$ correspond à une série discrète ni holomorphe
ni anti-holomorphe, alors $\Lambda=\lambda$.

On suppose d'abord que $\lambda$ correspond à une série discrète holomorphe $\pi_{\lambda}$ tel que $\lambda(H_{12})=n_1\in \N^+$ et $\lambda(H_{31})=n_3\in \N^+$. Notons l'orbite $\mathcal{O}_{\Lambda}=G.f_{\Lambda}$ avec $f_{\Lambda}=-i\Lambda\in \mathfrak{t}^*\subset \mathfrak{g}^*$. On vérifie facilement que $f_{\Lambda}(H_{12})=n_1-1$, $f_{\Lambda}(H_{31})=n_3+2$, $f_{\Lambda}(H)=n_1-1$ et $f_{\Lambda}(Z)=-(n_1+2n_3+3)$. Donc on a $f_{\Lambda}(Z)+f_{\Lambda}(H)<0$. Donc d'après la section 6.2, on déduit que $$\text{p}(\mathcal{O}_{\Lambda})=\bigcup_{\frac{n_3-n_1}{3}+1\leq r\leq \frac{2n_1+n_3}{3}}\Omega_{r,-}.$$

Puisque $\pi_{g_{m,\pm}}\cong\mathrm{T}_{m,\pm}$, la proposition 6.4 nous permet de conclure que \textbf{toute sous-représentation irréductible unitaire dans la décomposition de $\pi_{\lambda}\vert_{B}$ correspond à une (et une seule) $B$-orbite coadjointe fortement régulière et entière (au sens d'Auslander-Kostant) qui est contenue dans $\mathcal{O}_{\Lambda}=G.f_{\Lambda}$}.\\

\noindent \textbf{Remarque}. (1) Ceci n'est pas vrai pour $\mathcal{O}_{\lambda}=G.f_{\lambda}$ avec $f_{\lambda}=-i\lambda$, si l'on utilise la théorie d'Auslander-Kostant pour $B$. D'autre part, Ceci n'est pas vrai non plus pour $\mathcal{O}_{\Lambda}=G.f_{\Lambda}$ si l'on utilise la théorie de Duflo pour $B$.

 (2) Le nombre des orbites entières (au sens d'Auslander-Kostant)  dans $\text{p}(\mathcal{O}_{\Lambda})$ est \textbf{plus petit} que le nombre des orbites admissibles (au sens de Duflo) dans $\text{p}(\mathcal{O}_{\lambda})$.

 Pour les séries discrètes anti-holomorphes, on a les résultats analogues.\\

 Maintenant, on considère le cas les séries discrètes ni holomorphes ni anti-holomorphes. Donc supposons $\lambda(H_{12})=n_1$ et $\lambda(H_{13})=n_2$ avec $n_1>n_2$. Puisque dans ce cas $\lambda=\Lambda$,  on a $\mathcal{O}_{\Lambda}=\mathcal{O}_{\lambda}$. D'autre part, on a vu dans la section 6.7 que

 $$\text{p}(G.f_0)\cap \mathfrak{b}_{fr}^*=(\bigsqcup_{r\leq\frac{-n_1-n_2}{3}}\Omega_{r,-})\bigcup (\bigsqcup_{r\geq\frac{2n_1-n_2}{3}}\Omega_{r,+}). $$

 Donc le théorème 6.10 permet de conclure : \textbf{Toute sous-représentation irréductible unitaire dans la décomposition de $\pi_{\lambda}\vert_{B}$ correspond à une (et une seule) $B$-orbite coadjointe fortement régulière et entière (au sens d'Auslander-Kostant) qui est contenue dans la projection de $\mathcal{O}_{\Lambda}=\mathcal{O}_{\lambda}$}.
\newpage
\null
\newpage

\begin{center}\section{\'{E}tude de variétés réduites et formule pour
    la multiplicité}\end{center} Si $G$ est un groupe de Lie réductif
connexe, $K\subset G$ est un sous-groupe compact maximal, alors on a la célèbre {\it formule de Blattner}
qui est une identité combinatoire. Un problème intéressant est :
Est-ce que l'on peut exprimer les multiplicités à l'aide de formules
qui ne soient "visiblement " pas combinatoires? La réponse est : Du moins,
si $H\subset G$ est un sous-groupe compact ($G$ est toujours supposé
réductif), oui. En fait dans ce cas, Paradan [[24],[25]) a réussi à
obtenir une telle formule dans le cadre de la géométrie des orbites
coadjointes. Surtout il a réinterprété la formule de Blattner sous
forme non combinatoire. Parmi les objets principaux qu'il a utilisés
figurent les {\it variétés réduites} liées à des variétés
symplectiques d'orbites coadjointes.  Mais si $H\subset G$ n'est pas
un sous-groupe compact, les choses deviennent compliquées. Car il n'y
a pas encore une bonne théorie générale pour les groupes de Lie non
compacts agissant sur les variétés symplectiques. Cependant,
l'objectif de cette section est d'illustrer cette "philosophie" à
travers le groupe $G=SU(2,1)$ et ses sous-groupes non compacts $B_1$
et $B$.\\

Donc dans la suite de cette section, sauf indication contraire, on
garde toutes les notations concernant $G=SU(2,1)$ introduites dans les
chapitres et sections précédents.\\

\begin{subsection}{\'Etude des variétés réduites}
Soit $f_0\in \mathfrak{t}^*\subset \mathfrak{g}^*$ une forme fortement
régulière et $\mathcal{O}_{f_0}:=G.f_0$ sa $G$-orbite coadjointe. Muni
de la structure symplectique de Kirillov-Kostant-Souriau $\varpi$,
$(\mathcal{O}_{f_{0}},\varpi)\cong G/ \mathbb{T}$ est une variété
symplectique, où $\mathbb{T}$ est le tore maximal d'algèbre de Lie
$\mathfrak{t}$ ($\mathfrak{t}$ est défini dans le chapitre 3). $B_1$
(resp. $B$) agit dans $\mathcal{O}_{f_{0}}$ par l'action coadjointe,
de sorte que $(\mathcal{O}_{f_{0}},\varpi,B_1, \text{p}_1)$
(resp. $(\mathcal{O}_{f_{0}},\varpi,B, \text{p})$) devient un espace
hamiltonien, où l'application moment associée est la projection
$\text{p}_1$ (resp. $\text{p}$). On a vu que les $B_1$-orbites
(resp. $B$-orbites) fortement régulières dans $\mathfrak{b}_1^*$
(resp. $\mathfrak{b}^*$) sont $\Omega^{\pm}$
(resp. $\Omega_{r,\pm}$). Dans cette section, on va étudier la variété
réduite $\text{p}_1^{-1}(\Omega^{\pm})/B_1$ que l'on note $X_{\pm}$
(resp. $\text{p}^{-1}(\Omega_{r,\pm})/B$ que l'on note
$X_{r,\pm}$).

\begin{prop} Soit $f_0\in \mathfrak{t}^*\subset \mathfrak{g}^*$ une forme fortement régulière, et $\mathcal{O}_{f_0}:=G.f_0$. Alors\\

 \noindent (1) Soit $\Omega^{\pm}\subset \text{p}_1(\mathcal{O}_{f_0})$. Alors: Si $f_0$ est dans le cône holomorphe, la variété réduite $X_{\pm}$ est difféomorphe à la sphère compacte de dimension 2. Sinon, $X_{\pm}$ est une variété non compacte de dimension 2.\\

\noindent (2) Soit $\Omega_{r,\pm}\subset \text{p}(\mathcal{O}_{f_0})$. Alors, la variété réduite  $X_{r,\pm}$ est réduite à un point.

\end{prop}

\noindent \begin {demo}
Remarquons d'abord que, si $b$ est un élément quelconque de
$\Omega^{\pm}$ (resp. $\Omega_{r,\pm}$), on a $X_{\pm}\cong
\text{p}_1^{-1}(b)/B_{1}(b)$ (resp. $X_{r,\pm}\cong
\text{p}^{-1}(b)/B$.

Soit $S=K.f_{0}$ qui est l'orbite co-adjointe de $f_{0}$ dans
$\mathfrak{k}^{*}$, difféomorphe à la sphère de dimension $2$. Tout
d'abord on montre que $p_{1}(S)$ est une sous-variété compacte de
$\mathfrak{b}_{1}^{*}$ et que $p_{1}$ induit un difféomorphisme de $S$
sur $p_{1}(S)$~: cela résulte de ce qui est dit dans la section
5.3 et de la formule
$p_{1}(x_{1}H+x_{2}F+x_{3}V)=x_{1}E_{2}^*+x_{2}E_{1}^*+x_{3}{E'_{1}}^*$
 (remarquons que l'on peut
donner une autre démonstration de ce fait moins calculatoire et qui se
généralise, reposant sur le fait que le noyau de l'application
linéaire $p_{1}$ est, modulo identification de $\mathfrak{g}^{*}$ avec
$\mathfrak{g}$ au moyen de la forme de Killing,
$\mathfrak{m}\oplus\theta(\mathfrak{n})$).

Soit $f_{\pm}\in S$ tel que
$b_{1\pm}:=p_{1}(f_{\pm})\in\Omega^{\pm}$ et $S_{\pm}:=S\cap
p_{1}^{-1}(\Omega^{\pm})$. Alors
$X_{\pm}=p_{1}^{-1}(b_{1\pm})$ (car dans ce cas, le stabilisateur de $b_{1\pm}$ dans $\mathfrak{b}_1^*$ est trivial). L'application $x\mapsto x.b_{1\pm}$
est un difféomorphisme $M$-équivariant de $B_{1}$ sur $\Omega^{\pm}$
dont on désigne par $\sigma_{\pm}$ l'inverse. On pose
$\Sigma_{\pm}=\sigma_{\pm}\circ p_{1}(S_{\pm})$.

Alors

1 $\Sigma_{\pm}$ est une sous-variété $M$-invariante de $B_{1}$ difféomorphe via
$\sigma_{\pm}$ à $S_{\pm}$.

2 $S_{\pm}=S$ si et seulement si $p_{1}(\mathcal{O}_{f_{0}})\subset
\Omega^{\pm}$. Dans ce cas $S_{\mp}=\emptyset$. Sinon, $S_{\pm}$ est
une sous-variété ouverte non-compacte de la sphère $S$, qui est une
calotte sphérique.

3 Soit $\Delta_{\pm}=\{(g,\sigma_{\pm}(g))|g\in S_{\pm}\}$. Alors
l'application $g\mapsto(g,\sigma_{\pm}(g))$
(resp. $(g,\sigma_{\pm}(g))\mapsto\sigma_{\pm}(g)^{-1}.g$) est un
difféomorphisme $M$-équivariant de $S_{\pm}$ sur $\Delta_{\pm}$
(resp. de $\Delta_{\pm}$ sur $X_{\pm}$).

4 Un système de représentants des $M$-orbites dans $S$
(resp. $S_{\pm}$) est $\{f_{\theta}=\exp-\frac{\theta}{2}
V.f_{0}|\theta\in[0,\pi]\}$ (resp. $\{f_{\theta}=\exp-\frac{\theta}{2}
V.f_{0}|\theta\in[0,\pi] \mbox{ et }
\pm(f_{0}(H)\cos\theta+f_{0}(Z))>0\}$).

5 On pose $\epsilon(f_{0})=\frac{f_{0}(H)+f_{0}(Z)}{\vert f_{0}(H)+f_{0}(Z)\vert}$ et
$f_{\pm}=\pm\epsilon(f_{0})f_{0}$. Alors $f_{\pm}\in S$ si et
seulement si
$p_{1}(\mathcal{O}_{f_{0}})\cap\Omega^{\pm}\neq\emptyset$.

6 Le calcul montre que $p(f_{\theta})=
-\frac{1}{6}(f_{0}(Z)-3f_{0}(H)\cos\theta)W^{*}+(f_{0}(H)\cos\theta+f_{0}(Z))E_{2}^{*}
+(f_{0}(H)\sin\theta) E_{1}^{*}$. D'après la section 5.2, il en résulte que
$\langle\sigma_{\pm}(g_{\theta})^{-1}.g_{\theta},W\rangle =
\phi(\cos\theta)$ où $\phi$ est la fonction homographique $\phi(t)=\frac{1}{6}
(2f_{0}(Z)-3\frac{f_{0}(Z)^{2}-f_{0}(H)^{2}}{f_{0}(H)t+f_{0}(Z)})$. Il est alors clair que
les variétés réduites pour $B$ sont des points.

D'où les résultats.

\end{demo}

\noindent \textbf{Remarque}. Il est possible que l'image réciproque $\text{p}^{-1}(rW^*\pm E_2^*)$ est réduite à un point elle-même, par exemple lorsque $rW^*\pm E_2^*$ est dans la même orbite que $\text{p}(f_0)$.\\
\end{subsection}

\begin{subsection}{Multiplicité et variétés réduites}
On a vu que pour toutes les séries discrètes (holomorphes ou pas)
$\pi_{\lambda}$, on a $\pi_{\lambda}\vert_{B}$ est $B$-admissible et
toute sous-représentation irréductible (de $B$) qui intervient est de
multiplicité 1. Pourtant, il n'y a que exactement un tiers des
$B$-orbites fortement régulières et admissibles dans la projection
$\text{p}(\mathcal{O}_{\pi_{\lambda}})$ qui interviennent dans la
décomposition de $\pi_{\lambda}\vert_{B}$. En fait on a le résultat
suivant (on rappelle que les groupes $G$ et $B$ ont le même centre,
$Z_{G}$)~:
\begin{prop}
Soit $\pi_{\lambda}$ une série discrète de $G=SU(2,1)$ avec
$\lambda\in\mathfrak{t}^{*}$ et soit $f=-i\lambda$. Supposons que la
$B$-orbite fortement régulière et admissible
$\Omega_{\frac{m}{3}\pm\frac{1}{2},\pm}$ soit contenue dans
$p(G.f)$. Alors les assertions suivantes sont équivalentes~:

(i) la série discrète $T_{m,\pm}$ de $B$ intervient dans
$\pi_{\lambda}\vert_{B}$,

(ii) les  séries discrètes $\pi_{\lambda}$ et $T_{m,\pm}$ ont même
caractère central,

(iii) $m-\frac{3f(H)+f(Z)}{2}\in 3\Z$
\end{prop}
\begin{demo}
C'est une conséquence facile de
  l'assertion (4) remarque 1 du chapitre 4 et du corollaire 6.2 qui
  décrivent le caractère central des séries discrètes de $G$ et de
  $B$, ainsi que de la proposition 6.4 et du théorème 6.10.
\end{demo}\\

Dans la suite on va donner une interprétation de ce phénomène dans le
cadre de la quantification $\text{Spin}_{c}$. Cette
interprétation s'inspire des travaux de Paradan (théorème 2.16 dans
[26]) pour les groupes de Lie compacts agissant dans les variétés
symplectiques (non compactes) avec l'application moment propre. Ce qui
est nouveau dans notre cas est que le groupe $B$ n'est pas compact. De
plus pour toute série discrète (holomorphe ou pas), le résultat
s'applique, même si l'application moment (la projection) n'est pas
forcément propre (e.g. pour les séries discrètes ni holomorphes ni
anti-holomorphes).\\

Avant d'aller plus loin, nous devons expliquer quel est le fibré
associé à une orbite de type compact fortement régulière et admissible
dans le cadre de la quantification $\text{Spin}_{c}$. Soit $G$ un
groupe presque algébrique d'algèbre de Lie $\mathfrak{g}$ et
$f\in\mathfrak{g}^{*}$ une forme de type compact fortement régulière
et admissible. Nous supposons que $G$ et $G(f)$ sont connexes. Comme
$\text{Ad}\,G(f)$ est un tore compact, il existe un lagrangien complexe
$l$ et positif dans $\mathfrak{g}_{\C}$ qui soit $G(f)$-invariant. Par
suite, d'après le lemme 2.1, il existe un caractère $\chi_{f}$ de
$G(f)$ de différentielle $\rho_{l}+if\vert_{\mathfrak{g}(f)}$. Il
s'avère d'une part que ce caractère ne dépend pas du choix du
lagrangien positif $l$ et d'autre part que si $\chi_{l}$ désigne le
caractère de $G(g)^{\mathfrak{g}}$ de différentielle $\rho_{l}$, on a
$X(f)=\{\chi_{l}^{-1}\chi_{f}\}$. On définit alors le fibré associé à
l'orbite comme le fibré en droites
$\mathcal{L}=G\times_{G(f)}\C_{\chi_{f}}$.\\

Donnons deux exemples qui nous intéressent.  Soit $\pi_{\lambda}$ une
série discrète de $G=SU(2,1)$ avec $\lambda\in\mathfrak{t}^{*}$ et
soit $f=-i\lambda$ de sorte que
$\mathcal{O}_{\pi_{\lambda}}=G.f$. Dans ce cas, le fibré $\mathcal{L}$
de base $\mathcal{O}_{\pi_{\lambda}}$ est défini par le caractère
$\chi_{f}$ de $\mathbb{T}=G(f)$ de différentielle
le paramètre de Blattner $\Lambda=\lambda+\rho-2\rho_{K}$, où $\rho$
(resp. $\rho_{K}$) est la demi-somme de l'ensemble des racines
positives (resp. positives compactes) défini par $\lambda$.

Considérons la série discrète $T_{m,\pm}$ de $B$, associée à l'orbite
$\Omega_{\frac{m}{3}\pm\frac{1}{2},\pm}$ et soit $h_{m,\pm}$ l'élément
de cette orbite tel que $h_{m,\pm}=(\frac{m}{3}\pm\frac{1}{2})W^{*}\pm
E_{2}^{*}$ et dont le stabilisateur est $M$. En utilisant la
polarisation positive $\mathfrak{l}_{\pm}$ de la section 6.1, on voit
que dans ce cas le fibré $\mathcal{L}$ de base
$\Omega_{\frac{m}{3}\pm\frac{1}{2},\pm}$ est défini par le caractère
$\chi_{h_{m,\pm}}$ dont la différentielle vérifie
$d\chi_{h_{m,\pm}}(W)=i(\frac{m}{3}\pm1)$; par suite
$\chi_{h_{m,\pm}}$ à même restriction à $Z_{G}$ que le caractère
$\sigma_{m}$ de la section 6.1.

\begin{prop}
Soit $\pi_{\lambda}$ une série discrète de $G=SU(2,1)$ avec
$\lambda\in\mathfrak{t}^{*}$ et soit $f=i\lambda$. Supposons que la
$B$ orbite fortement régulière et admissible
$\Omega_{\frac{m}{3}\pm\frac{1}{2},\pm}$ soit contenue dans $p(G.f)$.
Alors, la multiplicité de la série discrète $T_{m,\pm}$ de $B$ dans
$\pi_{\lambda}\vert_{B}$ vaut $1$ si les caractères $\chi_{\lambda}$
et $\chi_{h_{m,\pm}}$ ont même restriction au sous-groupe $M(f)$ et
$0$ sinon.
\end{prop}
\noindent \begin {demo} On vérifie tout d'abord que $M(f)$ est le
  centre $Z_{G}$ du groupe $G$. D'autre part, il résulte de
  l'assertion (4) remarque 1 du chapitre 4 et du corollaire 6.2 que
  le caractère central de $\pi_{\lambda}$ (resp. $T_{m,\pm}$) est
  donné par $\chi_{\lambda}\vert_{Z_{G}}$
  (resp. $\chi_{h_{m,\pm}}\vert_{Z_{G}}$).
\end {demo}\\

Ce résultat est l'analogue du résultat démontré par Paradan dans ([26],
section 2.4 "Quantization of points", Theorem 2.16), dans lequel il
considère un groupe de Lie compact muni d'une action hamiltonienne sur
une variété symplectique avec application moment propre et calcule la
quantification de la variété réduite lorsqu'elle est réduite à un
point. En effet, d'après la proposition 8.1, la variété réduite d'une
orbite $\Omega_{\frac{m}{3}\pm\frac{1}{2},\pm}$ contenue dans $p(G.f)$
est réduite à un point.\\

Par analogie avec le résultat de Paradan, nous dirons que la
quantification $\text{Spin}_{c}$ $Q_{spin}(X_{\frac{m}{3}\pm
  \frac{1}{2},\pm})$ de la variété réduite $X_{\frac{m}{3}\pm
  \frac{1}{2},\pm}$ vaut $1$ si les caractères $\chi_{\lambda}$
et $\chi_{h_{m,\pm}}$ ont même restriction au sous-groupe $M(f)$ et
$0$ sinon (ce cas englobant celui où la variété réduite est vide,
i.e. l'orbite $\Omega_{\frac{m}{3}\pm\frac{1}{2},\pm}$ n'est pas
contenue dans $p(G.f)$).\\

Avec cette convention, on voit que l'assertion (iii) de la conjecture
de Duflo est confirmée pour $(G, \pi_{\lambda}, B)$ pour toutes les
séries discrètes. Puisque nous avons déjà confirmé les assertions (i)
et (ii) de la conjecture de Duflo pour $(G, \pi_{\lambda}, B)$ dans
le chapitre 6, \textbf{la conjecture de Duflo est confirmée pour
  $(G, \pi_{\lambda}, B)$ pour toutes les séries discrètes}.\\

\noindent \textbf{Remarque}. Le fait que toute variété réduite
$\text{p}^{-1}(h_{m,\pm})/M$ soit un point semble donner une explication
géométrique au fait que la multiplicité d'une série discrète de $B$
qui intervient dans $\pi_{\lambda}$ est au plus 1.
\\

D'après le Théorème 6.3 (de Rossi-Vergne) et le théorème 6.10, on sait
que $\pi_{\lambda}\vert_{B_1}$ est $B_1$-admissible si et seulement si
$\pi_{\lambda}$ est holomorphe (ou antiholomorphe). Supposons que
$\pi_{\lambda}$ est holomorphe, alors d'après la section 6.2, on a
$\text{p}_1(\mathcal{O}_{\pi_{\lambda}})=\Omega^-$.  Donc selon la
proposition 8.1, la variété réduite $\text{p}_1^{-1}(\Omega^{-})/B_1$
est une sous-variété de $\mathcal{O}_{\pi_{\lambda}}$ difféomorphe à
la sphère compacte de dimension 2. La théorie nous dit que c'est une
sous-variété symplectique de $\mathcal{O}_{\pi_{\lambda}}$ : on note
$\beta$ sa structure symplectique, qui provient par restriction de
celle (Kirillov-Kostant-Souriau) de
$\mathcal{O}_{\pi_{\lambda}}$. Puisque
$\text{p}_1^{-1}(\Omega^{-})/B_1$ est de dimension deux, alors la
forme volume de Liouville pour $\beta$ est $\frac{\beta}{2\pi}$.

\begin{prop}
Soit $\pi_{\lambda}$ une série discrète holomorphe de $G$, avec le
paramètre de Harish-Chandra $\lambda$. Soit $f_{0}=-i\lambda\in
\mathfrak{t}^*\subset\mathfrak{g}^*$ la forme linéaire associée (au
sens de Duflo) de sorte que $\mathcal{O}_{\pi_{\lambda}}=G.f_0$. Alors

$$\pi_{\lambda}\vert_{B_1}=(\int_{X_{-}}\frac{\beta}{2\pi})\mathrm{T}_{-}.$$

\end{prop}

\noindent \begin {demo} D'après le théorème 6.3, il s'agit de montrer
  que $$\int_{X_{-}}\frac{\beta}{2\pi}=f_0(H),$$ où $H\in\mathfrak{g}$
  est le même que celui dans le théorème 6.3.  Dans la suite on va
  calculer $\int_{X_{-}}\frac{d\beta}{2\pi}$ par deux méthodes
  différentes: La première méthode procède par calcul direct, la
  deuxième plus directe, utilise un argument de topologie
  différentielle.\\

\noindent Méthode 1: Dans la démonstration de la proposition 8.1, on
obtient un difféomorphisme de $K.f_0$ sur
$\text{p}_1^{-1}(\text{p}_1(f_0))$, on le note $\Phi$. Or
l'application $\sigma: (x_2,x_1,x_1')\mapsto
f_0(Z)Z^*+f_0(H)(x_2H^*-x_1F^*-x_1'V^*)$ induit un difféomorphisme de
$S^2=\{(x_2,x_1,x_1'): x_2^2+x_1^2+x_1'^2=1\}$ sur $K.f_0$. Par suite
$\Phi\circ\sigma$ est un difféomorphisme de $S^2$ sur
$\text{p}_1^{-1}(\text{p}_1(f_0))$. Soit ${S^2}'=\{(x_2,x_1,x_1')\in
S^2: x_1'\neq0 \ \text{ou}\ x_1>0 \}$. Donc $\Phi\circ\sigma({S^2}')$
est un ouvert de $\text{p}_1^{-1}(\text{p}_1(f_0))$ de complémentaire
négligeable pour toute densité sur $\text{p}_1^{-1}(\text{p}_1(f_0))$.

Soit $x_{\theta}=-\frac{f_0(H)\sin\theta}{f_0(H)\cos\theta+f_0(Z)}$ et
soit $t_{\theta}$ tel que
$e^{2t_{\theta}}=\frac{f_0(H)\cos\theta+f_0(Z)}{f_0(H)+f_0(Z)}$, où
$Z\in\mathfrak{g}$ est ce que l'on a défini le chapitre 3. Pour
$(\theta,\varphi)\in [0,\pi]\times[0,2\pi]$, on pose
$x(\theta,\varphi)=\exp\varphi W\exp t_{\theta}S\exp
x_{\theta}E_1'\exp-\frac{\theta}{2}V\in G$ et
$g(\theta,\varphi)=x(\theta,\varphi).f_{0}$. Alors
$\{g(\theta,\varphi):(\theta,\varphi)\in [0,\pi]\times
[0,2\pi]\}=\text{p}_1^{-1}(\text{p}_1(f_0))$, de plus l'application
$(\theta,\varphi)\mapsto g(\theta,\varphi)$ est un difféomorphisme de
$[0,\pi]\times [0,2\pi]$ sur $\Phi\circ\sigma({S^2}')$. En fait, $M$
stabilise $\text{p}(f_0)$ (dans $\mathfrak{b}^*$), donc agit dans
$\text{p}^{-1}(\text{p}(f_0))$, donc aussi agit dans
$\text{p}_1^{-1}(\text{p}_1(f_0))$ (car $\mathfrak{b}=\R
W\oplus\mathfrak{b}_1$ et $M=\exp\R W$). Soit $x\in B_1$ et $k\in K$,
tel que $\text{p}_1(x.k.f_0)=\text{p}_1(f_0)$. Alors si $m\in M$ on a
$m(x.k.f_0)=mxm^{-1}(m.k.f_0)$. On peut donc supposer que $k.f$ fait
partie d'un ensemble de représentants des $M$-orbites dans $K.f_0$.,
i.e. on peut supposer que $k.f_0=f_\theta=f_0(Z)Z^*+f_0(H)(\cos\theta
H^*-\sin\theta F^*)=\exp(-\frac{\theta}{2})V.f_0$ avec
$\theta\in[0,\pi]$, ici $F, V\in \mathfrak{k}$ sont définis dans la section 5.3. Dans ce cas, on peut vérifier directement (par les
informations d'après la proposition 5.2 de la section 5.2) que
$\text{p}_1(x.f_\theta)=\text{p}_1(f_0)$ si et seulement si $x=\exp
t_{\theta}S.\exp \exp x_{\theta}E_1'$.

Soit $\beta$ la forme symplectique sur $\Phi\circ\sigma({S^2}')$. Nous
allons calculer $g^*(\beta_{\theta})$. Tout d'abord, on a
$$\frac{\partial g}{\partial \theta}=x(\theta,\varphi)[-\frac{1}{2}V+\frac{dx_{\theta}}{d\theta}\text{Ad}(\exp\frac{\theta}{2}V).E_1'+\frac{dt_{\theta}}{d\theta}\text{Ad}(\exp\frac{\theta}{2}V\exp-x_{\theta}E_1').S].f\ \,$$

$$\frac{\partial g}{\partial \varphi}=x(\theta,\varphi)\text{Ad}(\exp\frac{\theta}{2}V\exp-x_{\theta}E_1').W.f\ .$$

Posons $$U(\theta,\varphi)=-\frac{1}{2}V+\frac{dx_{\theta}}{d\theta}\text{Ad}(\exp\frac{\theta}{2}V).E_1'+\frac{dt_{\theta}}{d\theta}\text{Ad}(\exp\frac{\theta}{2}V\exp-x_{\theta}E_1').S,$$

$$Z(\theta,\varphi)=\text{Ad}(\exp\frac{\theta}{2}V\exp-x_{\theta}E_1').W.$$

Alors on a

$g^*(\beta)=f_0([Z(\theta,\varphi),U(\theta,\varphi)])d\varphi\wedge d\theta$.

Or par des calculs directs, on peut obtenir que $$f_0([Z(\theta,\varphi),U(\theta,\varphi)])=\frac{1}{2}f_0(H)\sin\theta\frac{f_0(H)^2-f_0(Z)^2}{(f_0(Z)+f_0(H)\cos\theta)^2}.$$

On a donc $$g^*(\beta)=\frac{1}{2}f_0(H)\sin\theta\frac{f_0(Z)^2-f_0(H)^2}{(f_0(Z)+f_0(H)\cos\theta)^2}d\theta\wedge d\varphi.$$

Il vient donc $$\int_{X_{-}}\frac{\beta}{2\pi}=\frac{1}{4\pi}\mid\int_{]0,2\pi[\times]0,\pi[} f_0(H)\sin\theta\frac{f_0(Z)^2-f_0(H)^2}{(f_0(Z)+f_0(H)\cos\theta)^2} d\varphi d\theta\mid $$

$$=\frac{1}{4\pi}\mid\int_0^{2\pi}d\varphi.\int_0^{\pi}(f_0(Z)^2-f_0(H)^2)f_0(H)\frac{\sin\theta}{(f_0(Z)+f_0(H)\cos\theta)^2}d\theta\mid $$

$$=\frac{1}{2}(f_0(Z)^2-f_0(H)^2)\int_{-1}^{1}\frac{du}{(f_0(Z)+f_0(H)u)^2}$$
$$=\frac{1}{2}(f_0(Z)^2-f_0(H)^2)[\frac{1}{\mid f_0(Z)\mid-\mid f_0(H)\mid }-\frac{1}{\mid f_0(Z)\mid+\mid f_0(H)\mid }]$$
$$=\mid f_0(H)\mid =f_0(H),$$ d'où le résultat cherché.\\

\noindent Méthode 2: Notons $\Gamma:
G\longrightarrow\mathcal{O}_{\pi_{\lambda}}$ le fibré principal (avec
la fibre $G(f_0)$). Ici, il existe un unique caractère
$\tilde{\chi}_{f_{0}}$ de $G(f_{0})$ de différentille
$if_{0}\vert_{\mathfrak{g}(f_{0})}$. Soit
$\mathcal{L}=G\times_{G(f_{0})}\C_{\tilde{\chi}_{f_{0}}}$ le fibré en droite
de Kirillov-Kostant-Souriau correspondant. Pour simplifier on note
$X_{f_0}:=\text{p}_1^{-1}(\text{p}_1(f_0))$. Il est clair que le fibré
$\mathcal{L}_{X_{f_0}}$ restriction de $\mathcal{L}$ à $X_{f_0}$ est
$\Gamma^{-1}(X_{f_0})\times_{G(f_{0})}\C_{\tilde{\chi}_{f_{0}}}$. Donc via le
difféomorphisme $\Phi: K.f_0\longrightarrow X_{f_0}$,
$\mathcal{L}_{X_{f_0}}$ induit un fibré
$\Phi^{*}(\mathcal{L}_{X_{f_0}})$ de base $K.f_0\subset
\mathcal{O}_{\pi_{\lambda}}$. On peut vérifier directement que
$\Phi^{*}(\mathcal{L}_{X_{f_0}})$ est exactement la restriction de
$\mathcal{L}$ à $K.f_{0}$,
$\mathcal{L}_{K.f_0}:=\Gamma^{-1}(K.f_0)\times_{G(f_{0})}\C_{\tilde{\chi}_{f_{0}}}$.

Soit $\nabla$ la connexion de Kirillov-Kostant-Souriau pour
$\mathcal{L}$. Soient $\nabla_{K.f_0}$ et $\nabla_{X_{f_0}}$ les
connexions induites pour $\mathcal{L}_{K.f_0}$ et
$\mathcal{L}_{X_{f_0}}$ respectivement. Les formes courbures
correspondantes sont $\beta_{\mid K.f_0}$ et $\beta_{\mid
  X_{f_0}}$. D'autre part, $\Phi^{*}(\nabla_{X_{f_{0}}})$ est une
connexion sur le fibré
$\Phi^{*}(\mathcal{L}_{X_{f_0}})=\mathcal{L}_{K.f_0}$ dont la forme
courbure est $\Phi^{*}(\beta_{X_{f_0}})$. Puisque les connexions
$\Phi^{*}(\nabla_{X_{f_0}})$ et $\nabla_{K.f_0}$ proviennent du même
fibré, il est connu que les formes courbures correspondantes
$\Phi^{*}(\beta_{X_{f_0}})$ et $\beta_{K.f_0}$ sont dans la même
classe de cohomologie de de Rham (pour ceci, on peut consulter par
exemple le livre classique de Kobayashi et Nomizu intitulé
"Foundations of Differential Geometry ") . Donc on déduit
que $$\int_{X_{-}}\frac{\beta}{2\pi}=\int_{K.f_0}\frac{\Phi^{*}(\beta_{X_{f_0}})}{2\pi}=\int_{K.f_0}\frac{\beta_{K.f_0}}{2\pi}.$$

La mesure de densité définie par la forme volume apparaissant dans la
dernière intégrale est proportionnelle à la mesure invariante
canonique sur la sphère (ici on considère naturellement $K.f_0$ comme
une sphère de rayon 1). Pour calculer la constante de proportionalité,
il suffit de calculer la valeur de la deux-forme sur deux vecteurs
tangents indépendants en un point. On peut obtenir directement que la
constante concernée est $\frac{f_0(H)}{2}$. Donc
$\int_{K.f_0}\frac{\beta_{K.f_0}}{2\pi}=f_0(H)$,
d'où le résultat.
\end{demo}

Donc dans ce cas l'assertion (iii) de la conjecture de Duflo est
confirmée pour $(G, \pi_{\lambda}, B_1)$, et on peut confirmer de la
même manière l'assertion (iii) de la conjecture de Duflo pour
$(G, \pi_{\lambda}, B_1)$ pour $\pi_{\lambda}$ anti-holomorphe
. Donc en combinant les résultats que l'on a obtenus dans le chapitre
5, \textbf{la conjecture de Duflo est confirmée pour
  $(G, \pi_{\lambda}, B_1)$ pour toutes les séries discrètes}.
\end{subsection}
\pagebreak
\begin{center}\section {Appendice}\end{center}

\begin{subsection}{Sur la projection des orbites coadjointes fortement
  régulières d'une algèbre de Lie simple.}

Le résultat suivant et sa démonstration nous ont été communiqués par
Michel Duflo.

Rappelons que, si $\mathfrak{g}$ est une algèbre de Lie semi-simple,
la forme de Killing induit un isomorphisme de $\mathfrak{g}$-modules de
$\mathfrak{g}^{*}$ sur $\mathfrak{g}$. Nous dirons alors qu'une forme
linéaire sur $\mathfrak{g}$ est semi-simple si, sous cet isomorphisme,
elle correspond à un élément semi-simple de l'algèbre de Lie
$\mathfrak{g}$.

\begin{theo}
Soit $\mathfrak{g}$ une algèbre de Lie simple sur $\mathbb{C}$,
$\mathfrak{h}\subset\mathfrak{g}$ une sous-algèbre de Lie propre et
$p:\mathfrak{g}^{*}\rightarrow\mathfrak{h}^{*}$ la projection
naturelle. Soit $\Omega\subset\mathfrak{g}^{*}$ l'orbite sous l'action
naturelle du groupe adjoint de $\mathfrak{g}$ d'un élément \flqq
semi-simple\frqq~ régulier. Alors  $p(\Omega)$ contient un ouvert de
Zariski non vide de $\mathfrak{h}^{*}$.
\end{theo}

\begin{demo}
Pour démontrer le théorème, on s'appuie sur le lemme suivant
\begin{lem}
Soit $\mathfrak{g}$ une algèbre de Lie algébrique, $G$ un groupe
algébrique connexe d'algèbre de Lie $\mathfrak{g}$,
$\mathfrak{h}\subset\mathfrak{g}$ une sous-algèbre de Lie et
$p:\mathfrak{g}^{*}\rightarrow\mathfrak{h}^{*}$ la projection
naturelle. Soit $\Omega$ une $G$-orbite dans $\mathfrak{g}^{*}$. Alors
les assertions suivantes sont équivalentes~:

(i) $p(\Omega)$ contient un ouvert de Zariski non vide de
$\mathfrak{h}^{*}$.

(ii) Il existe $g\in\Omega$ tel que
$\mathfrak{g}(g)\cap\mathfrak{h}=\{0\}$.

(iii) Il existe $g\in\Omega$ tel que son orbite sous l'action du
sous-groupe analytique de $G$ d'algèbre de Lie $\mathfrak{h}$ soit de
dimension égale à celle de $\mathfrak{h}$.
\end{lem}

\begin{demo}
Les assertions (ii) et (iii) sont clairement équivalentes. L'assertion
(i) s'écrit quant à elle~: il existe $g\in\Omega$ tel que
$p(T_{g}(\Omega))=\mathfrak{h}^{*}$, où $T_{g}(\Omega)$ désigne
l'espace tangent en $g$ à $\Omega$. Ceci s'écrit aussi
$\dim\mathfrak{g}.g-\dim(\mathfrak{g}.g\cap\mathfrak{h}^{\perp})
=\dim\mathfrak{h}$, soit encore
$\dim(\mathfrak{g}.g\cap\mathfrak{h}^{\perp})=
\dim\mathfrak{g}-\dim\mathfrak{g}(g)-\dim\mathfrak{h}$. Considérant
l'orthogonal de $\mathfrak{g}.g\cap\mathfrak{h}^{\perp}$, on voit que
cette dernière relation équivaut à $\dim(\mathfrak{g}(g)+\mathfrak{h})
=\dim\mathfrak{g}(g)+\dim\mathfrak{h}$. D'où le lemme.
\end{demo}

Le résultat suivant est conséquence immédiate du lemme.
\begin{coro}
Sous les hypothèses du théorème, les assertions suivantes sont
équivalentes~:

(i) L'image par $p$ de toute orbite coadjointe \flqq semi-simple\frqq~
régulière de $\mathfrak{g}^{*}$ contient un ouvert de
Zariski non vide de $\mathfrak{h}^{*}$.

(ii) Il existe une sous-algèbre de Cartan $\mathfrak{t}$ de
$\mathfrak{g}$ telle que $\mathfrak{t}\cap\mathfrak{h}=\{0\}$.
\end{coro}

\medskip
Démontrer le théorème se ramène donc à démontrer l'assertion (ii) du
corollaire. Il suffit de le faire lorsque $\mathfrak{h}$ est une
sous-algèbre propre maximale de $\mathfrak{g}$, ce que l'on suppose
désormais. Soit $G$ le groupe adjoint de $\mathfrak{g}$. Alors le
normalisateur $H$ de $\mathfrak{h}$ dans $G$ est un sous-groupe
algébrique d'algèbre de Lie $\mathfrak{h}$. Soit $X=G/H$ qui est une
variété algébrique irréductible isomorphe à la classe de $G$-conjugaison de
$\mathfrak{h}$. Soit $\mathfrak{t}$ une sous-algèbre de Cartan de
$\mathfrak{g}$ et $T$ le sous-groupe de Cartan de $G$
correspondant. On fait agir $T$ sur $X$ et on étudie les sous-groupes
d'isotropie de cette action. Si $xH\in X$, le stabilisateur de $xH$
dans $T$ est $T^{xH}=T\cap xHx^{-1}$ et il a pour algèbre de Lie
$\mathfrak{t}(xH)=\mathfrak{t}\cap \text{Ad} x(\mathfrak{h})$. D'après
Richardson ([27]), il existe un ouvert de Zariski non vide $X'$
de $X$ tel que les sous-groupes d'isotropie $T(xH)$, $xH\in X'$,
soient deux à deux $T$-conjugués, donc tous égaux entre eux. Soit
$x_{0}H\in X'$ et $\mathfrak{a}=\mathfrak{t}\cap\text{Ad}
x_{0}(\mathfrak{h}) =\text{Ad} x_{0}(\text{Ad}
x_{0}^{-1}(\mathfrak{t})\cap\mathfrak{h})$. Alors, on a
$\mathfrak{a}\subset\mathfrak{i}$ où $\mathfrak{i}=\cap_{xH\in X'}\text{Ad}
x(\mathfrak{h})$. Or $\mathfrak{i}$ est un idéal propre de
$\mathfrak{g}$~: en effet si $Y\in\mathfrak{i}$, on a $\text{Ad}
x^{-1}Y\in\mathfrak{h}$ pour tout $xH\in X'$ et donc  pour tout $xH\in
X$, par continuité, de sorte que $\mathfrak{i}=\cap_{xH\in X}\text{Ad}
x(\mathfrak{h})$ est bien un idéal de $\mathfrak{g}$, propre car
contenu dans $\mathfrak{h}$. Comme $\mathfrak{g}$ est simple, on a
$\mathfrak{i}=\{0\}$ et donc $\text{Ad}
x_{0}^{-1}(\mathfrak{t})\cap\mathfrak{h}=\{0\}$. D'où le théorème
\end{demo}

\end{subsection}

\begin{subsection}{Sur le comportement asymptotique au voisinage d'un
    point singulier de première espèce des solutions d'un système
    différentiel linéaire ordinaire.}

Les résultats de cet appendice concernant le comportement asymptotique
au voisinage de $0$ des solutions systèmes étudiés dans la section 6.5
nous ont été communiqué par Claude Sabbah.\\

On considère un système différentiel de taille $\ell$
\[\tag*{$(*)$}
y'(z)=\frac{M(z)}{z}\cdot y(z),
\]
où $M(z)$ est holomorphe en $0$. On note $M_0=M(0)$ et $M_0^{(s)}$ la partie semi-simple de~$M_0$.

\begin{prop}
La solution fondamentale $Y(z)$ du système $(*)$ peut être mise sous la forme
\[
Y(z)=U(z)z^{\Delta_0}z^N,
\]
où $\Delta_0$ est diagonale et équivalente à $M_0^{(s)}$, $N$ est strictement triangulaire inférieure, et $U(z)$ est holomorphe inversible (c'est-à-dire $\det U(0)\neq0$).
\end{prop}

\noindent \textbf{Remarque}. Dans cette écriture, il se peut que $N$ ne commute pas à $\Delta_0$ lorsque certaines valeurs propres de $\Delta_0$ diffèrent d'un entier non nul. L'ordre de l'écriture est donc important.\\

\begin{demo}
Si on pose $\widetilde{Y}(z)=P(z)Y(z)$ avec $P$ holomorphe inversible, et si $Y(z)$ est une matrice fondamentale de $(*)$, alors $\widetilde{Y}(z)$ est une solution fondamentale du système analogue de matrice $\widetilde{M}(z)=PMP^{-1}+zP'P^{-1}$.

On effectue un premier changement de base constant $P^0$ de sorte que la nouvelle matrice $P^0M_0(P^0)^{-1}$, que je note encore $M_0$ pour simplifier, soit sous forme de Jordan triangulaire inférieure. On choisit aussi $P^0$ de sorte que les parties entières des éléments diagonaux de la diagonale $\Delta_0$ soient rangés en ordre croissant. Je les note $d_1\leq d_2\leq\cdots\leq d_\ell$. Pour simplifier, je note encore $M(z)$ la matrice obtenue et je note $D$ la matrice diagonale $\text{diag}(d_1,\dots,d_\ell)$. Enfin, je pose $d=d_\ell-d_1\geq0$. La matrice $D$ commute à $\Delta_0$ et $M_0$, et les valeurs propres de $M_0-D$ ont une partie réelle dans $[0,1[$, de sorte que la seule valeur propre entière de $\text{ad}(M_0-D)$ est $0$.

\begin{lem}
Il existe un changement de base holomorphe inversible $P^1(z)$ tel que $\widetilde{ M}=M_0+zB_1+\cdots+z^d B_d$, où les matrices constantes $B_k$ satisfont à $\text{ad}(D)(B_k)=kB_k$ (en particulier sont strictement triangulaires inférieures).
\end{lem}

On admet provisoirement ce lemme, et on écrit $\widetilde{ M}=z^DBz^{-D}$, avec $B=M_0+B_1+\cdots+B_d$. En effectuant le changement de variable (méromorphe) de matrice~$z^{-D}$, on trouve donc que $z^{-D}P^1P^0Y$ est solution fondamentale du système de matrice $B-D$. Maintenant, $B-D$ s'écrit $\Delta_0-D+M_0^{(n)}+B_1+\cdots+B_d$, où $M_0^{(n)}$ est la partie nilpotente de $M_0$, donc strictement triangulaire inférieure et commute à $\Delta_0$. Par conséquent, il existe une matrice constante inversible $P^2$ triangulaire inférieure telle que $P^2(B-\nobreak D)(P^2)^{-1}=\Delta_0-D+N$, avec $N$ strictement triangulaire inférieure commutant à $\Delta_0-D$ (mais peut-être pas à $\Delta_0$). Ainsi, on a, en posant $\widetilde{ Y}=P^2z^{-D}P^1P^0Y$,
$[
\widetilde{ Y}'=\frac{\Delta_0-D+N}z\cdot\widetilde{ Y}
]$
et donc $\widetilde{Y}=z^{\Delta_0-D}z^N$. Par suite,
\begin{align*}
Y&=(P^1P^0)^{-1}z^D(P^2)^{-1}z^{\Delta_0-D}z^N\\
&=[(P^1P^0)^{-1}z^D(P^2)^{-1}z^{-D}]z^{\Delta_0}z^N,
\end{align*}
et on pose $U(z)=(P^1P^0)^{-1}z^D(P^2)^{-1}z^{-D}$. La matrice $T:=(P^2)^{-1}$ est triangulaire inférieure. Par suite les entrées $z^{d_i-d_j}t_{ij}$ de $z^DTz^{-D}$ son nulles si $i<j$, et puisque la suite $(d_i)$ est croissante, la matrice $z^DTz^{-D}$ est holomorphe inversible. Finalement, $U(z)$ est holomorphe inversible, comme voulu.
\end{demo}

\begin{coro}
La dimension de l'espace des solutions de $(*)$ qui sont localement $L^2$ en $z=0$ pour la mesure $dz/z$ sur la demi-droite réelle $\R_+$ est égale au nombre de valeurs propres  de $M(0)$ ayant une partie réelle strictement positive.
\end{coro}

\begin{demo}
Une solution $y$ de $(*)$ est une combinaison linéaire à coefficients constants des vecteurs colonnes de $Y$. Elle est de la forme $U(z)\widetilde{y}$, où $\widetilde{y}$ est la combinaison linéaire correspondante des vecteurs colonnes de $\widetilde{ Y}$. De plus, $y$ est $L^2$ si et seulement si $\widetilde{y}$ l'est, puisque $U(z)$ est inversible. Si $\Delta_0=\text{diag}(\delta_1,\dots,\delta_\ell)$, les colonnes de $\widetilde{ Y}$ sont de la forme
\[
\begin{pmatrix}0\\\vdots\\0\\z^{\delta_j}\\n_{j+1,j}z^{\delta_{j+1}}\log z\\\vdots\\n_{\ell,j}z^{\delta_\ell}\dfrac{(\log z)^{\ell-j}}{(\ell-j)!}\end{pmatrix}
\]
et on voit qu'une combinaison linéaire à coefficients constants non nuls de colonnes est $L^2$ si et seulement si chaque $\text{Re}\delta_j$ correspondant est $>0$.
\end{demo}

\begin{demo}[Démonstration du lemme]
On cherche une matrice $P^1(z)=\text{Id}+zP^1_1+\cdots$ et des matrices $B_i$ comme dans le lemme, de sorte que l'on ait

\begin{equation}\label{eq:ccPB}
z\frac{dP^1}{dz}=(M_0+zB_1+\cdots+z^d B_d)P^1(z)-P^1(z)(M_0+zM_1+z^2M_2+\cdots).
\end{equation}
On commence par considérer les équations pour $j=1,\dots,d$. On écrit
\[
j P^1_j=M_0P^1_j+B_j-P^1_j M_0+\Psi_j(P^1_1,\dots,P^1_{j-1};B_1,\dots,B_{j-1};M_0,\dots,M_j),
\]
où $\Psi_j$ est connu par récurrence sur $j$, et on veut déterminer $P^1_j$ et $B_j$. On écrit ceci sous la forme
\[
(j\text{Id}-\text{ad} M_0)(P^1_j)=B_j+\Psi_j.
\]
L'endomorphisme $j\text{Id}-\text{ad} M_0$ est inversible si $j\geq d+1$: en effet, les valeurs propres de $\text{ad} M_0$ sont les différences des valeurs propres de $M_0$; si une telle différence est égale à l'entier $j$, autrement dit si $j\text{Id}-\text{ad}M_0$ a une valeur propre nulle, la différence des parties entières correspondantes est aussi égale à $j$, et donc $j\leq d$.

Puisque $\text{ad} M_0$ commute à $\text{ad}D$, on peut décomposer cette équation sur les sous-espaces propres de $\text{ad} D$. De plus, puisque la seule valeur propre entière de $\text{ad}(M_0-D)$ est $0$, l'endomorphisme $\text{ad}(M_0-D)+k\text{Id}$ est inversible pour tout entier $k\neq0$. Sa restriction à chaque espace propre de $\text{ad} D$ satisfait à la même propriété.

Ceci étant rappelé, on cherche donc à résoudre, pour tout entier $k$, l'équation
\begin{equation}\label{eq:ccPk}
(j\text{Id}-\text{ad} M_0)(P_j^{1(k)})=B_j^{(k)}+\Psi_j^{(k)},
\end{equation}
où $Q^{(k)}$ désigne la composante de la matrice $Q$ sur l'espace propre de valeur propre $k$ de $\text{ad} D$, qui satisfait donc $\text{ad} D(Q^{(k)})=kQ^{(k)}$.
\begin{enumeratea}
\item\label{cas1}
Si $k\neq j$, on a $B_j^{(k)}=0$ et $(j\text{Id}-\text{ad} M_0)$ coïncide sur cet espace propre avec l'endomorphisme $(j-k)\text{Id}-\text{ad}(M_0-D)$ qui, on l'a vu ci-dessus, est inversible. On peut donc trouver une solution (unique) à l'équation \eqref{eq:ccPk}.
\item\label{cas2}
Si $k=j$, on doit résoudre sur cet espace propre l'équation $\text{ad}(D-M_0)(P_j^{1(j)})-B_j^{(j)}=\Psi_j^{(j)}$ en déterminant $B_j^{(j)}$ par la même occasion.

Choisissons un supplémentaire de l'image de $\text{ad}(D-M_0)$ dans cet espace propre. Alors on peut décomposer $\Psi_j^{(j)}$ en somme d'un élément de l'image de $\text{ad}(D-M_0)$, ce qui donne $P_j^{1(j)}$ (de manière non unique) et d'un élément dans ce supplémentaire, qu'on baptise $B_j^{(j)}$.
\end{enumeratea}

On continue maintenant la récurrence lorsque $j\geq d+1$, tous les $B_i$ étant connus. On procède exactement comme dans le cas où $j\leq d$, mais il n'est plus nécessaire d'introduire des termes correctifs $B_j$ puisque $j\text{Id}-\text{ad} M_0$ est inversible, et on détermine $P^1_j$ comme dans le cas~\eqref{cas1}.

Reste à montrer la convergence de la série $P^1(z)$, ce qui est un résultat classique, puisque le système \eqref{eq:ccPB} est à pôle simple.
\end{demo}
\end{subsection}\\

\pagebreak

 Ancienne adresse: UMR 6086 CNRS, Université de Poitiers,
 Laboratoire de Mathématiques et Applications, Boulevard Marie et
 Pierre Curie, BP~30179,
 86962 Chasseneuil Cedex, France; liu@math.univ-poitiers.fr\\
 
 Adresse actuelle: Leibniz Universit\"{a}t Hannover, Institut f\"{u}r Analysis
, Welfengarten 1, 30167 Hannover Germany; liu@math.uni-hannover.de

\end{document}